\documentclass[11pt,reqno]{amsart}
\usepackage[T1]{fontenc}
\usepackage{lmodern,mathtools,amssymb,microtype}
\usepackage[margin=1in]{geometry}
\usepackage{xurl,xcolor}
\usepackage[colorlinks=true,linkcolor=blue!40!black,citecolor=blue!40!black,urlcolor=blue!40!black]{hyperref}
\newtheorem{theorem}{Theorem}[section]
\newtheorem{proposition}[theorem]{Proposition}
\newtheorem{lemma}[theorem]{Lemma}
\newtheorem{corollary}[theorem]{Corollary}
\theoremstyle{remark}
\newtheorem{remark}[theorem]{Remark}
\newtheorem{question}[theorem]{Question}
\newcommand{\Q}{\mathbb Q}
\newcommand{\Z}{\mathbb Z}
\newcommand{\F}{\mathbb F}
\DeclareMathOperator{\Gal}{Gal}
\DeclareMathOperator{\Li}{Li}
\DeclareMathOperator{\Res}{Res}
\DeclareMathOperator{\lc}{lc}
\DeclareMathOperator{\disc}{disc}
\DeclareMathOperator{\lcm}{lcm}
\DeclareMathOperator{\Sp}{Sp}
\DeclareMathOperator{\GSp}{GSp}
\DeclareMathOperator{\Hdg}{Hdg}
\DeclareMathOperator{\MT}{MT}
\DeclareMathOperator{\End}{End}
\DeclareMathOperator{\CH}{CH}
\newcommand{\doi}[1]{\href{https://doi.org/#1}{\nolinkurl{doi:#1}}}
\numberwithin{equation}{section}
\title[Arithmetic of Jacobi polynomials]{Arithmetic of Jacobi polynomials}
\author{Khai-Hoan Nguyen-Dang}
\thanks{ORCID: \href{https://orcid.org/0000-0002-6765-811X}{0000-0002-6765-811X}}
\date{}
\subjclass[2020]{Primary 11R09, 11R32; Secondary 11N13, 14K15, 20B15, 31B05, 33C45.}
\keywords{Jacobi polynomial, Legendre polynomial, Kaneko--Zagier polynomial, Galois group, resultant, Newton polygon}
\hypersetup{pdftitle={Arithmetic of Jacobi polynomials},pdfauthor={Khai-Hoan Nguyen-Dang}}
\begin{document}
\begin{abstract}
We study the arithmetic of Jacobi polynomials. As applications, we prove the Mahlburg--Ono conjecture on the irreducibility and Galois
groups of the four Kaneko--Zagier families, and the Cullinan--Hajir conjecture on the Galois groups of Legendre polynomials. Our argument also recovers Deng's irreducibility theorem for Legendre polynomials. Other applications concern common zeros of Jacobi polynomials, harmonic functions on cones, and hyperelliptic Jacobians.
\end{abstract}
\maketitle

\section{Introduction}\label{sec:introduction}

For rational parameters $\alpha,\beta>-1$, the zeros of the Jacobi polynomials
are algebraic numbers. Orthogonality places them in $(-1,1)$ and makes them
simple. Much less is known about their arithmetic.

We use the normalization
\begin{equation}\label{eq:jacobi-definition}
 P_n^{(\alpha,\beta)}(x)=2^{-n}\sum_{j=0}^n
 \binom{n+\alpha}{n-j}\binom{n+\beta}{j}
 (x-1)^j(x+1)^{n-j},
 \qquad n\ge0,
\end{equation}
where $\binom zj=z(z-1)\cdots(z-j+1)/j!$ and $\binom z0=1$.
For $\alpha,\beta>-1$, this is the degree-$n$ polynomial with
$P_n^{(\alpha,\beta)}(1)=\binom{n+\alpha}{n}$ that is orthogonal
to every polynomial of degree less than $n$ for the weight
$(1-x)^\alpha(1+x)^\beta$ on $(-1,1)$. The Legendre polynomials
are $P_n=P_n^{(0,0)}$, so $P_n(1)=1$. We also use Jacobi
polynomials outside this orthogonal range in the Galois arguments.
Standard identities are collected in \cite[Chapter~IV]{Szego} and
\cite[Chapter~18]{DLMF}.

For a nonconstant separable polynomial $f\in K[x]$, we write $\Gal(f/K)$ for
the Galois group of its splitting field, acting on its roots. The
polynomial is irreducible precisely when this action is transitive.

The irreducibility problem for Legendre polynomials goes back to Stieltjes. In his letter to Hermite
of 2 October 1890, he asked whether $P_{2m}$ and $P_{2m+1}/x$ are
irreducible over $\Q$ \cite[Letter~275]{Stieltjes}; see also
\cite[\S1]{CHL}. The division by $x$ removes the zero forced by
odd parity. Schur's work on Hermite and Laguerre polynomials showed
how coefficient valuations can determine factorization and Galois
groups \cite{Schur}. For Legendre polynomials, Wahab obtained
substantial infinite families of irreducible examples by Newton
polygons \cite{Wahab}. Related questions for Bessel polynomials
were settled through the work of Filaseta and Filaseta--Trifonov
\cite{Filaseta,FT}.

In a recent preprint, Deng gives a proof of Stieltjes' conjecture in every degree
\cite[Theorem~1.1]{Deng}. His argument compares the real size of a
factor resultant with the size permitted by its p-adic
valuations. This comparison is our starting point. We first give
it a common Jacobi form, then develop the arithmetic estimates
needed for the different parameter families. The Galois arguments
use another local input: ramification at suitable primes produces
long cycles on the roots. The interaction of these two methods is
also what gives the results on higher derivatives and number fields.

\subsection{Galois groups}

Our first examples come from modular forms. Kaneko and Zagier constructed polynomials over $\Q$ whose reductions at suitable
primes describe the supersingular $j$-invariants of elliptic curves
\cite{KZ}. After an affine change of variable and normalization, their Jacobi
factors give the four families
\[
 P_n^{(\lambda/3,\mu/2)},\qquad \lambda,\mu\in\{-1,1\}
\]
\cite[Eq.~(3)]{CHK}. Thus their rational factorization describes
the algebraic relations among the characteristic-zero lifts of
these modular zeros.

\begin{theorem}[The Mahlburg--Ono conjecture]\label{thm:kz-all}
For every $n\ge1$ and $\lambda,\mu\in\{-1,1\}$,
$P_n^{(\lambda/3,\mu/2)}$ is irreducible over $\Q$ and
\[
 \Gal(P_n^{(\lambda/3,\mu/2)}/\Q)=S_n.
\]
\end{theorem}

Mahlburg and Ono conjectured both assertions in 2004 and proved
irreducibility in infinitely many degrees
\cite[p.~40 and Theorem~1.1]{MO}. Cullinan and Hajir proved further
infinite families and established the implication from
irreducibility to $S_n$ through modular index $\ell\le10^9$
\cite[Theorems~1.8--1.9]{CHK}. Cullinan and Gajek-Leonard studied
the associated Newton polygons \cite{CGL}.
Theorem~\ref{thm:kz-all} treats all four families in every positive
degree. In Section~\ref{sec:modular-interpretation} we will give the precise
normalization and explain the relation with modular differential
equations. Appendix~\ref{sec:local-irred} gives stronger local
irreducibility statements in prime-power degrees and in degrees
$5\cdot3^a$.

For Legendre polynomials, the symmetry of the zeros determines the
largest possible group. Write
\[
 L_{2m+\delta}(x)=P_{2m+\delta}(x)/x^\delta,
 \qquad m\ge1,\quad\delta\in\{0,1\}.
\]
The $2m$ roots occur in opposite pairs. Every Galois automorphism
preserves these pairs, so the group is contained in
\[
 W_m=C_2\wr S_m=\F_2^m\rtimes S_m.
\]
Here $S_m$ permutes the pairs and $\F_2^m$ independently exchanges
the two members of each pair. The irreducibility in Theorem~\ref{thm:all-degrees} gives transitivity of the group action. Our next result determines all these permutations.

\begin{theorem}[The Legendre Galois groups]\label{thm:legendre-all}
For every $m\ge1$ and $\delta\in\{0,1\}$,
\[
 \Gal(P_{2m+\delta}(x)/x^\delta\,/\Q)=W_m.
\]
\end{theorem}

This proves both parts of \cite[Conjecture~1.5]{CHL}: the induced
group on the pairs is $S_m$, and the kernel of this action is all
of $\F_2^m$. Cullinan and Hajir proved the first part, conditional
on irreducibility, through reduced degree $10^9$
\cite[Theorem~1.9(a)]{CHL}. Our argument removes that bound and
determines the kernel in every degree.

\subsection{Integral parameters}

Now we consider parameters $(\alpha,\beta)$ of Jacobi polynomials to be integers.

\begin{theorem}\label{thm:main}
Fix distinct nonnegative integers $a,b$. There is an effectively
computable integer $N(a,b)$ such that, for every $n\ge N(a,b)$,
$P_n^{(a,b)}$ is irreducible over $\Q$ and has Galois group $S_n$.
\end{theorem}

Moreover, for every fixed number field $E$, the splitting field of $P_n^{(a,b)}$ is linearly disjoint from $E$ over $\Q$ for all sufficiently large $n$, effectively in $a,b,E$. Consequently the polynomial remains irreducible over $E$ with group $S_n$.

The restriction to sufficiently large degree is necessary. For
example, $P_2^{(4,1)}(x)=3x(3x+2)$. Negative integral parameters give endpoint factors. After removing them, the same conclusion
holds for each fixed integral pair with unequal absolute values
(Corollary~\ref{cor:integral-parameters}). In our statements the
degree varies and the parameters are fixed, while the Jacobi families
studied in \cite{CHS} have a different variation: the degree is
fixed and a rational parameter varies.

Equal parameters connect the problem to derivatives. Define the
Gegenbauer polynomials by
\[
 (1-2xt+t^2)^{-\lambda}
 =\sum_{j\ge0}C_j^\lambda(x)t^j.
\]
For an integer $r\ge0$, the polynomial $C_n^{r+1/2}$ is a nonzero
rational multiple of $P_n^{(r,r)}$ \cite[(18.7.1)]{DLMF}, and
\begin{equation}\label{eq:derivative-gegenbauer}
 \frac{d^r}{dx^r}P_{n+r}(x)
 =(2r-1)!!\,C_n^{r+1/2}(x),
 \qquad (2r-1)!!=\prod_{j=1}^r(2j-1).
\end{equation}
The empty product is $1$. A single superscript $P_\ell^{(r)}$
will denote a derivative; two superscripts denote Jacobi parameters.
After removing the zero forced by parity, put
\begin{equation}\label{eq:reduced-gegenbauer}
 F_{m,r,\delta}(x)=\frac{C_{2m+\delta}^{r+1/2}(x)}{x^\delta},
 \qquad m\ge1,\quad\delta\in\{0,1\}.
\end{equation}
This polynomial has degree $2m$ in $x$. All irreducibility
assertions for $F_{m,r,\delta}$ refer to that full degree.

\begin{theorem}\label{thm:symmetric}
For integers $r\ge0$, $m\ge1$, and $\delta\in\{0,1\}$,
the polynomial $F_{m,r,\delta}$ is irreducible over $\Q$ whenever
\[
 m\ge\max\{100r,1000\}.
\]
There is also an effective absolute integer $M_0$ such that the same
conclusion holds whenever $m\ge\max\{8r,M_0\}$.
\end{theorem}

The conclusion is uniform in the parameter: the derivative order
may grow with $m$. In particular, the second assertion covers
every $0\le r\le m/8$ once $m$ exceeds one absolute threshold.
For derivative orders $r=0$ and $r=1$, the conclusion holds
without any restriction on $m$.

\begin{theorem}\label{thm:all-degrees}
For every $m\ge1$, the polynomials
\[
 P_{2m}(x),\qquad P_{2m+1}(x)/x,\qquad
 P'_{2m+1}(x),\qquad P'_{2m+2}(x)/x
\]
are irreducible over $\Q$.
\end{theorem}

Deng proved the first two assertions in~\cite[Theorem~1.1]{Deng}. We recover them here from the same comparison that treats the derivatives. The last two say that the nonzero critical points of each Legendre polynomial $P_\ell$,
$\ell\ge3$, form a single rational Galois orbit. Rolle's theorem
places these critical points between consecutive zeros; the
arithmetic assertion makes them conjugate over $\Q$. It also gives
a common-zero theorem: distinct members of
$\{P_\ell,P_\ell':\ell\ge1\}$ have no common nonzero root
(Corollary~\ref{cor:common-roots}). This extends Deng's common-zero
result for the Legendre polynomials themselves \cite{DengRoots}.

\subsection{Rational parameters}

We consider parameters $(\alpha,\beta)$ to be rational with certain conditions. One can observe that local ramification determines
the Galois group on almost all roots even when irreducibility of the whole polynomial is not known. The condition below concerns the reduced denominators of the parameters and their sum. It
ensures that the relevant Jacobi coefficients have one long Newton edge at suitably chosen primes.

\begin{theorem}\label{thm:jacobi}
Fix $\alpha=a/u$ and $\beta=b/v$ in lowest terms, with $u,v>0$,
and let $D$ be the reduced denominator of $\alpha+\beta$. Assume
\[
 D\text{ is even},\qquad D\ge\max\{u,v\},\qquad
 \alpha,\beta\notin\{-1,-2,\ldots\}.
\]
For all sufficiently large $n$, $f_n=P_n^{(\alpha,\beta)}$ has a
unique monic irreducible factor $h_n$ of degree $M_n>n/2$, with
\[
 \Gal(h_n/\Q)=S_{M_n},\qquad
 n-M_n=O_{\alpha,\beta,A}\!\left(\frac{n}{(\log n)^A}\right)
 \quad\text{for every fixed }A>0.
\]
Write $f_n=\kappa_nh_ng_n$, with $g_n$ monic and
$\kappa_n\in\Q^\times$. The splitting fields of $h_n$ and $g_n$
are linearly disjoint over $\Q$, and restriction gives
\begin{equation}\label{eq:general-product}
 \Gal(f_n/\Q)=S_{M_n}\times\Gal(g_n/\Q).
\end{equation}
For every fixed $c\in\Q^\times$ and all sufficiently large $n$,
$h_n(1+cx^2)$ is irreducible with group $W_{M_n}$. Its splitting
field is linearly disjoint from that of $g_n(1+cx^2)$, and
\begin{equation}\label{eq:general-quadratic-product}
 \Gal(f_n(1+cx^2)/\Q)=W_{M_n}\times\Gal(g_n(1+cx^2)/\Q).
\end{equation}
The splitting fields of $h_n$ and $h_n(1+cx^2)$ are also each
linearly disjoint from every fixed number field once $n$ is
sufficiently large.
\end{theorem}

The threshold may depend on the fixed parameters, on $c$, and on
the number field; in the degree estimate it may also depend on
$A$. A constant polynomial has splitting field $\Q$ and trivial
Galois group. The theorem separates a large symmetric part from
a small complementary part and proves that their splitting fields
have no common extension of $\Q$.

No orthogonality hypothesis is needed here. Some restriction on
the parameters is necessary: the specialization
$\alpha=\beta=-1/2$ gives Chebyshev polynomials
\cite[(18.7.3)]{DLMF}, whose zeros lie in cyclotomic fields and
whose Galois groups are abelian. The small complementary factor
$g_n$ in the theorem is not known to be constant so eventual
irreducibility under the stated denominator condition remains a
natural question.

The integral irreducibility theorem and the rational-parameter
Galois theorem meet through a classical quadratic transformation.
For $\delta=0,1$, the polynomial $F_{m,r,\delta}(x)$ is a nonzero
rational multiple of
\[
 P_m^{(\delta-1/2,r)}(1-2x^2).
\]
The parameters $\delta-1/2$ and $r$ satisfy the denominator
condition with $D=2$. Write
$F_{m,r,\delta}(x)=Q_{m,r,\delta}(x^2)$, where
$Q_{m,r,\delta}\in\Q[t]$ has degree $m$.

\begin{corollary}
\label{cor:derivative-galois}
Fix an integer $r\ge0$ and a number field $E$. For all sufficiently
large $m$, and for both $\delta\in\{0,1\}$, the polynomials
$F_{m,r,\delta}$ and $Q_{m,r,\delta}$ are irreducible over $E$, with
\[
 \Gal(F_{m,r,\delta}/E)=W_m,\qquad
 \Gal(Q_{m,r,\delta}/E)=S_m.
\]
Their splitting fields over $\Q$ are each linearly disjoint from
$E$ over $\Q$.
\end{corollary}

The two groups describe the two natural sets of roots: $S_m$ acts
on their squares, while $W_m$ also records the independent changes
of sign. Here $r$ and $E$ are fixed. For $r=0$ and $E=\Q$,
Theorem~\ref{thm:legendre-all} gives every degree. For $r=1$,
irreducibility holds in every degree, whereas the full
wreath-product conclusion is proved for sufficiently large degree.

\subsection{The resultant method}
The idea is to estimate the resultant in two ways: the real geometry of the roots gives a lower
bound, while their local arithmetic gives an upper bound.
Section~\ref{sec:comparison} develops the common Jacobi form of
this argument. It extends Deng's comparison of real and
$p$-adic factor resultants for Legendre polynomials
\cite[\S1.2]{Deng}. 

Fix rational parameters $\alpha,\beta>-1$ and put
$y_n(t)=P_n^{(\alpha,\beta)}(1-2t)$. This change of variable
preserves irreducibility over $\Q$ and places the roots in $(0,1)$.
Suppose that a primitive integral multiple $H=\kappa y_n$ admits
a factorization
\[
 H=AB,\qquad A,B\in\Z[t]\ \text{primitive},\qquad
 1\le d=\deg A\le e=\deg B.
\]
The roots are simple by orthogonality
\cite[\S18.2(vi)]{DLMF}, so the resultant $R=\Res(A,B)$ is a nonzero integer.
Consequently,
\begin{equation}\label{eq:intro-product-formula}
 \log|R|=\sum_p v_p(R)\log p,
\end{equation}
where $v_p(p)=1$. We seek incompatible bounds for the two sides,
uniformly in the possible factor degree $d$.

The common input is an auxiliary polynomial constructed from
orthogonality. Let $I$ denote integration against the Jacobi
weight $t^\alpha(1-t)^\beta\,dt$ on $(0,1)$, normalized so that
$I(1)=1$. Its moments are
\[
 \mu_j=I(t^j)
 =\frac{(\alpha+1)_j}{(\alpha+\beta+2)_j}\in\Q,
\]
where $(u)_0=1$ and $(u)_j=u(u+1)\cdots(u+j-1)$ for $j\ge1$;
see \eqref{eq:jacobi-moments}. Define
\[
 Q_H(t)=I_s\!\left(\frac{H(t)-H(s)}{t-s}\right),
 \qquad \sigma(t)=t(1-t).
\]
The divided difference is a polynomial, so $Q_H\in\Q[t]$.
Orthogonality and the Jacobi differential equation give the
identity of Lemma~\ref{lem:jacobi-wronskian}, which at a root $z$
of $H$ becomes
\[
 \sigma(z)H'(z)Q_H(z)=C_H,
 \qquad C_H=\kappa^2C_n\in\Q_{>0}.
\]
Here $C_n$ is the explicit constant in
\eqref{eq:jacobi-wronskian}. Because the moments and the constant
are rational, this identity can be used both over $\mathbb R$
and over $\Q_p$.

At the real place, the Gaussian quadrature weight $\omega_z$
attached to a root $z$ satisfies
$Q_H(z)=H'(z)\omega_z$. These weights are positive and sum to one
\cite[\S3.5(v), equations (3.5.18) and (3.5.20.1)]{DLMF}. Since $0<z(1-z)\le1/4$, the root identity gives
$|H'(z)|\ge2\sqrt{C_H}$; this is the derivative estimate in
Lemma~\ref{lem:jacobi-wronskian}. For integral parameters, energy
estimates from the differential equation give sharper bounds,
as shown in the proofs of
Theorems~\ref{thm:criterion} and~\ref{thm:symmetric-exact}.

To pass from derivatives to the resultant, write $a=\lc A$ and
let $\xi_1,\ldots,\xi_d$ be the roots of $A$. Differentiating
$H=AB$ at these roots gives
\begin{equation}\label{eq:intro-resultant-identity}
 |R|=
 \frac{|a|^{e-d}}
      {\prod_{1\le i<j\le d}|\xi_i-\xi_j|^2}
 \prod_{i=1}^d|H'(\xi_i)|.
\end{equation}
See also \cite[Lemma~5.3]{Deng}. The denominator is the square
of a Vandermonde determinant. The Chebyshev--Hadamard estimate
of \cite[Lemma~5.5]{Deng}, recalled in
Lemma~\ref{lem:interval-vandermonde}, bounds it from above.
Together with the derivative estimate and the information
supplied by the endpoint values of $A$, this gives the real
lower bounds in Lemma~\ref{lem:real-resultant} and
Proposition~\ref{prop:arch-resultant}.

At a finite prime, the same root identity controls how closely
the roots can lie. The moment denominators bound the
coefficients of $Q_H$, and hence bound $v_p(H'(z))$ when $z$
is integral. For a nonintegral root, we pass to the reciprocal
polynomial. Lemma~\ref{lem:projective-moment-budget} combines
these two cases by using projective root distances.
Their nonnegativity allows us to sum over the roots of the
smaller factor, giving
\[
 v_p(R)\le d\,\mathcal E_p,
 \qquad
 \mathcal E_p=
 v_p(C_H)-\min_{0\le j\le2n-1}v_p(\mu_j).
\]
The exact expression of $v_p(R)$ as a sum of projective
distances is \eqref{eq:projective-resultant}. This formulation
also accounts for primes dividing the leading coefficient,
at which some roots may have negative valuation.

A further step restricts the primes that can contribute.
Reduction of the differential equation excludes multiple
reduced roots away from certain distinguished points.
Bounds for the endpoint, central, and leading coefficients
then prevent the two factors from meeting at the remaining
points of the projective line. The resulting exclusions are
proved in \eqref{eq:asymmetric-prime-support},
\eqref{eq:symmetric-prime-support}, and
Lemma~\ref{lem:good-prime-resultant}. They make the local sum
small enough for the comparison to succeed in the ranges
stated in our theorems.

The estimates needed to complete the argument depend on the
parameters. For fixed unequal nonnegative integers $a,b$,
Lemma~\ref{lem:center} expresses
$P_{2m+\delta}^{(a,b)}(0)$, with $\delta\in\{0,1\}$, as a
central binomial factor times a rational function of $m$
whose numerator is a nonzero polynomial of fixed degree.
Outside its finitely many zeros, this controls the additional
large primes that may divide the resultant. After division
by $d$, their contribution is $O_{a,b}(\log n)$; see
\eqref{eq:asymmetric-upper}. For equal parameters, we remove
the zero forced by odd parity, as in
\eqref{eq:reduced-gegenbauer}. The local bound then contains
a factorial ratio whose valuations are retained with their
signs; see \eqref{eq:symmetric-row}. Preserving this
cancellation gives the sharper uniform estimates when the
parameter grows proportionally to the degree
(Proposition~\ref{prop:asymptotic-range}).

For the Kaneko--Zagier families, passing to the primitive
polynomial is a substantial part of the argument.
After denominators are cleared, the coefficients have
a common divisor whose logarithm grows linearly with the
degree; see \eqref{eq:primitive-normalization} and
Proposition~\ref{prop:arithmetic-asymptotics}.
We estimate this content together with the moment
denominators, counting divisibility at each prime power
so as to preserve the cancellation between them
(Lemma~\ref{lem:carry-budget}). At $2$ and $3$, Taylor
recurrences from the differential equation control the
separation of roots, while digit sums associated with the
Newton polygons control the accumulated error
(Proposition~\ref{prop:bad-prime-resultant}).
Combining these estimates with the real bound gives the
positive gap in \eqref{eq:positive-gap}, uniformly in the
degree of a possible factor. The explicit cutoff and the
remaining finite checks are given in
Propositions~\ref{prop:explicit-kz-cutoff}
and~\ref{prop:finite-kz}, completing the argument in every
degree.

\subsection{Arithmetic applications}

To determine the Galois groups, we combine irreducibility with
the ramification detected by Newton polygons. At suitable primes,
the reciprocal polynomial has an Eisenstein factor, while
squarefree reduction shows that the remaining roots lie in an
unramified extension. In the tame cases, inertia therefore
contains a single cycle fixing every other root
(Lemma~\ref{lem:local-cycle} and
Proposition~\ref{prop:jacobi-cycle}).
On the transitive root set under consideration, a cycle moving
more than half the roots forces every possible block size to
divide its length (Lemma~\ref{lem:blocks}).
The formula for the cycle length would then confine all the
selected primes to one residue class for a larger modulus.
Explicit estimates for primes in arithmetic progressions
\cite[Corollary~1.7]{BMOR} and the Brun--Titchmarsh inequality
\cite[Theorem~2, Eq.~(1.10)]{MV} exclude such a block system;
see Proposition~\ref{prop:prime-criterion} and its extension
in Section~\ref{sec:general}.
The action is therefore primitive. Since the constructed
cycles leave at least three roots fixed, Jones's theorem
\cite[Corollary~1.3]{Jones}, recalled in
Theorem~\ref{thm:jones}, gives the alternating group.
A cycle of even length then gives the symmetric group.

When the roots occur in opposite pairs, we first determine
the permutation group on the pairs. A suitable power of an
inertia element interchanges the two roots in an odd, nonempty
proper set of pairs. Its conjugates generate all independent
sign changes, giving the full wreath product
(Lemma~\ref{lem:sign-kernel}).
The unequal integral case uses a different argument:
wild inertia supplies a prime-length cycle, to which Jordan's
theorem applies \cite[Theorem~1.1]{Jones}, and an odd valuation
of the discriminant supplies the odd permutation needed to
obtain the symmetric group
(Lemma~\ref{lem:galois-prime}).

The same ramification also controls intersections of splitting
fields. Let $L/\Q$ be one of the splitting fields under
consideration, and let $\widetilde E$ be the Galois closure
of a fixed number field $E$. The intersection
$L\cap\widetilde E$ is Galois over $\Q$. At a prime unramified
in $\widetilde E$, inertia acts trivially on this intersection.
The inertia subgroups used in our constructions normally
generate $\Gal(L/\Q)$, so $L\cap\widetilde E=\Q$;
see Lemma~\ref{lem:normal-generation} and, for unequal integral
parameters, the proof of Theorem~\ref{thm:main} in
Section~\ref{sec:galois}.
Since the selected primes tend to infinity with the degree,
they eventually avoid the finite set of primes ramified in
$\widetilde E$. This gives the assertions over every fixed
number field, including Corollary~\ref{cor:fixed-field}.
The disjointness of the large and complementary factors in
Theorem~\ref{thm:jacobi} follows from the same principle:
there the chosen inertia elements fix every root of the
complementary factor.

The modular origin of the Kaneko--Zagier families gives one
reason to study these Galois groups. These Jacobi polynomials
occur as polynomial factors of modular forms. When the
corresponding modular index is prime, the affine normalization
described in Section~\ref{sec:modular-interpretation}, together
with the prescribed factors at $j=0$ and $j=1728$, reduces to
the polynomial of supersingular $j$-invariants
\cite[Theorem~1]{KZ}.
Theorem~\ref{thm:kz-all} determines the arithmetic of the
characteristic-zero Jacobi factors in every positive degree.

A second application comes from spherical harmonics. The
dependence on the polar angle of a mode of degree $\ell$
and order $r$, with $0\le r\le\ell$, is, up to normalization,
\[
 (\sin\theta)^r P_\ell^{(r)}(\cos\theta),
\]
where $P_\ell^{(r)}$ denotes the $r$th derivative of the
Legendre polynomial
\cite[Eqs.~14.7.8 and~14.30.1]{DLMF}.
The circular harmonic cone $x^2+y^2-2z^2=0$ meets the unit
sphere at $\cos\theta=\pm1/\sqrt3$. Thus such a mode vanishes
on the cone precisely when $3x^2-1$ divides $P_\ell^{(r)}$.
Proposition~\ref{prop:fixed-factor} retains the arithmetic
of this prescribed factor in the resultant comparison.
It yields an effective absolute threshold beyond which,
uniformly for $0\le r\le\ell$,
\[
 P_\ell^{(r)}(1/\sqrt3)=0
 \quad\Longrightarrow\quad
 23r>4(\ell-1)
\]
(Corollary~\ref{cor:angular}).
This gives a simultaneous restriction on degree and order
in Agranovsky's question concerning homogeneous harmonic
polynomials divisible by this cone
\cite[\S6, Question~4]{Agranovsky}.
It complements the theorem of Mangoubi and Weller Weiser~\cite[Theorem~2]{MWW}:
the order must be even, and each fixed order admits at most
one degree.
The new estimate restricts the possible modes without giving
a complete classification.

The third application concerns hyperelliptic Jacobians.
The roots of a separable polynomial $f$ are the finite branch
points of the double cover with affine equation $y^2=f(x)$.
Their permutation action determines the Galois action on
the $2$-torsion of the Jacobian, connecting the arithmetic
of the polynomial with that of an abelian variety.
Let $J_f$ be the Jacobian of the smooth projective model.
If $f\in\Q[x]$ has degree $d\ge9$ and
$\Gal(f/\Q)=S_d$, Zarhin's theorem \cite[Theorem~2.4]{Zarhin02} gives
\[
 \End_{\overline{\Q}}(J_f)=\Z.
\]
His results also establish the Hodge and Tate conjectures
on every power of $J_f$, with the corresponding classes
generated by divisor classes, and the Mumford--Tate
conjecture for $J_f$ and its powers
\cite[Theorems~8.1 and~9.1, and Remark~9.2]{Zarhin02}.
Corollary~\ref{cor:jacobians} records these consequences
for our families. In particular, Theorem~\ref{thm:kz-all}
makes them available for the Jacobians associated with
all four Kaneko--Zagier families in every degree at least $9$.

\subsection{Organization}

Section~\ref{sec:comparison} develops the common resultant
comparison. Section~\ref{sec:integral} treats integral parameters
and derivatives, and Section~\ref{sec:kz-irreducibility} treats
Kaneko--Zagier irreducibility. Section~\ref{sec:galois-groups}
proves the Galois and splitting-field results.
Section~\ref{sec:applications} gives the applications.

The appendices contain local irreducibility theorems,
the explicit cutoff, and the finite calculations needed for the
all-degree statements. The finite Galois inputs are the published
conditional results in \cite{CHL,CHK}, used with their stated
degree and modular-index ranges.

\subsection*{Acknowledgements}

The author thanks the AxiomProver team for the original Lean formalization of this paper and is especially grateful to Ken Ono for sending it to him and drawing his attention to the team's corrections. The team observed that the prime-interval estimate used in the proof of Lemma~\ref{lem:odd-legendre-cycle} requires an integer argument and that Lemma~\ref{lem:finite-factor-degrees} requires positive degree. The formalization takes results cited from the literature as hypotheses.

For this revision, the author added supplementary Lean proofs for the Legendre polynomials and their first derivatives, including all $398$ finite cases in Lemma~\ref{lem:first-derivative-finite}. The author also supplied two independent programs that verify the finite inequalities in Proposition~\ref{prop:finite-kz}. The Lean proof of Theorem~\ref{thm:kz-all} assumes these inequalities. The proofs and verification programs are available at \url{https://github.com/KhaiHoanMath/arithmetic-of-jacobi-polynomials}.

\subsection*{Disclosure of AI assistance}

This work was developed by the author with assistance from ChatGPT (OpenAI). The author developed the overall approach and methods, building on the prior work cited in the manuscript. During the research and preparation of the manuscript, ChatGPT was used to discuss and revise mathematical arguments, improve the exposition, and search for relevant literature. It also assisted with all computations presented in the appendices.

\section{A common resultant comparison}\label{sec:comparison}

We now establish the estimates used in the irreducibility proofs.
The aim is to bound a factor resultant from below at the real place
and from above through its prime divisors. Both estimates must be
proportional to the smaller factor degree, so that the comparison
excludes every proper factor at once.

The common input is a Jacobi identity obtained from orthogonality.
We first use it to bound derivatives at the real roots and distances
between roots over a $p$-adic field. We then specialize to integral
weights, where the moment denominators admit a particularly simple
bound. The final two lemmas turn these estimates into a form ready
for the arithmetic comparison.

For a prime $p$, we normalize the valuation by $v_p(p)=1$, extend
it to $\overline{\Q}_p$, and set $v_p(0)=+\infty$. The Gauss norm
is the maximum of the absolute values of the coefficients. A
polynomial in $\Z_p[t]$ is primitive when its Gauss norm is $1$.
Globally, a primitive integral polynomial has coefficient gcd $1$.
All logarithms are natural, and empty products equal $1$.
Write $(z)_j=z(z+1)\cdots(z+j-1)$, with $(z)_0=1$.

\subsection{The Jacobi identity}

Fix real $\alpha,\beta>-1$, put $c=\alpha+\beta$, and let
$y_n(t)=P_n^{(\alpha,\beta)}(1-2t)$, where $n\ge1$.
The endpoint expansion is
\begin{equation}\label{eq:irred-coefficients}
 y_n(t)=\sum_{j=0}^n(-1)^j
 \binom{n+\alpha}{n-j}\binom{n+c+j}{j}t^j.
\end{equation}
Let $I$ denote integration against the probability measure
\[
 d\omega(t)=\frac{t^\alpha(1-t)^\beta}{B(\alpha+1,\beta+1)}\,dt,
 \qquad 0<t<1.
\]
Here $B(s,t)$ is Euler's beta function. Thus
\begin{equation}\label{eq:jacobi-moments}
 I(t^j)=\mu_j=\frac{(\alpha+1)_j}{(c+2)_j}.
\end{equation}
The roots of $y_n$ are simple and lie in $(0,1)$, and
$I(y_nq)=0$ for $\deg q<n$. These standard Jacobi properties follow
from Rodrigues' formula; see \cite[\S\S18.2--18.5]{DLMF}.
Define
\[
 Q_n(t)=I_s\left(\frac{y_n(t)-y_n(s)}{t-s}\right),\qquad
 \sigma(t)=t(1-t),\qquad \rho(t)=\alpha-ct.
\]
The difference quotient is a polynomial of degree at most $n-1$.
It gives a second polynomial whose value at a root can be compared
with the derivative of $y_n$. Positivity of the measure makes the
quadrature weights positive, and the probability normalization makes
their sum one. For rational
parameters the moments are rational, so $I$ also has a coefficientwise
interpretation over every field of characteristic zero. This is why
the following identity can be used at both the real and finite places.

\begin{lemma}\label{lem:jacobi-wronskian}
One has the polynomial identity
\begin{equation}\label{eq:jacobi-wronskian}
 (c+1)y_n^2+\sigma(y_n'Q_n-y_nQ_n')+\rho y_nQ_n=C_n,
 \qquad
 C_n=\frac{(\alpha+1)_n(\beta+1)_n}{n!(c+2)_{n-1}}>0.
\end{equation}
At every root $t_i$ of $y_n$,
\begin{equation}\label{eq:root-derivative-lower}
 |y_n'(t_i)|\ge2\sqrt{C_n},
 \qquad
 C_n\longrightarrow
 \frac{\Gamma(c+2)}{\Gamma(\alpha+1)\Gamma(\beta+1)}>0.
\end{equation}
\end{lemma}
\begin{proof}
Write $\tau=\sigma'+\rho=\alpha+1-(c+2)t$ and
$\Lambda=n(n+c+1)$. The Jacobi equation is
\begin{equation}\label{eq:unit-interval-jacobi-ode}
 \sigma y_n''+\tau y_n'+\Lambda y_n=0,
\end{equation}
Writing $d\omega=w(t)\,dt$, the density satisfies
$(\sigma w)'=\tau w$.
For real $t>1$, put
\[
 m(t)=I_s((t-s)^{-1}),\qquad
 q(t)=I_s(y_n(s)/(t-s)).
\]
Integration by parts gives
\[
 \sigma m'-\rho m=c+1,
 \qquad
 \sigma q''+(\sigma'-\rho)q'+(\Lambda+c)q=0.
\]
The endpoint terms vanish because $\alpha,\beta>-1$ and
$\sigma w$ vanishes at both endpoints. Substituting
$Q_n=y_nm-q$ gives
\[
 \sigma Q_n''+(\sigma'-\rho)Q_n'+(\Lambda+c)Q_n
       =2(c+1)y_n'.
\]
Differentiating the left side of \eqref{eq:jacobi-wronskian},
and using these two differential equations, shows that it is constant.
Orthogonality gives
\[
 q(t)=\frac{I(y_n^2)}{\operatorname{lc}(y_n)}t^{-n-1}
          +O(t^{-n-2}).
\]
Consequently the constant is $(2n+c+1)I(y_n^2)$, which equals
$C_n$ by the Jacobi norm formula. The identity is polynomial, so
its validity for $t>1$ proves it everywhere.

For a root $t_i$, let
\[
 \ell_i(t)=\frac{y_n(t)}{(t-t_i)y_n'(t_i)}.
\]
Since $\ell_i^2-\ell_i$ is divisible by $y_n$, with quotient of
degree at most $n-2$, orthogonality yields
\[
 \omega_i:=I(\ell_i)=I(\ell_i^2)>0.
\]
Also $\sum_i\ell_i=1$, so $\omega_i\le1$. The definition of
$Q_n$ gives $Q_n(t_i)=y_n'(t_i)\omega_i$. Thus
\[
 C_n=t_i(1-t_i)y_n'(t_i)^2\omega_i
       \le\tfrac14y_n'(t_i)^2.
\]
Finally, the stated limit follows by expressing the rising factorials
in $C_n$ as gamma quotients.
\end{proof}

The numbers $\omega_i$ in the proof are the Gaussian quadrature
weights. Indeed, if $\deg q\le2n-1$, division by $y_n$ gives
$q=y_ns+r$ with $\deg s,\deg r<n$. Orthogonality and interpolation
then give
\[
 I(q)=I(r)=\sum_{i=1}^n\omega_i r(t_i)
             =\sum_{i=1}^n\omega_i q(t_i).
\]
Thus a familiar formula for numerical integration supplies the
derivative estimate needed for the arithmetic argument; see also
\cite[\S3.5(v)]{DLMF}.

\subsection{Projective root distances}

A primitive polynomial over $\Z_p$ need not have integral roots:
its leading coefficient may be divisible by $p$. To include these
roots, we measure distances on the projective line. The quantity
$\delta_p(z,w)$ defined below is nonnegative and unchanged by
inversion. For integral roots it equals $v_p(z-w)$; a large value
means that the roots are close in the projective $p$-adic metric.
The valuation of a resultant is the sum of these quantities over
pairs of roots in different factors. We bound the sum from a single
root to all the others, then sum over the smaller factor.

\begin{lemma}[Projective moment bound]\label{lem:projective-moment-budget}
Let $p$ be a prime and let $H\in\mathbb Z_p[t]$ be primitive
(at least one coefficient is a unit) of
degree $n\ge1$. Let $I$ be a $\mathbb Q_p$-linear functional on the
polynomials of degree at most $2n-1$, with moments
$\mu_j=I(t^j)$, such that
\[
 I(Hq)=0\qquad(\deg q<n).
\]
Put
\[
 Q_H(t)=I_s\!\left(\frac{H(t)-H(s)}{t-s}\right),
 \qquad \sigma(t)=t(1-t).
\]
Suppose that some $C\in\mathbb Q_p^\times$ satisfies
\[
 \sigma(z)H'(z)Q_H(z)=C
 \qquad\text{at every root }z\text{ of }H.
\]
Write
\[
 \eta_p=\min_{0\le j\le2n-1}v_p(\mu_j),\qquad
 \mathcal E_p=v_p(C)-\eta_p.
\]
For distinct roots put
\[
 \delta_p(z,w)=v_p(z-w)-\min\{v_p(z),0\}-\min\{v_p(w),0\}.
\]
Every root satisfies the projective row bound
\begin{equation}\label{eq:projective-row-bound}
 \sum_{H(w)=0,\ w\ne z}\delta_p(z,w)\le\mathcal E_p.
\end{equation}
Consequently, for every primitive factorization $H=AB$ in
$\mathbb Z_p[t]$, with $k=\deg A$ and $R=\operatorname{Res}(A,B)$,
\[
 v_p(R)\le k\mathcal E_p.
\]
\end{lemma}

\begin{proof}
The root identity implies that the roots are simple and avoid $0,1$.
The projective distances $\delta_p(z,w)$ are nonnegative.
Multiplicativity of the Gauss norm gives, for every primitive
$F\in\mathbb Z_p[t]$,
\[
 v_p(\operatorname{lc}F)
       =-\sum_{F(w)=0}\min\{v_p(w),0\}.
\]
Thus the full projective row at an integral root $z$ of $H$ is
$v_p(H'(z))$. If $z$ is nonintegral, inversion preserves
$\delta_p$, and the same observation applied to the primitive
polynomial $\widehat H(X)=X^nH(1/X)$ identifies its row with
$v_p(\widehat H'(z^{-1}))$. The same formula gives
\begin{equation}\label{eq:projective-resultant}
 v_p(R)=\sum_{A(z)=0}\sum_{B(w)=0}\delta_p(z,w).
\end{equation}

The coefficients of $Q_H$ are $\mathbb Z_p$-linear combinations of the
moments. At an integral root, therefore,
$v_p(Q_H(z))\ge\eta_p$ and $v_p(\sigma(z))\ge0$, so the root
identity gives $v_p(H'(z))\le\mathcal E_p$.

For a nonintegral root put $b=z^{-1}$ and $r=v_p(b)>0$.
Extend $I$ linearly to the local field containing $b$. The polynomial
\[
 q(t)=\frac{H(t)}{1-bt}
\]
has integral coefficients in that field: if $H_j$ and $q_j$ denote
the coefficients of $H$ and $q$, respectively, coefficient comparison
gives $q_j=H_j+bq_{j-1}$ for $0\le j<n$, with $q_{-1}=0$.
For $0\le j<n$, orthogonality gives
$I(t^jq)=bI(t^{j+1}q)$, and hence
\[
 I(q)=b^nI(t^nq),\qquad
 \widehat Q(b)=-b^nI(q)=-b^{2n}I(t^nq),
 \quad \widehat Q(X)=X^{n-1}Q_H(1/X).
\]
Since $\deg(t^nq)=2n-1$, it follows that
$v_p(\widehat Q(b))\ge2nr+\eta_p$.
Reversing the root identity gives
\[
 (1-b)\widehat H'(b)\widehat Q(b)=Cb^{2n-1},
\]
so $v_p(\widehat H'(b))\le\mathcal E_p-r$.
Finally, all projective distances are nonnegative. Summing the $k$
full rows belonging to the roots of $A$ bounds the cross-factor sum
in \eqref{eq:projective-resultant} by $k\mathcal E_p$.
\end{proof}

The reflection $t=-X$ gives the same bound for root identities
with $\sigma(X)=X(X+1)$: the change of sign in the constant does
not affect its valuation. The estimate applies to subsets of roots
as well, since all omitted projective distances are nonnegative.

\subsection{Integral weights}

For integral parameters, we move the endpoints to $-1,0$ and
normalize the polynomial to have value one at $0$. The weight is
then a polynomial with integer coefficients, so its moment
denominators are controlled by a short interval of integers. The
constant in the root identity becomes a factorial ratio. Its
valuation may have either sign; retaining that sign will be
essential when the parameter is allowed to grow with the degree.

Fix nonnegative integers $a,b$ and an integer $n\ge1$. Put
\[
 B=\binom{n+a}{a},\quad S=2n+a+b+1,\quad
 T=\prod_{j=1}^{a+b}(n+j),
\]
and write
\[
 J(X)=P_n^{(a,b)}(2X+1)
 =\sum_{j=0}^n\binom{n+a}{n-j}\binom{n+a+b+j}{j}X^j,
 \qquad Q=J/B.
\]
Thus $J\in\Z[X]$ and $Q(0)=1$. With
$p_0(X)=X(X+1)$, define
\[
 I(g)=\int_{-1}^0g(t)(-t)^a(t+1)^b\,dt,
 \qquad
 V(X)=I\left(\frac{Q(X)-Q(t)}{X-t}\right).
\]
The quotient inside the integral is a polynomial. The moments of $I$
are rational, so we extend $I$ coefficientwise to polynomials over
any field of characteristic zero. This is the meaning of $I$ in the
local arguments below. Orthogonality gives $I(t^jQ)=0$ for $0\le j<n$.
Expanding the weight, which has integer coefficients, also gives
\begin{equation}\label{eq:moments}
 v_p(I(t^j))\ge-\lfloor\log_p(j+a+b+1)\rfloor
 \qquad(j\ge0).
\end{equation}
Indeed, the integral is an integral linear combination of reciprocals
of integers at most $j+a+b+1$.

\begin{lemma}\label{lem:weighted}
At every root $\alpha$ of $Q$,
\begin{equation}\label{eq:weighted}
 p_0(\alpha)Q'(\alpha)V(\alpha)=\mathcal C_n,
 \qquad
 \mathcal C_n=-\frac{(a!)^2n!(n+b)!}{(n+a)!(n+a+b)!}.
\end{equation}
\end{lemma}

\begin{proof}
Apply Lemma~\ref{lem:jacobi-wronskian} with parameters $a,b$,
then make the substitution $t=-X$. The present functional is
$B(a+1,b+1)$ times the probability functional in that lemma,
and $Q(X)=y_n(-X)/\binom{n+a}{a}$. Thus the root constant is
\[
 -\frac{B(a+1,b+1)}{\binom{n+a}{a}^2}
 \frac{(a+1)_n(b+1)_n}{n!(a+b+2)_{n-1}}
 =-\frac{(a!)^2n!(n+b)!}{(n+a)!(n+a+b)!}.
\]
\end{proof}

\begin{lemma}\label{lem:local}
Put
\[
 d_p=-\min_jv_p([X^j]Q),\qquad L_p=\lfloor\log_pS\rfloor.
\]
For every root $\alpha$ of $Q$,
\[
 \sum_{\substack{Q(\gamma)=0\\\gamma\ne\alpha}}
 \delta_p(\alpha,\gamma)\le L_p+2d_p+v_p(\mathcal C_n).
\]
In particular, this sum is at most
\begin{equation}\label{eq:local-coarse}
 \mathcal D_p=L_p+v_p(T).
\end{equation}
For $a=b=r$ one has the sharper bound
\begin{equation}\label{eq:local-symmetric}
 D_p=L_p+v_p\left(\frac{(n+r)!^2}{n!(n+2r)!}\right).
\end{equation}
Each bound applies at $p=2$ as well. The resultant of two factors normalized to Gauss norm $1$
and with disjoint root sets has valuation at most the smaller
degree times the corresponding bound.
\end{lemma}

\begin{proof}
Since $Q(0)=1$ and $BQ\in\Z[X]$, one has $0\le d_p\le v_p(B)$.
The polynomial $u=p^{d_p}Q$ has Gauss norm $1$, and its root
identity has constant $p^{2d_p}\mathcal C_n$. The moment bound
\eqref{eq:moments}, through degree $2n-1$, is at least $-L_p$.
The reflected form of Lemma~\ref{lem:projective-moment-budget}
therefore gives the first assertion, including nonintegral roots.
For $U_a=\prod_{j=1}^a(n+j)$ and
$V_a=\prod_{j=1}^a(n+b+j)$,
\[
 2d_p+v_p(\mathcal C_n)\le2v_p(B)+v_p(\mathcal C_n)
 =v_p(U_a)-v_p(V_a)\le v_p(T).
\]
When $a=b=r$, the middle expression is the valuation of
$(n+r)!^2/(n!(n+2r)!)$. Summing over the smaller root set
bounds the resultant; removing roots only decreases the sums.
\end{proof}

\subsection{The real estimate and the prime 2}

We use the Chebyshev--Hadamard estimate appearing in
\cite[Lemma~5.5]{Deng}.

\begin{lemma}\label{lem:interval-vandermonde}
For an integer $k\ge1$ and $t_1,\ldots,t_k\in[0,1]$,
\[
 \prod_{i<j}|t_i-t_j|^2\le k^k2^{-2(k-1)^2}.
\]
\end{lemma}
\begin{proof}
Let $T_j(\cos\theta)=\cos(j\theta)$ be the Chebyshev polynomial.
For $j\ge1$, the polynomial $q_j(t)=2^{1-2j}T_j(2t-1)$ is
monic and has absolute value at most $2^{1-2j}$ on $[0,1]$;
put $q_0=1$. The determinant of $(q_{j-1}(t_i))_{i,j=1}^k$
is the Vandermonde determinant. Hadamard's inequality bounds its
absolute value by
\[
 k^{k/2}\prod_{j=1}^{k-1}2^{1-2j}
 =k^{k/2}2^{-(k-1)^2}.
\]
Squaring proves the assertion.
\end{proof}

We combine the derivative bound with this Vandermonde estimate.
The values of an integral factor at the two endpoints also bound
its leading coefficient from below. The next lemma records both
contributions and separates the prime $2$. We keep track of the
powers of $2$ introduced by $x=2X+1$, since their cancellation
determines the strength of the bound for the odd part of the
resultant.

\begin{lemma}\label{lem:real-resultant}
Let $h\in\Z[1/2][x]$ have degree $N$, positive leading coefficient,
and distinct nonzero roots in $(-1,1)$. Suppose $b_h=h(1)$ is a
positive integer and $|h'(\alpha)|>\gamma>0$ at every root.
Assume that the projective row sums for $h(2X+1)$ at $2$ are at most
$D_2$. Choose $\kappa>0$ so that $H=\kappa h$ is primitive integral,
and suppose
\[
 H=FG,\qquad 1\le d=\deg F\le e=\deg G,\qquad d+e=N,
\]
with $F,G\in\Z[x]$ primitive and with positive leading coefficients. Put
\[
 R=\Res(F,G),\quad \mathcal R=|R|/2^{v_2(R)},\quad
 u=\lc F,\quad c_F=\min_jv_2([X^j]F(2X+1)),
\]
and $V_F=\prod_{i<j}(\alpha_i-\alpha_j)^2$ for the roots of $F$.
Then $u>2^{c_F}$ and
\begin{equation}\label{eq:master-resultant}
 \frac{\log\mathcal R}{d}>
 e\log2-\frac{\log V_F}{d}+\log\frac{\gamma}{b_h}
 -D_2\log2+\frac{e-d}{d}\log\frac{u}{2^{c_F}}.
\end{equation}
In particular,
\begin{equation}\label{eq:master-coarse}
 \frac{\log\mathcal R}{d}>
 (N-3)\log2-\log d+\log\frac{\gamma}{b_h}-D_2\log2.
\end{equation}
\end{lemma}

\begin{proof}
The derivative and resultant formulas give
\[
 |R|=\frac{u^{e-d}}{V_F}\prod_{F(\alpha)=0}|H'(\alpha)|
 >\frac{u^{e-d}}{V_F}(\kappa\gamma)^d.
\]
Define $c_G$ as for $c_F$. The affine resultant formula and the
normalized row bound give
\[
 v_2(R)\le ec_F+dc_G-de+dD_2.
\]
Both $F(1)$ and $F(-1)$ are nonzero integers divisible by $2^{c_F}$.
Therefore
\[
 u^2=\frac{|F(1)F(-1)|}{\prod_i(1-\alpha_i^2)}>2^{2c_F}.
\]
After subtracting the binary estimate, the first two displays yield
the right side of \eqref{eq:master-resultant}, with $-\log b_h$
replaced by $\log\kappa-(c_F+c_G)\log2$.

Write $\kappa=2^k/c$, with $c$ positive and odd. Since
$h\in\Z[1/2][x]$ and $H$ is primitive, $c$ divides the odd part of
$b_h$. If $t_2$ is the minimum binary coefficient valuation of
$h(2X+1)$, Gauss multiplicativity gives $c_F+c_G=k+t_2$, whereas
$t_2\le v_2(b_h)$. Thus
\[
 \log\kappa-(c_F+c_G)\log2
 =-\log c-t_2\log2\ge-\log b_h.
\]
This proves \eqref{eq:master-resultant}.

Finally, rescaling Lemma~\ref{lem:interval-vandermonde} from
$[0,1]$ to $[-1,1]$ gives
\[
 V_F\le d^d2^{-(d-1)(d-2)},\qquad
 -\frac{\log V_F}{d}\ge(d-3)\log2-\log d.
\]
The same inequality holds for $d=1$. Dropping the nonnegative last
term in \eqref{eq:master-resultant} proves \eqref{eq:master-coarse}.
\end{proof}

\subsection{A bound for the prime powers}

The local moment bound contains $\lfloor\log_p S\rfloor$.
After summing over primes, these terms are naturally expressed
through least common multiples. The following observation groups
the higher prime powers into a single smaller interval. For an
integer $B\ge1$ and a real number $S\ge1$, write
\[
 \Psi_B(S)=\sum_{p\le B}\lfloor\log_pS\rfloor\log p.
\]

\begin{lemma}\label{lem:lcm}
If $B\ge9$ and $S\le5B$, then
\begin{equation}\label{eq:lcm-bound}
 \Psi_B(S)\le
 \log\lcm(1,\ldots,B)+\log\lcm(1,\ldots,\lfloor\sqrt S\rfloor)
 <(B+\sqrt S)\log3.
\end{equation}
\end{lemma}

\begin{proof}
For $p\le B$, let $l=\lfloor\log_pB\rfloor$ and
$k=\lfloor\log_pS\rfloor$. Since $B\ge9$, one has $p^l\ge5$,
including $p=2,3$. Hence
\[
 S\le5B<5p^{l+1}\le p^{2l+1},
\]
so $k\le2l$ and
$k\le l+\lfloor k/2\rfloor=l+\lfloor\log_p\sqrt S\rfloor$.
Multiplying these primewise inequalities proves the first bound.
The second follows from Hanson's estimate
$\lcm(1,\ldots,t)<3^t$ for integers $t\ge1$ \cite{Hanson}.
\end{proof}
\section{Integral parameters and derivatives}\label{sec:integral}

The local bounds of Section~\ref{sec:comparison} involve primes
up to $2n+O(1)$ when the parameters are fixed. Summing over all
these primes gives an upper bound too large for the real estimate
to contradict. We must therefore show that most primes above the
half-degree $m=\lfloor n/2\rfloor$ do not divide a factor resultant.

For unequal integral parameters, the exceptional primes above
$m+O_{a,b}(1)$ divide values of fixed polynomials. Their total
logarithmic contribution is only $O_{a,b}(\log n)$. For equal
parameters, in the range of the criterion below, parity excludes
every odd prime above $m$.
In that case, retaining the exact valuations of a factorial ratio
also allows the parameter to grow with the degree. These two
arguments prove the irreducibility assertions of
Theorems~\ref{thm:main} and~\ref{thm:symmetric}, respectively.

\subsection{Unequal parameters}\label{sec:asymmetric}

By reflection,
$P_n^{(b,a)}(-x)=(-1)^nP_n^{(a,b)}(x)$, we may suppose $a>b\ge0$.
Write $q=a-b$, $s=a+b$ and $n=2m+\delta$, where
$\delta\in\{0,1\}$, and put $b_m=4^{-m}\binom{2m}{m}$.
The value at zero supplies the main arithmetic input. For each
parity, Lemma~\ref{lem:center} will express it as
$(-1)^mb_m$ times a rational function $A_\delta(m)/D_\delta(m)$
whose numerator and denominator have degrees depending only on
$a,b$. The numerator identifies the additional large primes that
must be retained in the resultant estimate, while the denominator
enters the real lower bound. Both therefore cost only
$O_{a,b}(\log n)$ in the final
comparison, away from the finitely many zeros of the numerator.
The following formulas make this dependence explicit:
\begin{align}
 H_\delta&=\left\lfloor\frac{q-\delta}{2}\right\rfloor,
 &K_\delta&=b+\delta+H_\delta,\label{eq:HK}\\
 D_\delta(z)&=2^\delta\prod_{j=1}^{K_\delta}(z+j),\label{eq:D}\\
 A_\delta(z)&=\left(\prod_{j=1}^{b+\delta}(2z+j)\right)
 \sum_{h=0}^{H_\delta}(-1)^h\binom q{\delta+2h}
 z^{\underline h}\prod_{j=b+\delta+h+1}^{K_\delta}(z+j),\label{eq:A}
\end{align}
where $z^{\underline h}=z(z-1)\cdots(z-h+1)$ and
$z^{\underline0}=1$. The polynomials $A_\delta,D_\delta$ have integer
coefficients and degree at most $K_\delta$, and $A_\delta(0)>0$.
Recall from Section~\ref{sec:comparison} that
\[
 B=\binom{n+a}{a},\qquad S=2n+s+1,\qquad
 T=\prod_{j=1}^{s}(n+j).
\]
Put
\begin{equation}\label{eq:asymmetric-constants}
 M=m+s+2,\qquad
 c_{a,b}=\frac{(a+b+1)^{a+b+1}}
 {(2a+1)^{a+1/2}(2b+1)^{b+1/2}}.
\end{equation}

\begin{theorem}\label{thm:criterion}
Suppose $m\ge2s+10$ and $A_\delta(m)\ne0$. If
\begin{equation}\label{eq:criterion}
\begin{split}
 n\log2>{}&(M+\sqrt S)\log3+3\log T+\log|A_\delta(m)|\\
 &+\log\left(\frac{8(m+1)SBD_\delta(m)}{c_{a,b}\sqrt m}\right),
\end{split}
\end{equation}
then $P_n^{(a,b)}$ is irreducible over $\Q$.
\end{theorem}

For fixed $a,b$, the right side is
$(n/2)\log3+O_{a,b}(\sqrt n)$, while the left side is $n\log2$.
The positive difference $\log2-\tfrac12\log3$ proves eventual
irreducibility once the exceptional primes have been controlled.

\subsubsection{Central values and reduction}

We first prove the central-value formula. Its exact numerator is
needed to determine the exceptional primes; a lower bound for its
absolute value will also give a lower bound for derivatives at roots.

\begin{lemma}\label{lem:center}
For $m\ge0$,
\begin{equation}\label{eq:center}
 P_{2m+\delta}^{(a,b)}(0)
 =(-1)^mb_m\frac{A_\delta(m)}{D_\delta(m)}.
\end{equation}
If $m\ge1$ and $A_\delta(m)\ne0$, then
\[
 |P_{2m+\delta}^{(a,b)}(0)|
 \ge\frac1{2\sqrt m\,D_\delta(m)}.
\]
In particular, each parity has only finitely many zero central values.
\end{lemma}

\begin{proof}
Rodrigues' formula \cite[(4.3.1)]{Szego} gives
\[
 P_n^{(a,b)}(0)=(-1)^n2^{-n}[t^n](1-t)^{n+a}(1+t)^{n+b}.
\]
Expanding $(1-t)^q(1-t^2)^{n+b}$ yields
\[
 P_n^{(a,b)}(0)=(-1)^m2^{-2m-\delta}
 \sum_{h=0}^{H_\delta}(-1)^h\binom q{\delta+2h}
 \binom{2m+\delta+b}{m-h}.
\]
The identity
\[
 \frac{\binom{2m+\delta+b}{m-h}}{\binom{2m}{m}}
 =\frac{\prod_{j=1}^{b+\delta}(2m+j)\,m^{\underline h}}
 {\prod_{j=1}^{b+\delta+h}(m+j)}
\]
proves \eqref{eq:center}. It also holds for $h>m$, with both sides
zero. At $m=0$, the summands with $h>0$ vanish and the remaining
summand is positive, so $A_\delta$ is nonzero. Finally,
$|A_\delta(m)|\ge1$ when nonzero, and $b_m\ge1/(2\sqrt m)$ for
$m\ge1$. The latter inequality follows by induction from $b_1=1/2$
and $b_{m+1}/b_m=(2m+1)/(2m+2)$.
\end{proof}

For example, when $(a,b)=(1,0)$ the formulas reduce to
\[
 A_0(z)=D_0(z)=1,\qquad
 A_1(z)=2z+1,\qquad D_1(z)=2(z+1).
\]
Thus the central value never vanishes in either parity in this case.

The next lemma controls roots whose reduction is nonzero.
For an odd prime $p$, a bar denotes coefficientwise reduction modulo
$p$, and roots of reduced polynomials are taken in an algebraic
closure of $\F_p$.

\begin{lemma}\label{lem:reduction}
Suppose $m\ge2s+10$, and let $p>M$ be an odd prime with $p\nmid T$.
Write $n=\epsilon p+j$, with $\epsilon\in\{0,1\}$ and $0\le j<p$.
Then $j+s<p$ and
\[
 \overline{P_n^{(a,b)}}=x^{\epsilon p}\overline{P_j^{(a,b)}}.
\]
Every nonzero root of this reduction is simple, and the values at
$1$ and $-1$ are nonzero.
\end{lemma}

\begin{proof}
We have $n+s<2p$. If $j+s\ge p$, a multiple of $p$ would occur
between $n+1$ and $n+s$, contradicting $p\nmid T$.
For $\epsilon=0$ the congruence is immediate. For $\epsilon=1$,
Lucas's binomial congruence applied to \eqref{eq:jacobi-definition} leaves only
the ranges $0\le k\le j$ and $p\le k\le p+j$. Indeed, an
intermediate surviving term would require both $k\ge p-a$ and
$k\le j+b$, contrary to $j+s<p$. The two ranges sum to
\[
 2^{-p}\{(x+1)^p+(x-1)^p\}\overline{P_j^{(a,b)}}
 =x^p\overline{P_j^{(a,b)}}.
\]
The endpoint values of $P_j^{(a,b)}$ are the units
$\binom{j+a}{a}$ and $(-1)^j\binom{j+b}{b}$. Its reduction has
degree less than $p$ and satisfies
\begin{equation}\label{eq:integral-jacobi-ode}
 (1-x^2)y''+\{b-a-(a+b+2)x\}y'+j(j+s+1)y=0.
\end{equation}
A root away from the endpoints with multiplicity $t\ge2$ would
produce a nonzero lowest-order term proportional to $t(t-1)$,
because $t<p$. Thus every such root is simple. Multiplication by
$x^p$ can introduce a multiple root only at zero.
\end{proof}

\subsubsection{Comparison and effectivity}

\begin{proof}[Proof of Theorem~\ref{thm:criterion}]
Put $f=P_n^{(a,b)}$. Suppose $H=\kappa f=FG$ is a primitive integral
factorization, with positive leading coefficients and
$1\le d=\deg F\le e=\deg G$. Write
$R=\Res(F,G)$ and $\mathcal R=|R|/2^{v_2(R)}$.

We first estimate the odd part from above. For odd $p$, the substitution
$x=2X+1$ preserves the Gauss norm and resultant valuation, so
Lemma~\ref{lem:local} gives $v_p(R)\le d(L_p+v_p(T))$.
Moreover,
\begin{equation}\label{eq:asymmetric-prime-support}
 p>M,\quad p\nmid TA_\delta(m)\quad\Longrightarrow\quad v_p(R)=0.
\end{equation}
To see this, $f(1)=B$ is a unit and $f$ is integral at $p$, so
$\kappa$ is a unit. Lemma~\ref{lem:reduction} makes every unit root
simple modulo $p$. Since $p^2>2n+s$, the leading coefficient
$2^{-n}\binom{2n+s}{n}$ has valuation $0$ or $1$. Formula
\eqref{eq:center} gives the same alternatives for $f(0)$, because
$p>m+K_\delta$ and $p\nmid A_\delta(m)$.

For a polynomial $U$ of Gauss norm $1$ with nonzero constant term,
Gauss multiplicativity gives
\begin{equation}\label{eq:root-masses}
 v_p(U(0))=\sum_\alpha\max\{v_p(\alpha),0\},\qquad
 v_p(\lc U)=-\sum_\alpha\min\{v_p(\alpha),0\}.
\end{equation}
For each primitive integral factor, both sums are nonnegative
integers. A total of at most $1$ therefore cannot be positive in
both $F$ and $G$. Thus the two factors cannot both have a root
reducing to zero, and they cannot both have a root reducing to
infinity. Every other residue class contains at most one root of $f$.
Every cross projective distance vanishes, proving
\eqref{eq:asymmetric-prime-support}.

The terms $v_p(T)$ contribute at most $\log T$ after division by $d$.
For $p\le M$, Lemma~\ref{lem:lcm} bounds the remaining contribution
by $(M+\sqrt S)\log3$. For $p>M$ we have $p^2>S$, and a surviving
prime divides $TA_\delta(m)$. Hence
\begin{equation}\label{eq:asymmetric-upper}
 \frac{\log\mathcal R}{d}
 <(M+\sqrt S)\log3+2\log T+\log|A_\delta(m)|.
\end{equation}

We now need a lower bound for $|f'(\alpha)|$ at every root
$\alpha$ of $f$. The central-value formula gives information at
one point; the differential equation extends it to the whole
interval. Set $\lambda=n(n+s+1)$ and define the energy
\[
 E(x)=(1-x^2)f'(x)^2+\lambda f(x)^2.
\]
At a root $\alpha$, this is $(1-\alpha^2)f'(\alpha)^2$, so a
positive lower bound for $E$ gives the required derivative bound.
The differential equation gives
\[
 E'(x)=2\{q+(s+1)x\}f'(x)^2.
\]
The minimum on $[-1,1]$ occurs at $x_0=-q/(s+1)$. On $[x_0,0]$,
\[
 0\le E'(x)/E(x)\le\frac{2\{q+(s+1)x\}}{1-x^2}.
\]
The integral of the last expression shows that
$\sqrt{E(x_0)}\ge c_{a,b}\sqrt{E(0)}$, with the explicit constant
in \eqref{eq:asymmetric-constants}. Lemma~\ref{lem:center} yields
\[
 \sqrt{E(x_0)}\ge c_{a,b}\sqrt\lambda\,|f(0)|
 \ge\frac{c_{a,b}\sqrt m}{D_\delta(m)}=:\gamma.
\]
At every root $\alpha$, $\alpha\ne0$ and $|f'(\alpha)|>\gamma$.
Apply Lemma~\ref{lem:real-resultant} with $b_h=B$ and
$D_2=\mathcal D_2$. Since $d\le m+1$ and
$\mathcal D_2\log2\le\log S+\log T$, it gives
\begin{equation}\label{eq:asymmetric-lower}
 \frac{\log\mathcal R}{d}>
 n\log2-\log T-
 \log\left(\frac{8(m+1)SBD_\delta(m)}{c_{a,b}\sqrt m}\right).
\end{equation}
The two bounds contradict \eqref{eq:criterion}.
\end{proof}

\begin{proof}[Irreducibility in Theorem~\ref{thm:main}]
Fix $a>b$. A polynomial root bound gives an effective point beyond
which neither $A_0(m)$ nor $A_1(m)$ vanishes. Their coefficients give
explicit constants such that
\[
 |A_\delta(m)|\le C(m+1)^{K_\delta},\quad
 D_\delta(m)\le2(m+K_\delta)^{K_\delta},\quad
 T\le(n+s)^s,\quad B\le(n+a)^a.
\]
For sufficiently large $m$, the difference between the left and right
sides of \eqref{eq:criterion} is at least
\[
 m\log(4/3)-C_1\sqrt m-C_2\log(m+1)-C_3,
\]
with effective nonnegative constants depending on $a,b$. After
division by $m$, every error term decreases to zero. Successive
integer doubling therefore gives an effective bound valid for all
subsequent $m$. The constants and logarithms can be bounded by rational
interval estimates; take the larger bound for the two parities.
Reflection treats $b>a$.
\end{proof}

\subsection{Symmetric parameters}\label{sec:symmetric}

Equal parameters give an additional restriction on the reductions.
At an odd prime, parity pairs every nonzero root with its distinct
negative. A multiple root away from the central point and the
endpoints must therefore occur in an opposite pair. The differential equation forces these
multiplicities to be large, and the degree then excludes them at
the primes under consideration. This is what removes the exceptional
large primes present in the unequal-parameter argument.

Let $r,m$ be integers with $r\ge0$, $m\ge1$, and put
$n=2m+\delta$ for $\delta\in\{0,1\}$. We use the derivative
normalization from \eqref{eq:derivative-gegenbauer}:
\[
 f(x)=P_{n+r}^{(r)}(x)=(2r-1)!!\,C_n^{r+1/2}(x).
\]
Set
\begin{equation}\label{eq:symmetric-constants}
 S=2n+2r+1,\quad M=m+r,\quad
 A=\frac{(n+r)!^2}{n!(n+2r)!},\quad
 K=\binom{m+r}{r},\quad B=\binom{n+r}{r},
\end{equation}
and put $A_{\le m}=\prod_{p\le m}p^{v_p(A)}$.
We retain the sign of every exponent $v_p(A)$: cancellation between
the numerator and denominator will determine the uniform range in
$r$. Although $A\le1$, the truncated product $A_{\le m}$ need not
be at most $1$. For example, $r=1$, $m=4$, and $\delta=0$ give
$A=9/10$ but $A_{\le m}=9/2$.

\subsubsection{The exact criterion}

\begin{theorem}\label{thm:symmetric-exact}
Suppose $m\ge\max\{2r,4\}$. If
\begin{equation}\label{eq:symmetric-exact}
 m\log4>\Psi_m(S)+\log(8\sqrt m)+\log K+\log A_{\le m},
\end{equation}
then $C_{2m+\delta}^{r+1/2}(x)/x^\delta$ is irreducible over $\Q$.
When also $M\ge9$, the following elementary condition is sufficient:
\begin{equation}\label{eq:symmetric-elementary}
 m\log4>(M+\sqrt S)\log3+\log(8\sqrt m)
 +\log K+2\log B-\log\binom{2r}{r}.
\end{equation}
\end{theorem}

\begin{proof}
Apply Lemma~\ref{lem:local} to
$Q(X)=B^{-1}P_n^{(r,r)}(2X+1)$. It gives the row bound
\begin{equation}\label{eq:symmetric-row}
 D_p=L_p+v_p(A),\qquad L_p=\lfloor\log_pS\rfloor.
\end{equation}
This number is nonnegative, since it bounds a nonnegative row sum.
Removing the zero forced by odd parity removes only nonnegative
distances, so the same bound holds for $h=f/x^\delta$ after the
shift. Equivalently, the removed factor $2X+1$ has Gauss norm $1$
at every prime, including $2$.

Differentiation of the Legendre coefficient formula gives
\begin{equation}\label{eq:gegenbauer-coefficients}
 f(x)=2^{-n-r}\sum_{k=0}^m
 \frac{(-1)^k(2n+2r-2k)!}{k!(n+r-k)!(n-2k)!}\,x^{n-2k}.
\end{equation}
Its nonzero central coefficient is
\begin{equation}\label{eq:gegenbauer-center}
 [x^\delta]f(x)=\frac{(-1)^m}{2^mm!}
 (2m+2r+2\delta-1)!!.
\end{equation}
The other roots are simple and lie in $(-1,1)$. Choose $\kappa>0$
so that $H=\kappa h$ is primitive integral, and suppose
$H=FG$ with $1\le d=\deg F\le e=\deg G$, $d+e=2m$.
As before, put $R=\Res(F,G)$ and $\mathcal R=|R|/2^{v_2(R)}$.

We claim that
\begin{equation}\label{eq:symmetric-prime-support}
 p>m\quad\Longrightarrow\quad v_p(R)=0.
\end{equation}
Such a prime is odd, exceeds $2r$, and satisfies $p^2>S$.
The primitive polynomial $x^\delta H=\kappa f$ has parity $n$ and
satisfies the Gegenbauer equation. Its nonzero reduction has degree
at most $n<2p$. A multiple root away from $0,1,-1$ would have
multiplicity at least $p$, by the lowest-order term of the differential
equation. Parity gives the same multiplicity at its distinct negative,
forcing degree at least $2p$. Thus no such multiple root exists.

For the four remaining residue classes $0,1,-1,\infty$, each total
root-valuation sum is at most $1$. Indeed, before primitive
normalization, the central coefficient in \eqref{eq:gegenbauer-center}
has valuation at most $1$: its odd numerator factors are smaller than
$3p$, and $p>m$. The leading coefficient is
\[
 \lc(f)=\frac{(2n+2r-1)!!}{n!}.
\]
Its odd numerator factors are smaller than $5p$. If $p\le n$, at
most $p$ and $3p$ occur and $n!$ contributes one $p$; if $p>n$,
the numerator factors are smaller than $3p$ and $n!$ is a unit.
The valuation is at most $1$ in either case. Finally,
\begin{equation}\label{eq:symmetric-endpoint}
 f(1)=\frac{(n+2r)!}{2^rr!n!}
 =\binom{n+2r}{2r}(2r-1)!!,
 \qquad |f(-1)|=f(1).
\end{equation}
The numerator product contains $2r<p$ consecutive integers, all
smaller than $p^2$, so its valuation is at most $1$.
Primitive normalization can only decrease these odd-prime valuations.
Formula \eqref{eq:root-masses}, and its translates at $1$ and $-1$,
now applies. Each primitive factor contributes a nonnegative integer
to each sum, so no one of the four classes is shared by $F$ and $G$.
All cross distances vanish. This proves
\eqref{eq:symmetric-prime-support}.

The row bound at the remaining odd primes gives
\begin{equation}\label{eq:symmetric-upper}
 \frac{\log\mathcal R}{d}
 \le\Psi_m(S)+\log A_{\le m}-D_2\log2.
\end{equation}
For the real estimate, the equation for $f$ is
\[
 (1-x^2)f''-(2r+2)xf'+n(n+2r+1)f=0.
\]
With $\lambda=n(n+2r+1)$ and
$E=(1-x^2)f'^2+\lambda f^2$, it gives
$E'=2(2r+1)xf'^2$. Hence $E(x)\ge E(0)=\gamma^2$, where
$\gamma^2=\lambda f(0)^2+f'(0)^2>0$.
At each nonzero root $\alpha$,
$|h'(\alpha)|=|f'(\alpha)|/|\alpha|^\delta>\gamma$.
The exact central and endpoint formulas imply
\begin{equation}\label{eq:symmetric-energy}
 \frac{\gamma}{f(1)}=
 \begin{cases}
 \sqrt{2m(2m+2r+1)}\,b_m/K,&\delta=0,\\
 (2m+1)b_m/K,&\delta=1,
 \end{cases}
 \qquad \frac{\gamma}{f(1)}>\frac{\sqrt m}{K}.
\end{equation}
Lemma~\ref{lem:real-resultant}, with $N=2m$ and $d\le m$, therefore
gives
\begin{equation}\label{eq:symmetric-lower}
 \frac{\log\mathcal R}{d}>
 m\log4-\log(8\sqrt m)-\log K-D_2\log2.
\end{equation}
Comparing \eqref{eq:symmetric-upper} and
\eqref{eq:symmetric-lower} proves \eqref{eq:symmetric-exact}.
In particular the binary contributions cancel exactly.

For the elementary criterion, retain the nonnegativity of $D_p$.
Every prime occurring in $A$ divides
$\prod_{j=1}^{2r}(n+j)$. A prime $p>M$ in this product satisfies
$2p>n+2r$ and therefore belongs to $(n,n+2r]$. Thus
\begin{equation}\label{eq:symmetric-enlarge}
 \Psi_m(S)+\log A_{\le m}
 \le\Psi_M(S)+\sum_{n<p\le n+2r}\log p+\log A.
\end{equation}
The difference is a sum of nonnegative terms $D_p\log p$.
Every prime in $(n,n+2r]$ divides $\binom{n+2r}{2r}$, and
\[
 \log\binom{n+2r}{2r}+\log A
 =2\log B-\log\binom{2r}{r}.
\]
Since $S\le5M$, Lemma~\ref{lem:lcm} proves
\eqref{eq:symmetric-elementary}.
\end{proof}

\subsubsection{Uniform linear ranges}

The elementary criterion first gives an explicit range by bounding
each additional term relative to $m$. We then use the exact
factorial ratio to obtain a sharper asymptotic range.

\begin{proof}[The explicit assertion of Theorem~\ref{thm:symmetric}]
Suppose $m\ge\max\{100r,1000\}$. Then $S\le4.05m$.
Using $\binom uv\le(eu/v)^v$ and monotonicity in $r$, we have
\[
 \frac{\log K+2\log B}{m}
 \le\frac{3+\log101+2\log201.1}{100}<0.1824.
\]
For $r=0$ the left side is zero. Also
\[
 \frac{(M-m)\log3}{m}<0.011,\qquad
 \frac{\sqrt S\log3}{m}<\frac{1.1\cdot2.1}{31}<0.0746,\qquad
 \frac{\log(8\sqrt m)}m<0.006.
\]
These bounds hold at $m=1000$ and decrease thereafter where required.
Discarding $-\log\binom{2r}{r}$ in
\eqref{eq:symmetric-elementary}, the excess over $m\log3$, divided
by $m$, is less than
\[
 0.1824+0.011+0.0746+0.006=0.274<\frac27<\log(4/3).
\]
The last inequality follows from the power series for
$2\operatorname{arctanh}(1/7)$. The criterion applies.
\end{proof}

When $r$ is proportional to $m$, the terms involving $K$ and $A$
are also of order $m$, so they affect the leading constant in the
comparison. For $r\ge1$, write $t=m/r$. The next lemma shows that the arithmetic
side of \eqref{eq:symmetric-exact}, divided by $m$, is
$1+g(t)/t+o(1)$, whereas the real side is $\log4$.
Thus $1+g(t)/t<\log4$ supplies the required asymptotic gap.
The same calculation will later be used to
exclude a prescribed factor. Define
\begin{equation}\label{eq:g-definition}
 g(t)=(4t+2)\log(t+1/2)-3t\log t-(t+1)\log(t+1).
\end{equation}

\begin{lemma}\label{lem:uniform-asymptotics}
Uniformly for integers $1\le r\le m/2$ and $\delta\in\{0,1\}$,
as $m\to\infty$,
\begin{equation}\label{eq:uniform-asymptotics}
 \Psi_m(S)+\log K+\log A_{\le m}
 =m+r\,g(m/r)+o(m).
\end{equation}
The error has an effective bound independent of $r$ and $\delta$.
Moreover, $g(t)/t$ is strictly decreasing for $t\ge2$.
\end{lemma}

\begin{proof}
The first prime powers contribute $\vartheta(m)=m+o(m)$.
The higher powers contribute $O(\sqrt m\log m)=o(m)$, uniformly,
since $S\le5m+3$. For $p>m$, one has $p^2>n+2r$ for large $m$,
and $3p>n+2r$. The only possible multiples in
\[
 A=\frac{\prod_{j=1}^r(n+j)}{\prod_{j=r+1}^{2r}(n+j)}
\]
are $p$ and $2p$. Consequently
\begin{align}
 \log A-\log A_{\le m}
 ={}&2\vartheta(n+r)-\vartheta(n)-\vartheta(n+2r)\notag\\
 &+2\vartheta((n+r)/2)-\vartheta(n/2)
 -\vartheta((n+2r)/2).\label{eq:large-A}
\end{align}
Every argument is between fixed positive multiples of $m$. The prime
number theorem gives a uniform $o(m)$ error, and the linear terms
cancel. An explicit prime number theorem, for example
\cite[Theorem~5.2]{Dusart10}, makes this error effective.

The elementary factorial estimate
$\log(k!)=k\log k-k+O(\log(k+1))$, obtained by integrating
$\log x$, gives uniformly
\[
 \log K+\log A=r\,g(m/r)+O(\log m).
\]
Indeed, for $t=m/r$, the two main terms divided by $r$ are
\[
 (t+1)\log(t+1)-t\log t
\]
and
\[
 2(2t+1)\log(2t+1)-2t\log(2t)-(2t+2)\log(2t+2).
\]
Their sum is $g(t)$; changing $\delta$ changes the error by
$O(\log m)$. Combining the estimates proves
\eqref{eq:uniform-asymptotics}. Finally,
\[
 tg'(t)-g(t)=\log\frac{t+1}{(t+1/2)^2}<0\qquad(t\ge2).
\]
\end{proof}

\begin{proposition}\label{prop:asymptotic-range}
Let $c\ge2$ satisfy
\begin{equation}\label{eq:asymptotic-condition}
 c(\log4-1)>g(c).
\end{equation}
For all sufficiently large $m$, uniformly for $0\le r\le m/c$,
both polynomials $C_{2m+\delta}^{r+1/2}(x)/x^\delta$ are irreducible
over $\Q$. For rational $c$, the threshold is effectively computable.
In particular, $c=8$ is admissible.
\end{proposition}

\begin{proof}
For $r\ge1$, put $t=m/r\ge c$. Lemma~\ref{lem:uniform-asymptotics}
and monotonicity give
\[
 \Psi_m(S)+\log K+\log A_{\le m}
 \le m\left(1+\frac{g(c)}c\right)+o(m).
\]
The strict margin in \eqref{eq:asymptotic-condition} absorbs this
error and $\log(8\sqrt m)$ in Theorem~\ref{thm:symmetric-exact}.
For $r=0$, $K=A=1$ and the conclusion follows directly from that
criterion. All errors above are effective; for rational $c$ the strict
margin can be bounded below by rational approximation. Finally,
\[
 8(\log4-1)-g(8)>0.0097.
\]
An exact rational verification is recorded in Appendix~\ref{app:finite-derivatives}.
This also proves the remaining assertion of
Theorem~\ref{thm:symmetric}.
\end{proof}

\subsubsection{Legendre polynomials and their first derivatives}

\begin{proof}[Proof of Theorem~\ref{thm:all-degrees}]
For $P_{2m}$ and $P_{2m+1}/x$, specialize
Theorem~\ref{thm:symmetric-exact} to $r=0$. Then $K=A=1$,
and Lemma~\ref{lem:lcm} gives the sufficient inequality
\[
 m\log(4/3)>\sqrt{4m+3}\log3+\log(8\sqrt m).
\]
This holds for $m\ge100$: after division by $m$, the right side
decreases, and at $m=100$ it is less than
$0.275<2/7<\log(4/3)$.

It remains to treat $N=2m+\delta$ with $2\le N\le199$.
By Holt and Wahab, $L_N$ is irreducible if $p-4\le N\le p+3$
for some odd prime $p$ \cite[Corollary~3.4(b)]{CHL}. The only indices
in this range outside these intervals are
\[
 117,118,119,120,121,122,143,144,185,186.
\]
For each of them $m\ge4$. Put
$Z=\prod_{p\le m}p^{\lfloor\log_p(2N+1)\rfloor}$, with the product over
primes. Since $\log Z=\Psi_m(2N+1)$, condition~\eqref{eq:symmetric-exact}
is just $16^m>64m Z^2$. Direct multiplication gives this inequality in
all ten cases. This proves the first two assertions.

For the first derivatives, take $r=1$ and write $n=2m+\delta$.
The parameters in \eqref{eq:symmetric-constants} become
\[
 K=m+1,\qquad A=\frac{n+1}{n+2},\qquad S=2n+3.
\]
Since $A_{\le m}\le n+1\le2m+2$, it suffices that
\[
 m\log(4/3)>
 \sqrt{4m+5}\log3+\log(8\sqrt m)+\log\bigl(2(m+1)^2\bigr).
\]
After division by $m$, the right side decreases. At $m=200$
it is less than $(29\cdot1.1+5+12)/200=0.2445<2/7$,
so the inequality holds for every $m\ge200$.
Lemma~\ref{lem:first-derivative-finite} treats $1\le m\le199$
by exact integer comparisons and reduction modulo primes.
This proves both first-derivative assertions in every degree.
\end{proof}
\section{The Kaneko--Zagier families}\label{sec:kz-irreducibility}
\label{sec:resultant-estimates}

We turn to irreducibility of the four Kaneko--Zagier families.
The comparison of Section~\ref{sec:comparison} requires two further
estimates here: one for the common divisor introduced by clearing
denominators, and one for the clustering of roots at $2$ and $3$.
Both contribute to the leading constants.  We calculate them below,
then combine the bounds to obtain a positive gap uniform in the
possible factor degree.  Appendix~\ref{sec:explicit-cutoff} makes
the comparison explicit and completes the finite range.

\subsection{Primitive normalization}
\label{sec:irred-prelim}

In this section and Appendix~\ref{sec:explicit-cutoff}, unless
stated otherwise, $n\ge1$ and we fix
\[
 \alpha=\lambda/3,\quad \beta=\mu/2,\quad
 \lambda,\mu\in\{-1,1\},\quad c=\alpha+\beta,
 \qquad y_n(t)=P_n^{(\alpha,\beta)}(1-2t).
\]
All constants implicit in asymptotic notation are uniform in these
four choices. The functional $I$, moments $\mu_j$, and constant
$C_n$ are those of \eqref{eq:jacobi-moments} and
\eqref{eq:jacobi-wronskian}, evaluated at these parameters. Put
\begin{equation}\label{eq:primitive-normalization}
 \kappa_n=2^{n+v_2(n!)}3^{n+v_3(n!)},\qquad
 g_n=\operatorname{cont}(\kappa_n y_n),\qquad
 H_n=\frac{\kappa_n}{g_n}y_n.
\end{equation}
Here $\operatorname{cont}$ is the positive gcd of the integral coefficients.

The factor $\kappa_n$ clears the denominators at $2$ and $3$.
Removing the content $g_n$ is necessary before applying Gauss's
lemma and the local resultant bounds.  It also changes the real
derivative estimate by the factor $g_n^{-1}$.  We shall therefore
estimate $g_n$ as carefully as the local contributions to the
resultant.

\begin{lemma}\label{lem:primitive-kz}
The polynomial $H_n$ is primitive and integral. Its leading coefficient
is prime to $6$, and $(g_n,6)=1$. Moreover,
\[
 \log\kappa_n=
  \left(2\log2+\frac32\log3\right)n+O(\log n).
\]
\end{lemma}
\begin{proof}
For $p\ge5$, both generalized binomial coefficients in
\eqref{eq:irred-coefficients} are $p$-integral. At $2$, the first
is integral and the second has valuation $-j-v_2(j!)$.
At $3$, the coefficient has valuation
\[
 -n-v_3(j!)-v_3((n-j)!)
 =-n-v_3(n!)+v_3\binom nj.
\]
The minimum at either prime is attained by the leading coefficient.
This proves the integral and unit assertions. The last formula follows
from $v_p(n!)=(n-s_p(n))/(p-1)$, where $s_p(n)$ is the sum of the
base-$p$ digits of $n$.
\end{proof}

\subsection{The resultant away from \texorpdfstring{$2$ and $3$}{2 and 3}}
\label{sec:good-resultant}

Suppose that
\[
 H_n=A B,\qquad A,B\in\Z[t]\text{ primitive},\qquad
 k=\deg A\le e=\deg B,
\]
and put
\begin{equation}\label{eq:prime-to-six-resultant}
 R=|\operatorname{Res}(A,B)|>0,
 \qquad R^{(6)}=\frac{R}{2^{v_2(R)}3^{v_3(R)}}.
\end{equation}
The main estimates in this section are
\begin{equation}\label{eq:two-arithmetic-asymptotics}
 \log g_n=(\mathfrak c+o(1))n,
 \qquad \log R^{(6)}\le(\mathfrak u+o(1))kn,
\end{equation}
where
\begin{equation}\label{eq:arithmetic-constants}
 \mathfrak c=\frac{\pi}{2\sqrt3},
 \qquad \mathfrak u=\frac{43}{20}-\frac{\pi}{\sqrt3}.
\end{equation}

\subsubsection{A valuation calculation}

The content and the moment denominators come from closely related
products.  Estimating them separately loses the cancellation needed
in the resultant comparison.  We instead count divisibility by
each prime power $q$ at the same time.  The indicator $m_q$ records
the common divisibility forced on all coefficients at that level;
$e_q$ records the contribution that remains in the local bound after
the primitive normalization.  The crucial fact will be
$m_q+e_q\le1$.

Fix $p\ge5$ and $q=p^h$. Let $a_q,b_q$ be the least residues of
$\alpha,\beta$ modulo $q$, and put
\[
 c_q=(a_q+b_q)\bmod q,\qquad d_q=\lfloor(a_q+b_q)/q\rfloor,
 \qquad s=n\bmod q.
\]
Thus $1\le a_q,b_q\le q-1$ and $d_q\in\{0,1\}$. Define
\[
 E_a=\mathbf1_{s+a_q\ge q},\quad
 E_b=\mathbf1_{s+b_q\ge q},\quad
 E_c=\mathbf1_{s+c_q\ge q},
\]
\[
 m_q=\mathbf1_{E_a=E_b=1,\ d_q+E_c=1},\quad
 \eta_q=\mathbf1_{c_q=q-1},\quad
 b_q^*=\begin{cases}0,&c_q=q-1,\\1-d_q,&c_q\ne q-1,\end{cases}
\]
and
\begin{equation}\label{eq:local-budget-definition}
 e_q=\eta_q+E_a+E_b-E_c+b_q^*-2m_q.
\end{equation}
The notation $b_q^*$ is a numerical bound for moments, distinct
from the residue $b_q$.

\begin{lemma}\label{lem:carry-budget}
For every $q=p^h$, $p\ge5$, one has $e_q\in\{0,1\}$ and
$m_q+e_q\le1$.
Writing $M=12n+5$ and
\[
 E_p(n)=\sum_{\substack{h\ge1\\p^h\le M}}e_{p^h},
 \qquad m_p=v_p(g_n),
\]
one has
\begin{equation}\label{eq:moment-content-budget}
 v_p(C_n)-2m_p-\min_{0\le j\le2n-1}v_p(\mu_j)
 \le E_p(n)\le\lfloor\log_p M\rfloor.
\end{equation}
If $p>\sqrt M$, then $m_p=m_p^{\rm lev}$, where
$m_p^{\rm lev}$ means $m_q$ at the single level $q=p$.
\end{lemma}
\begin{proof}
For a rational $p$-adic integer $z$, the contribution at level
$q=p^h$ to $v_p\binom zj$ is
\[
 \mathbf1_{j\bmod q>z\bmod q}.
\]
This follows by counting multiples of $q$ in the numerator and
denominator. Hence the contribution to the coefficient of $t^j$
in \eqref{eq:irred-coefficients}, with $u=j\bmod q$, is
\begin{equation}\label{eq:coefficient-level-indicators}
 \mathbf1_{(s-u\bmod q)>(s+a_q\bmod q)}
 +\mathbf1_{u+(s+c_q\bmod q)\ge q}.
\end{equation}
Its minimum over residues $u$ is $m_q$. To see this, if $E_a=0$
the value at $u=0$ is zero. If $E_a=1$, the first indicator is
zero on $[q-a_q,s]$, and the second is zero on
$[0,q-(s+c_q\bmod q)-1]$. These intervals are disjoint exactly
when $E_b=1$ and $d_q+E_c=1$. The same minimum is attained for
$0\le j\le n$: if $n<q$, the relevant interval $[q-a_q,s]$
already lies in $[0,n]$. Therefore
\begin{equation}\label{eq:content-level-lower}
 m_p\ge\sum_{h\ge1}m_{p^h}.
\end{equation}

Counting the four products in $C_n$ shows that its contribution
at level $q$ is $\eta_q+E_a+E_b-E_c$. For the moments
\eqref{eq:jacobi-moments}, complete periods in the numerator and
denominator cancel. If $c_q\ne q-1$, the first denominator factor
divisible by $q$ precedes the first numerator factor exactly when
$d_q=0$. If $c_q=q-1$, it does not precede it. Thus every moment
has contribution at least $-b_q^*$ at that level.

It remains to check \eqref{eq:local-budget-definition}. If
$c_q\ne q-1$ and $d_q=0$, then
$E_c\ge\max(E_a,E_b)$; the case $E_a=E_b=1$ has $m_q=1$.
If $d_q=1$, then $E_c\le\min(E_a,E_b)$; now $m_q=1$ exactly
when $E_a=E_b=1$ and $E_c=0$. These cases give $e_q=0$ or $1$.
If $c_q=q-1$, then $d_q=0$. For $s=0$ the value is one,
and for $s\ge1$ it is
$E_a+E_b-2\mathbf1_{E_a=E_b=1}$, again zero or one.

All numerator factors relevant to the coefficients, $C_n$, and
the moments through index $2n-1$ are at most $M$, apart from
powers of $2$ and $3$. Summing the level inequalities proves
\eqref{eq:moment-content-budget}. For $p>\sqrt M$ there is only
one level, so the minimum in
\eqref{eq:coefficient-level-indicators} is the actual content
valuation.
\end{proof}

The moment bound alone allows primes as large as $12n+5$.
For primes exceeding the degree, the differential equation gives
more information.  A repeated reduced root can occur only at an
endpoint or at infinity, and the endpoint and leading-coefficient
valuations are too small for two primitive factors to meet there.
For $n\ge11$, this removes every prime $p>n$ with $p\ge5$
from the resultant.

\begin{lemma}
\label{lem:projective-cluster-exclusion}
Let $H\in\mathbb Z_p[t]$ be primitive of degree $n<p$ and satisfy
\[
 t(1-t)H''+\tau(t)H'+\Lambda H=0,
 \qquad \tau\in\mathbb Z_p[t],\quad \Lambda\in\mathbb Z_p.
\]
If
\[
 v_p(H(0))\le1,\qquad v_p(H(1))\le1,\qquad
 v_p(\operatorname{lc}H)\le1,
\]
then every primitive factorization $H=AB$ in $\mathbb Z_p[t]$
satisfies $v_p(\operatorname{Res}(A,B))=0$.
\end{lemma}

\begin{proof}
The nonzero reduction $\overline H$ has no multiple root outside
$0,1$: a root $\xi$ of multiplicity $2\le d\le\deg\overline H<p$
would give a nonzero lowest coefficient
$\xi(1-\xi)d(d-1)$ times the remaining nonvanishing factor in the
reduced differential equation.

Homogenize $A,B$ using their original degrees. A common reduced
finite root must therefore be $a=0$ or $1$, forcing
$v_p(A(a))+v_p(B(a))\ge2$, a contradiction. A common reduced root
at infinity would similarly force both leading coefficients to be
divisible by $p$. Thus the reduced homogeneous polynomials have
no common projective root. Their resultant is nonzero, and reduction
commutes with the resultant.
\end{proof}

\begin{lemma}\label{lem:good-prime-resultant}
For $p\ge5$,
\begin{equation}\label{eq:good-prime-bound}
 v_p(R)\le kE_p(n).
\end{equation}
If $n\ge11$ and $p>n$, then $v_p(R)=0$.
\end{lemma}
\begin{proof}
Fix $p\ge5$. Orthogonality and \eqref{eq:jacobi-wronskian} allow
Lemma~\ref{lem:projective-moment-budget} to be applied to $H_n$
with moments \eqref{eq:jacobi-moments} and
$C=(\kappa_n/g_n)^2C_n$. Since $\kappa_n$ is a $p$-adic unit,
Lemma~\ref{lem:carry-budget} gives $\mathcal E_p\le E_p(n)$,
proving \eqref{eq:good-prime-bound}.

Now suppose $n\ge11$ and $p>n$. Each of
\[
 y_n(0)=\frac{(\alpha+1)_n}{n!},\qquad
 y_n(1)=(-1)^n\frac{(\beta+1)_n}{n!},\qquad
 \operatorname{lc}(y_n)=(-1)^n\frac{(n+c+1)_n}{n!}
\]
has valuation at most one: $n!$ is a unit, and each cleared
numerator product has length $n<p$, step prime to $p$, and factors
at most $12n+5<(n+1)^2\le p^2$. Multiplication by $\kappa_n$
preserves these valuations, and removing content cannot increase them.
Thus the corresponding valuations for $H_n$ are at most one.
Equation~\eqref{eq:unit-interval-jacobi-ode} has $p$-integral
coefficients, so Lemma~\ref{lem:projective-cluster-exclusion}
gives $v_p(R)=0$.
\end{proof}

\subsubsection{The two asymptotic constants}

For primes exceeding $\sqrt{12n+5}$ there is only one valuation
level.  The preceding indicators then depend on the fractional
part of $n/p$ and on the residue class of $p$ modulo $6$.
Their supports are short intervals.  Summing the lengths of these
intervals gives the two constants below; smaller primes and higher
prime powers contribute only an error of smaller order than $n$.

\begin{proposition}\label{prop:arithmetic-asymptotics}
The estimates \eqref{eq:two-arithmetic-asymptotics} hold uniformly
in the four families and in the rational factorization.
\end{proposition}
\begin{proof}
For $p>\sqrt{12n+5}$ there is only one valuation level.
In the class $p\equiv r\pmod6$, where $r=2\lambda+3\mu$,
one has $a_p/p=1/3+O(1/p)$; in the other reduced class,
$a_p/p=2/3+O(1/p)$. In both classes
$b_p/p=1/2+O(1/p)$.
The indicator calculation in Lemma~\ref{lem:carry-budget} gives
the following intervals, up to endpoint changes of size $O(1/p)$:
\[
\begin{array}{c|c|c}
 a_p/p&\{n/p\}\text{ for }e_p=1&\{n/p\}\text{ for }m_p=1\\ \hline
 1/3&[0,1/6)\cup[1/2,2/3)&[2/3,1)\\
 2/3&[1/3,1/2)\cup[5/6,1)&[1/2,5/6).
\end{array}
\]
Small primes and higher prime powers contribute $O(\sqrt n\log n)$.
For the content this follows also from
$v_p(g_n)\le v_p(\operatorname{lc}(y_n))
\le\lfloor\log_p(12n+5)\rfloor$.

Apply the prime number theorem in the two classes modulo $6$
to the intervals in the table; the fixed-modulus consequence of \cite[Corollary~1.7]{BMOR}
suffices. For each fixed bound on
$\ell=\lfloor n/p\rfloor$ only finitely many intervals occur.
The remaining primes contribute $O(n/\ell)$, so the resulting
convergent sums may be obtained by first fixing this bound and
then letting it tend to infinity. Endpoint perturbations contribute
$o(n)$; after removing $p\le\sqrt{12n+5}$ there are only
$O(\sqrt n)$ endpoints, each shifted by $O(1)$.

Lemma~\ref{lem:good-prime-resultant} permits restriction to $p\le n$.
The upper bound for the resultant is therefore
\begin{align*}
 \frac1n\sum_{5\le p\le n}E_p(n)\log p\longrightarrow
 \frac12\sum_{\ell\ge1}\bigg(&\frac1\ell-\frac1{\ell+1/6}
 +\frac1{\ell+1/3}-\frac1{\ell+2/3}\\
 &+\frac1{\ell+5/6}-\frac1{\ell+1}\bigg)
 =\frac{43}{20}-\frac{\pi}{\sqrt3}.
\end{align*}
For the content all primes are included, so the corresponding sum
starts at zero:
\[
 \frac{\log g_n}{n}\longrightarrow
 \frac12\sum_{\ell\ge0}\left(
 \frac1{\ell+2/3}-\frac1{\ell+1}
 +\frac1{\ell+1/2}-\frac1{\ell+5/6}\right)
 =\frac{\pi}{2\sqrt3}.
\]
These evaluations follow from the partial-fraction expansion of
$\pi\cot(\pi z)$, or from the reflection and duplication formulas
for the logarithmic derivative of $\Gamma$
\cite[Eqs.~5.5.4 and 5.5.8]{DLMF}.
\end{proof}

\subsection{The primes \texorpdfstring{$2$ and $3$}{2 and 3}}\label{sec:bad-resultant}

The monic endpoint expansions place every root in the ball
$v_2(t)\ge1$ at $2$, and in the ball $v_3(1-t)\ge1$ at $3$.
The endpoint values of a factor measure the total valuations of
its roots in these balls.  If two roots have different endpoint
valuations, their difference has the smaller valuation.  The only
additional contribution comes from roots with equal valuations.
We bound it by expanding the differential equation at one root.

To sum the resulting errors, we use the Newton polygon.  For a
polynomial $f(X)=\sum_j c_jX^j$ over $\Q_p$ with nonzero constant
term, this is the lower convex hull of the points $(j,v_p(c_j))$,
omitting zero coefficients.  A side of slope $-r$ and horizontal
length $d$ corresponds to exactly $d$ roots in $\overline{\Q}_p$
of valuation $r$, counted with multiplicity.  Thus inequalities for coefficient valuations
can bound the total number and size of exceptional root valuations.

For the factorization $H_n=AB$, put
\[
 D_{2,A}=v_2(A(0)),\quad D_{2,B}=v_2(B(0)),\qquad
 D_{3,A}=v_3(A(1)),\quad D_{3,B}=v_3(B(1)).
\]
The leading coefficients of both factors are units at $2$ and $3$.

\begin{proposition}\label{prop:bad-prime-resultant}
Let $s_p(m)$ denote the sum of the base-$p$ digits of $m$, and put
$\delta=(\mu+1)/2$. Then
\begin{align*}
 v_2(R)&\le\min(eD_{2,A},kD_{2,B})
                  +\frac{ke}{2}+\frac{k}{2}s_2(n),\\
 v_3(R)&\le\min(eD_{3,A},kD_{3,B})
                  +\frac{ke}{3}+\frac{k}{6}s_3(2n+\delta).
\end{align*}
These bounds hold for every factorization; no restriction on $k$
is needed. In particular, the error beyond the displayed bilinear
terms is $O(k\log n)$.
\end{proposition}

\begin{proof}
The proof has two parts.  First, the Taylor coefficients at a root
give a uniform upper bound for the extra closeness of another root
with the same endpoint valuation.  The error in that bound measures
how far the root valuation lies below $2$ at the prime $2$, or below
$3/2$ at the prime $3$.  Second, the coefficient valuations and the
Newton polygon show that the total of these deficits is controlled
by a digit sum, hence is $O(\log n)$.

At $p=2$ use the endpoint coordinate $z=t$ and
$b=\alpha$. At $p=3$ use $z=1-t$ and $b=\beta$. In either case
the corresponding Jacobi polynomial $P(z)$ satisfies
\[
 \sigma P''+\tau P'+\Lambda P=0,
 \qquad
 \sigma=z(1-z),\quad \tau=b+1-(c+2)z,
 \quad \Lambda=n(n+c+1).
\]
Every root $z$ has $v_p(z)\ge1$. Indeed, the monic endpoint
coefficients are given by
\begin{equation}\label{eq:monic-endpoint}
 \frac{n!P_n^{(a,b)}(2X-1)}{(n+a+b+1)_n}
 =\sum_{j=0}^n(-1)^{n-j}\binom nj
   \frac{(b+j+1)_{n-j}}{(n+a+b+j+1)_{n-j}}X^j.
\end{equation}
Here $b\in\mathbb Z_p$, while every factor in the denominator
has valuation $-1$. Thus their valuations are at least $n-j$.

Fix a root $z$, put $r=v_p(z)\ge1$, and let
\[
 h=\min(1,r-1),\qquad
 u_j=z^{j-1}\frac{P^{(j)}(z)}{P'(z)}.
\]
Thus $u_0=0$ and $u_1=1$. Differentiating the differential equation
$m$ times gives
\begin{equation}\label{eq:taylor-valuation-recurrence}
 \begin{gathered}
 u_{m+2}=-a_mu_{m+1}-b_mu_m,\\
 a_m=\frac{b+1+m-(c+2+2m)z}{1-z},\qquad
 b_m=\frac{z(n-m)(n+m+c+1)}{1-z}.
 \end{gathered}
\end{equation}
Since $v_p(c+2+2m)=v_p(n+m+c+1)=-1$, we have
$v_p(b_m)\ge r-1\ge h$ and $v_p(a_m)\ge0$.
At $2$, $\alpha+1$ is even in $\mathbb Z_2$, so
$v_2(a_m)\ge h$ whenever $m$ is even.
At $3$, $\beta+1$ is congruent to $2$ if $\delta=0$ and to $0$
if $\delta=1$, so $v_3(a_m)\ge h$ whenever
$m\equiv1-\delta\pmod3$. Induction in
\eqref{eq:taylor-valuation-recurrence} therefore gives
\begin{equation}\label{eq:taylor-coefficient-gain}
 v_2(u_j)\ge\left\lfloor\frac j2\right\rfloor h,
 \qquad
 v_3(u_j)\ge\left\lfloor\frac{j+\delta}{3}\right\rfloor h.
\end{equation}
For example, in the $3$-adic induction the first term gains an
extra $h$ precisely when the target floor increases, and the
second term gains at least $h$ at every step. The initial values
satisfy the stated bounds.

For $p=2$ or $3$, $j\ge2$, and $0\le h\le1$, Legendre's formula
implies
\begin{equation}\label{eq:factorial-radius-bound}
 \frac{v_p(j!)-\lfloor j/p\rfloor h}{j-1}
 \le\frac1{p-1}-\frac hp.
\end{equation}
To verify it, write $j=pq+s$, $0\le s<p$. After rearrangement,
the required inequality is
\[
 \frac{1-s_p(j)}{p-1}+\frac{h(s-1)}p\le0.
\]
For $s=0,1$ this is immediate; for $s\ge2$ use
$s_p(j)\ge s$ and $h\le1$.
If $w\ne z$ is another root, the Taylor quotient
\[
 \frac{P(w)}{(w-z)P'(z)}
 =1+\sum_{j=2}^n\frac{u_j}{j!}
                  \left(\frac{w-z}{z}\right)^{j-1}
\]
cannot vanish when
$v_p((w-z)/z)>1/(p-1)-h/p$. Thus, writing
$\rho_p=p/(p-1)$, we obtain the convenient weaker bound
\begin{equation}\label{eq:paired-prime-separation}
 v_p(w-z)-\min(v_p(w),v_p(z))
 \le\frac1p+\frac{(\rho_p-v_p(z))_+}{p}
 \quad\text{if }v_p(w)=v_p(z).
\end{equation}
The excess is zero when the root valuations differ. Here
$(x)_+=\max(x,0)$; the bound follows from
\[
 \frac1{p-1}-\frac{\min(1,r-1)}p
 \le\frac1p+\frac{(\rho_p-r)_+}{p}.
\]

It remains to control the total endpoint deficits.
Let $c_j$ be the coefficients of the monic normalization of
$y_n(t)$. At $2$ their exact valuations are
\[
 v_2(c_j)=2(n-j)+\Psi(j),\qquad
 \Psi(j)=s_2(j)-s_2(n)+v_2\binom{n+\alpha}{n-j}.
\]
The generalized binomial coefficient is $2$-integral, so
$\Psi(j)\ge-s_2(n)$, and $\Psi(n)=0$.
If $W$ is the lower convex hull of the points $(j,\Psi(j))$,
then the Newton-polygon theorem gives
\[
 \sum_i(2-v_2(t_i))_+
 =\int_0^n(W'(x))_+\,dx
 =W(n)-\min W\le s_2(n).
\]
The middle equality follows from convexity: a convex function
first decreases and then increases.

At $3$, the monic polynomial in $z=1-t$ has coefficient valuations
\[
 \frac32(n-j)+\frac12\Phi(j),\qquad
 \Phi(j)=s_3(2j+\delta)+s_3(n-j)-s_3(2n+\delta).
\]
Here $\Phi(j)\ge-s_3(2n+\delta)$ and $\Phi(n)=0$. Applying
the same positive-variation argument gives
\[
 \sum_i\left(\frac32-v_3(1-t_i)\right)_+
 \le\frac12s_3(2n+\delta).
\]

The leading coefficients of $A,B$ are units at $2,3$, so the
valuation of their resultant is the sum of the cross-factor
root difference valuations. The sum of the minimum valuations
is at most both $eD_{p,A}$ and $kD_{p,B}$. Sum
\eqref{eq:paired-prime-separation}, using the root of $B$ for the
deficit term. The sum is at most $ke/p$ plus $k/p$ times the
total endpoint deficit among the roots of $B$, and hence among
all roots of $H_n$. The two asserted inequalities follow.
\end{proof}

\subsection{The archimedean comparison}
\label{sec:irred-proof}

The real roots lie in $(0,1)$.  Their small Vandermonde product
therefore makes the resultant large once the derivatives are bounded
below.  The leading coefficient of a factor contributes another
positive term, because its values at $0$ and $1$ are nonzero
integers.  These are the same endpoint values that occur in the
local bounds at $2$ and $3$.  Combining the bounds before optimizing
eliminates the unknown way in which the two factors share the
endpoint valuations.  Retaining the Vandermonde term makes the final
estimate uniform in the factor degree.

Write $\gamma_n=\kappa_n/g_n$, and set
\begin{align}
 D_{2,n}&=n+\sum_{j=1}^n v_2(3j+\lambda),&
 T_{2,n}&=D_{2,n}\log2/n,\label{eq:exact-two-total}\\
 D_{3,n}&=n+\sum_{j=1}^n v_3(2j+\mu),&
 T_{3,n}&=D_{3,n}\log3/n.\label{eq:exact-three-total}
\end{align}
These are the endpoint valuations of $H_n$, and
\begin{equation}\label{eq:endpoint-asymptotics}
 D_{2,n}=2n+O(\log n),\qquad
 D_{3,n}=\tfrac32n+O(\log n).
\end{equation}
Define
\begin{equation}\label{eq:arch-constants-finite}
 \begin{split}
 L_n&=\log(e^{T_{2,n}}+e^{T_{3,n}}),\qquad
 E=\tfrac12\log2+\tfrac13\log3,\\
 J_n&=\tfrac12s_2(n)\log2+\tfrac16s_3(2n+\delta)\log3,
 \qquad \delta=(\mu+1)/2.
 \end{split}
\end{equation}

\begin{proposition}\label{prop:arch-resultant}
For every factorization $H_n=AB$ with $1\le k=\deg A\le n/2$,
\begin{equation}\label{eq:exact-arch-resultant}
 \begin{split}
 \frac{\log R^{(6)}}k\ge{}&
 \log(2\gamma_n\sqrt{C_n})-J_n
 +(n-2k)L_n\\
 &-(n-k)(T_{2,n}+T_{3,n}+E)
 +\frac{2(k-1)^2}{k}\log2-\log k.
 \end{split}
\end{equation}
Consequently, uniformly in all such factorizations,
\begin{equation}\label{eq:arch-resultant-lower}
 \log R^{(6)}\ge(\mathfrak l-o(1))kn,
 \quad
 \mathfrak l=\log(4+3\sqrt3)-\tfrac12\log2-\tfrac13\log3
                         -\frac{\pi}{2\sqrt3}.
\end{equation}
\end{proposition}
\begin{proof}
Let $a=|\operatorname{lc}(A)|$, $e=n-k$, and let
$t_1,\ldots,t_k$ be the roots of $A$. Differentiating $H_n=AB$
at these roots gives the exact identity
\[
 R=a^{e-k}\frac{\prod_{i=1}^k|H_n'(t_i)|}
                    {\prod_{i<j}|t_i-t_j|^2}.
\]
Lemmas~\ref{lem:jacobi-wronskian} and
\ref{lem:interval-vandermonde} therefore give
\[
 R\ge a^{n-2k}(2\gamma_n\sqrt{C_n})^k
                   k^{-k}2^{2(k-1)^2}.
\]
By H\"older's inequality,
\[
 \left(\prod_i t_i\right)^{1/k}
 +\left(\prod_i(1-t_i)\right)^{1/k}\le1.
\]
Hence, with $U=D_{2,A}\log2/k$ and $V=D_{3,A}\log3/k$,
\begin{equation}\label{eq:leading-endpoint-bound}
 a^{1/k}\ge|A(0)|^{1/k}+|A(1)|^{1/k}\ge e^U+e^V.
\end{equation}

Put $x=k/n$. Since $D_{p,n}=D_{p,A}+D_{p,B}$,
\[
 \frac{\min(eD_{2,A},kD_{2,B})\log2}{kn}
   =\min(U,T_{2,n})-xU,
\]
and the analogous identity holds at $3$. This exact expression
allows either factor to contain an arbitrary share of the endpoint
valuations. For $0\le x\le1/2$, the function
\[
 (1-2x)\log(e^U+e^V)+x(U+V)
       -\min(U,T_2)-\min(V,T_3)
\]
has its minimum at $(U,V)=(T_2,T_3)$. Indeed, with $V$ fixed,
its derivative is nonpositive below $T_2$ and nonnegative above:
\[
 (1-2x)\frac{e^U}{e^U+e^V}+x-1\le0,
 \qquad
 (1-2x)\frac{e^U}{e^U+e^V}+x\ge0.
\]
The same argument applies to $V$. Subtracting the two local
resultant bounds of Proposition~\ref{prop:bad-prime-resultant}
and applying this minimum proves \eqref{eq:exact-arch-resultant}.

To obtain the uniform asymptotic form, divide by $n$ and use
\[
 \frac{2(k-1)^2}{kn}\log2
 =2x\log2-\frac{4\log2}{n}+\frac{2\log2}{kn}.
\]
Since $J_n=O(\log n)$ and $\log k\le\log n$, the result is
\begin{align*}
 \frac{\log R^{(6)}}{kn}\ge{}&
 \frac{\log\gamma_n}{n}+L_n-T_{2,n}-T_{3,n}-E\\
 &+x\bigl(T_{2,n}+T_{3,n}+E-2L_n+2\log2\bigr)-o(1),
\end{align*}
with an error independent of $k$. The coefficient of $x$ is
eventually positive. In fact, $L_n\le\max(T_{2,n},T_{3,n})+\log2$
shows that it is at least $E-|T_{2,n}-T_{3,n}|$, whose limit is
\[
 \frac52\log2-\frac76\log3
 =\frac16\log\frac{2^{15}}{3^7}>0.
\]
We may discard this term. Now use
Lemma~\ref{lem:primitive-kz}, Proposition~\ref{prop:arithmetic-asymptotics},
\eqref{eq:root-derivative-lower}, and \eqref{eq:endpoint-asymptotics}.
The remaining expression tends to $\mathfrak l$, proving
\eqref{eq:arch-resultant-lower}.
\end{proof}

The two uniform bounds already imply eventual irreducibility, since
\begin{equation}\label{eq:positive-gap}
 \mathfrak l-\mathfrak u
 =\log(4+3\sqrt3)-\tfrac12\log2-\tfrac13\log3
                    -\frac{43}{20}+\frac{\pi}{2\sqrt3}
 >\frac14.
\end{equation}
For reference, $\mathfrak l=0.5991078130\ldots$,
$\mathfrak u=0.3362006358\ldots$, and their difference is
$0.2629071772\ldots$.  For an elementary verification of the
rational bound, use
\[
 \frac{1732}{1000}<\sqrt3<\frac{1733}{1000},\qquad
 \pi>\frac{157}{50},\qquad
 \log3>\frac{549}{500},\qquad \log2<\frac7{10}.
\]
Since $4+3\sqrt3>9\cdot2299/2250$ and
$\log(1+t)>t/(1+t)$ for $t>0$, these bounds give
\[
 \mathfrak l-\mathfrak u>
 \frac53\frac{549}{500}-\frac7{20}-\frac{43}{20}
       +\frac{49}{2299}+\frac{1570}{1733}>\frac14.
\]
The logarithmic bounds follow, for example, from the positive-term
series used in Appendix~\ref{sec:finite-verification}.
Thus a proper factorization is impossible once $n$ is sufficiently
large, with no restriction on the degree of a possible factor.
To obtain an all-degree theorem, we retain the finite terms in the
same comparison.

\subsection{An exact finite criterion}\label{sec:finite-criterion}

The resultant estimates give a finite criterion
that does not require a factorization of the Jacobi polynomial.
Use the quantities in \eqref{eq:exact-two-total}--\eqref{eq:arch-constants-finite}
and the residue bounds of Lemma~\ref{lem:carry-budget}, and put
\[
 S_n=\sum_{\substack{p\text{ prime}\\5\le p\le n}}E_p(n)\log p.
\]
For $z\ge1$, define
\begin{equation}\label{eq:finite-criterion-function}
 \begin{split}
 \mathcal H_n(z)={}&\log(2\gamma_n\sqrt{C_n})-J_n
  +(n-2z)L_n-(n-z)(T_{2,n}+T_{3,n}+E)\\
 &+\frac{2(z-1)^2}{z}\log2-\log z.
 \end{split}
\end{equation}

\begin{proposition}\label{prop:finite-criterion}
Let $n\ge11$. If
\[
 \min_{1\le k\le\lfloor n/2\rfloor}\mathcal H_n(k)>S_n,
\]
then $y_n$ is irreducible over $\Q$.
More generally, the minimum may be restricted to any set known to
contain all possible degrees of a nonconstant factor of degree at
most $n/2$.
\end{proposition}
\begin{proof}
A nontrivial factorization would give
$S_n\ge\log R^{(6)}/k\ge\mathcal H_n(k)$ by
Lemma~\ref{lem:good-prime-resultant} and
\eqref{eq:exact-arch-resultant}, a contradiction.
\end{proof}

The minimization is elementary. Set
\[
 a_n=T_{2,n}+T_{3,n}+E-2L_n+2\log2.
\]
The part of $\mathcal H_n(z)$ that depends on $z$ is
$a_nz+2\log2/z-\log z$, which is strictly convex for $z>0$.
If $a_n\le0$ it is decreasing. If $a_n>0$, its minimum is at
\[
 q_n=\frac{1+\sqrt{1+8a_n\log2}}{2a_n}.
\]
Thus on an interval of integer degrees it suffices to check the
integers adjacent to the clipped value of $q_n$.

The criterion depends only on endpoint valuations, content, and
finite sums over primes.  Appendix~\ref{sec:finite-verification}
uses a slightly weaker form in which the possible factor degree
has been eliminated and every comparison reduces to an inequality
between linear forms in logarithms.

\subsection{Completion in every degree}

Appendix~\ref{sec:explicit-cutoff} proves irreducibility for
$n\ge250\,000$.  The degrees $1\le n\le8332$ are covered by the
published computation of Cullinan and Hajir
\cite[Remark~4.6]{CHK}; the conversion from their index to the
degree is given in Appendix~\ref{sec:finite-verification}.
For $8333\le n\le250\,000$, two explicit inequalities from the
resultant comparison are verified by integer arithmetic.  Together
these ranges give irreducibility in every degree.

\begin{proof}[Irreducibility in Theorem~\ref{thm:kz-all}]
Proposition~\ref{prop:explicit-kz-cutoff} treats $n\ge250\,000$,
and Proposition~\ref{prop:finite-kz} treats $1\le n\le250\,000$.
The Galois assertion will follow from Theorem~\ref{thm:kz}.
\end{proof}
\section{Galois groups and splitting fields}
\label{sec:galois-groups}

Irreducibility gives transitivity on the roots. The remaining task
is to find enough permutations in the Galois group to determine it.
We obtain these permutations from ramification at suitable primes.
The first argument is common to the Legendre and Kaneko--Zagier
families and extends to fixed rational parameters. For unequal
integral parameters, a different local factorization gives a
prime-length cycle. We also keep track of which roots each inertia
element fixes. This support information will prove the
splitting-field disjointness statements.

\subsection{Inertia cycles and permutation groups}\label{sec:cycles}

We call a permutation a \emph{single $d$-cycle} if it cyclically
permutes $d$ points and fixes every other point; for $d=1$ this
means the identity. All valuations $v_p$ are normalized by
$v_p(p)=1$.

\subsubsection{A local criterion}

Under the hypotheses of the next lemma, the Newton polygon has
a final edge of slope $1/d$: exactly $d$ roots have valuation
$-1/d$, and all other roots have valuation zero. Passing to the
reciprocal polynomial isolates the first set in an Eisenstein
factor.

In a finite Galois extension of $\Q_p$, the inertia subgroup
consists of the automorphisms acting trivially on the residue
field. We want its action to be a single cycle on the isolated
roots. The degree bound $d<p$ makes this inertia action tame,
and squarefree reduction puts all remaining roots in an
unramified extension. Thus inertia fixes those roots. The lemma
also records the effect of a quadratic substitution.

\begin{lemma}\label{lem:local-cycle}
Let $f(X)=\sum_{j=0}^{N}a_jX^j\in\Z_p[X]$ have degree $N$.
Assume that $a_0$ is a unit, $v_p(a_N)=1$, and the reduction
$\overline f$ is squarefree of degree $N-d$, where $1\le d<p$.
Then an inertia subgroup of the Galois group of $f$ over $\Q_p$
contains a single $d$-cycle. The other $N-d$ roots lie in an
unramified extension of $\Q_p$.

If $p$ is odd, $c\in\Z_p^\times$, and $2d<p$, inertia for
$f(cX^2)$ contains a single $2d$-cycle. Its image on the opposite
root pairs is a single $d$-cycle, and its $d$th power interchanges
exactly those $d$ pairs.
\end{lemma}

\begin{proof}
The normalized reciprocal
\[
 R(Y)=\frac{Y^Nf(1/Y)}{a_0}
\]
is monic and integral. Its reduction is
$Y^d\overline B(Y)$, where $\overline B(0)\ne0$ and
$\overline B$ is squarefree. Hensel factorization gives
$R=AB$ with $A,B$ monic, $\deg A=d$, and
$A\equiv Y^d\pmod p$. Since $B(0)$ is a unit and
$v_p(R(0))=1$, the polynomial $A$ is Eisenstein. The roots of
$B$ lift over unramified extensions.

The polynomial $A$ remains Eisenstein over every finite unramified
extension. In particular, it is irreducible over the fixed field
of inertia in the local splitting field. Thus inertia acts
transitively on its $d$ roots. Wild inertia is a $p$-group. It fixes the roots of $B$, and its faithful
action on the remaining roots embeds it in $S_d$. Since $d<p$,
this wild inertia group is trivial. The finite inertia group is
therefore cyclic, and a generator is a single $d$-cycle. Inverting
the roots proves the first assertion. The argument also covers
$B=1$.

The reduction of $f(cX^2)$ is squarefree: its constant term is a
unit, $\overline f$ is squarefree, and $p\ne2$. Its degree drops
by $2d$ and its leading coefficient has valuation one. The first
assertion gives a single $2d$-cycle $\tau$. Negation commutes
with every field automorphism and preserves the support of $\tau$.
The centralizer of a full cycle on its support is cyclic. Hence
negation, a fixed-point-free involution there, equals $\tau^d$.
This proves both assertions about the opposite pairs.
\end{proof}

\subsubsection{Jacobi coefficients and reduction}

The normalization \eqref{eq:jacobi-definition} satisfies
\begin{align}
 P_n^{(\alpha,\beta)}(2X+1)
 &=\sum_{j=0}^{n}
   \binom{n+\alpha}{n-j}
   \binom{n+\alpha+\beta+j}{j}X^j,\label{eq:jacobi-shift}\\
 P_n^{(\alpha,\beta)}(1)&=\binom{n+\alpha}{n},\qquad
 P_n^{(\alpha,\beta)}(-1)=(-1)^n\binom{n+\beta}{n},
 \label{eq:jacobi-endpoints}
\end{align}
and its differential equation is
\begin{equation}\label{eq:jacobi-ode}
 (1-X^2)y''+[\beta-\alpha-(\alpha+\beta+2)X]y'
             +n(n+\alpha+\beta+1)y=0.
\end{equation}
These identities follow from the binomial expansion and hold as
polynomial identities in the parameters; see also
\cite[Eqs.~18.5.7--18.5.8 and Table~18.8.1]{DLMF}.

Write $\alpha=a/u$, $\beta=b/v$, where $a,b\in\Z$ and $u,v$ are
positive integers, and choose a positive integer $D$ such that
\begin{equation}\label{eq:D-r}
 r=D(\alpha+\beta)\in\Z.
\end{equation}
Equation~\eqref{eq:jacobi-shift} gives the integral normalization
\begin{align}
 K_n^{\alpha,\beta}(X)
 &:=u^n n!P_n^{(\alpha,\beta)}\left(1+\frac{2D}{u}X\right)\notag\\
 &=\sum_{j=0}^{n}\binom nj
       \prod_{k=j+1}^{n}(uk+a)
       \prod_{k=1}^{j}(Dn+r+Dk)X^j.\label{eq:integral-Jacobi}
\end{align}

This normalization separates the factors involving $\alpha$ from
those involving $\alpha+\beta$. We choose $p$ so that the first
product consists of units, while the second product encounters
$p$ exactly once. From that coefficient onward, every coefficient
has valuation one. The degree of the reduction therefore gives
the length of the inertia cycle. The Jacobi differential equation
will show that this reduction is squarefree, completing the local
criterion.

\begin{proposition}\label{prop:jacobi-cycle}
Let $p$ be an odd prime satisfying
\begin{gather}
 p>n,\qquad p\nmid uvD,\qquad
 0<|uk+a|<p,\quad0<|vk+b|<p\quad(1\le k\le n),
 \label{eq:local-small-factors}\\
 0<D(n+1)+r,\qquad 2Dn+r<2p,\qquad p\equiv r\pmod D.
 \label{eq:local-large-factors}
\end{gather}
Suppose
\[
 1\le k_0:=\frac{p-r}{D}-n\le n,
 \qquad d:=n-k_0+1.
\]
Then an inertia subgroup at $p$ of
$\Gal(P_n^{(\alpha,\beta)}/\Q)$ contains a single $d$-cycle.
If $c\in\Q^\times$, $v_p(c)=0$, and $2d<p$, an inertia
subgroup at $p$ of
$\Gal(P_n^{(\alpha,\beta)}(1+cX^2)/\Q)$ contains a single
$2d$-cycle, with the action on opposite pairs described in
Lemma~\ref{lem:local-cycle}.
\end{proposition}

\begin{proof}
In \eqref{eq:integral-Jacobi}, the binomial coefficients and the
first product are $p$-adic units. Every factor in the second product
lies in $(0,2p)$. Exactly one is divisible by $p$, namely the factor
with index $k_0$, and it equals $p$. The coefficient valuations are
therefore zero for $j<k_0$ and one for $j\ge k_0$. In particular,
the constant coefficient is a unit, the leading coefficient has
valuation one, and the reduction has degree $k_0-1=n-d$.

It remains to prove that this reduction is squarefree. The Jacobi
polynomial and \eqref{eq:jacobi-ode} are integral at $p$, and both
endpoint values in \eqref{eq:jacobi-endpoints} are units by
\eqref{eq:local-small-factors}. Its reduction has degree at most
$n<p$. If a root $z\ne\pm1$ had multiplicity $h\ge2$, write
$\overline y=(X-z)^h g(X)$ with $g(z)\ne0$. Since $2\le h<p$,
the coefficient of $(X-z)^{h-2}$ in the reduced differential
equation would be
\[
 (1-z^2)h(h-1)g(z)\ne0,
\]
a contradiction. Thus the reduction is squarefree. The scalar and
affine transformation in \eqref{eq:integral-Jacobi} are units at
$p$, so the reduction of $K_n^{\alpha,\beta}$ is squarefree as well.

Lemma~\ref{lem:local-cycle} now gives a single $d$-cycle. A local
Galois group identifies with a decomposition subgroup of the global
group, so this is also a global Galois element. For the quadratic
assertion, use
\[
 P_n^{(\alpha,\beta)}(1+cX^2)
   =(u^n n!)^{-1}K_n^{\alpha,\beta}
                          \left(\frac{uc}{2D}X^2\right)
\]
and apply the second part of that lemma.
\end{proof}

\subsubsection{Blocks and opposite pairs}

A block system for a transitive group is a partition into equally
sized sets permuted by the group. It is nontrivial when the block
size is strictly between $1$ and the degree; a transitive group is
primitive if it has no such system.

A cycle moving more than half the roots places a simple
restriction on a nontrivial block system: every block meeting the
support of the cycle must lie entirely in that support.
The block size must therefore divide the cycle length.

\begin{lemma}\label{lem:blocks}
Let a transitive group of degree $N$ preserve a nontrivial block
system with blocks of size $b$. If it contains a single $d$-cycle
with $d>N/2$, then $b\mid d$.
\end{lemma}

\begin{proof}
We have $b\le N/2<d$. A block containing a fixed point of the
cycle is fixed setwise, because it intersects its image. If it
also contained a moved point, it would contain the entire orbit
of that point, which has size $d>b$. Thus the support of the
cycle is a union of whole blocks, and $b\mid d$.
\end{proof}

Once primitivity is known, a single cycle with at least three
fixed points forces the alternating group. We use Jones's theorem
in this form; the cycle length need not be prime.

\begin{theorem}[Jones {\cite[Corollary~1.3]{Jones}}]\label{thm:jones}
A primitive permutation group of degree $N$ containing a
nontrivial single cycle with at least three fixed points contains
$A_N$.
\end{theorem}

This form of Jones's theorem uses the classification of finite
simple groups. The prime-length case used later is Jordan's
classical theorem.

For roots occurring in opposite pairs, the largest possible group has two
parts: permutations of the pairs and interchanges within each
pair. Even after the group on pairs is known to be $S_N$, the
interchanges may satisfy a relation. A single $2d$-cycle with
odd $d<N$ removes this obstruction. Its $d$th power flips a
nonempty proper set of an odd number of pairs, and permutations
of the pairs generate every sign vector from this one vector.

\begin{lemma}\label{lem:sign-kernel}
Let $\Gamma\le\F_2^N\rtimes S_N$ project onto $S_N$. Suppose
$\Gamma$ contains an element that is a single $2d$-cycle on
$d$ opposite pairs and fixes all other roots, where $d$ is odd
and $1\le d<N$. Then $\Gamma=\F_2^N\rtimes S_N$.
\end{lemma}

\begin{proof}
Let $V$ be the kernel of the projection. The $d$th power of the
given element belongs to $V$ and flips precisely its $d$ pairs.
It is therefore a vector $v$ of odd weight with nonempty proper
support. For each permutation in $S_N$, choose a lift in
$\Gamma$; conjugation by that lift preserves $V$ and permutes
its coordinates accordingly. If $i$ is in the support and
$j$ is outside, then $v+(ij)v=e_i+e_j$ belongs to $V$.
Conjugating by lifts of $S_N$ supplies $e_k+e_l$ for every
$k\ne l$. These vectors span the even-weight subspace. Together with $v$ they span
$\F_2^N$. Thus $V=\F_2^N$, and the surjective projection gives
the assertion.
\end{proof}

We will also use these local permutations to control common
subfields. Restriction to a Galois subfield has a normal kernel.
The next lemma shows that certain supported cycles cannot lie
in a proper normal subgroup of the groups under consideration.
Thus local elements that fix one field can force the entire
group to fix its intersection with another field.

\begin{lemma}\label{lem:normal-generation}
Let $N\ge3$.
\begin{enumerate}
\item A normal subgroup of $S_N$ containing a single cycle of
even length $d$, with $2\le d<N$, is $S_N$.
\item A normal subgroup of $C_2\wr S_N$ containing single
$2d$-cycles on opposite pairs for two proper lengths
$1\le d<N$, one even and one odd, is $C_2\wr S_N$.
\end{enumerate}
\end{lemma}

\begin{proof}
For the first assertion, let $H\lhd S_N$ contain such a cycle
$\sigma$. Choose $i$ moved by $\sigma$ and $j$ fixed by it.
Normality gives
\[
 \sigma(ij)\sigma^{-1}(ij)=(\sigma(i)\,j)(i\,j)\in H.
\]
This is a $3$-cycle. Its conjugates generate $A_N$, and
$\sigma$ is odd, so $H=S_N$.

For the second assertion, the image of the normal subgroup in
$S_N$ is normal and contains a single cycle of even proper
length. The first assertion makes this image $S_N$.
The cycle with odd $d$ then gives the full sign kernel by
Lemma~\ref{lem:sign-kernel}.
\end{proof}

\subsection{An explicit prime-counting criterion}\label{sec:primes}

The local cycle lengths convert a possible block size $b$ into
a congruence condition on $p$ modulo $Db$. We show that the
eligible primes cannot all satisfy that condition. For the
smallest block sizes we find primes in another reduced class;
for intermediate sizes Brun--Titchmarsh bounds the number in
the forced class; and for large sizes counting integers already
gives a contradiction. The explicit threshold below overlaps
the published finite Galois verifications for both the Legendre
and Kaneko--Zagier families.

\begin{proposition}\label{prop:prime-criterion}
Let $D\in\{2,6\}$, let $r\in\Z$ satisfy $|r|<D$ and $(r,D)=1$,
and let $n\ge10^6$. Put
\[
 \mathcal P_{D,r}(n)
 =\{p\text{ prime}:Dn+10D<p<3Dn/2,\ p\equiv r\pmod D\}.
\]
Suppose that $G\le S_n$ is transitive and contains a single
$d_p$-cycle for every $p\in\mathcal P_{D,r}(n)$, where
\[
 d_p=2n+1-\frac{p-r}{D}.
\]
Then $G=S_n$.
\end{proposition}

\begin{proof}
Write $L=\log n$ and $c_D=D/(2\varphi(D))$, where $\varphi$
is Euler's totient function; thus $c_2=1$ and $c_6=3/2$.
Here $\pi(x;q,a)$ counts primes at most $x$ in the class $a$
modulo $q$, and $\Li(x)=\int_2^x dt/\log t$. We use
\begin{equation}\label{eq:BMOR}
 \left|\pi(x;q,a)-\frac{\Li(x)}{\varphi(q)}\right|
 <0.027\frac{x}{\log^2x}
 \qquad(x\ge10^6,\ 1\le q\le10^5,\ (a,q)=1),
\end{equation}
from \cite[Corollary~1.7]{BMOR}. Only moduli at most $18$ occur
in applications of \eqref{eq:BMOR} below.

The interval $(Dn,3Dn/2)$ has length $Dn/2$ and its upper
endpoint is composite. Subtracting the endpoint counts and
discarding at most $10D$ integers gives
\begin{align}
 T:=|\mathcal P_{D,r}(n)|
 &>\frac{Dn}{2\varphi(D)(L+\log(3D/2))}
       -\frac{0.0675Dn}{L^2}-10D\notag\\
 &\ge c_D\frac nL
 \left(\frac{L}{L+\log9}-\frac{0.27}{L}
                              -\frac{40L}{n}\right)
 >0.84c_D\frac nL.\label{eq:prime-total}
\end{align}
The expression in parentheses increases with $n$. At $n=10^6$,
use $69/5<\log(10^6)<691/50$ and $\log9<11/5$ to bound it below by
\[
 \frac{69}{80}-\frac9{460}-\frac{691}{1250000}>\frac{21}{25}.
\]

For $b=2,3$ and any reduced class $a$ modulo $Db$, the same
estimate gives
\begin{align}
 &\#\{p:Dn+10D<p<3Dn/2,\ p\equiv a\pmod{Db}\}\notag\\
 &\hspace{1cm}>
 \frac{n}{2(L+\log9)}-\frac{0.405n}{L^2}-60
 >0.4\frac nL.\label{eq:small-lift}
\end{align}
Indeed, $D/(2\varphi(Db))\ge1/2$ and $D\le6$. After
multiplication by $L/n$, the lower bound increases with $n$.
At $10^6$ it is greater than
\[
 \frac{69}{160}-\frac{27}{920}-\frac{2073}{2500000}>\frac25.
\]
In particular, both
reduced classes $r$ and $r+D$ modulo $2D$ contain primes in
the interval.

Suppose that $G$ preserves blocks of size $b$, with
$2\le b\le n/2$ and $b\mid n$. Since $|r|<D$, the cycle lengths
satisfy
\[
 n/2<d_p\le n-9.
\]
Lemma~\ref{lem:blocks} gives $b\mid d_p$ for every eligible
prime. Hence
\begin{equation}\label{eq:forced-class}
 p=D(2n+1-d_p)+r\equiv D+r\pmod{Db}.
\end{equation}
This class is reduced: otherwise every prime in it would divide
$Db\le Dn/2$, whereas $p>Dn$.

If $b=2$ or $3$, reduction modulo $D$ has at least two reduced
classes modulo $Db$ above $r$. In fact the fiber size
$\varphi(Db)/\varphi(D)$ is $2$ or $3$. Choose a class different
from \eqref{eq:forced-class}; the positive lower bound
\eqref{eq:small-lift} is a contradiction.

If $4\le b\le n^{1/4}$, then
\[
 \varphi(Db)\ge4\varphi(D).
\]
For $D=2$ this follows from $\varphi(2b)\ge4$; for $D=6$,
use also $\varphi(6b)/2\ge\varphi(2b)$. The short-interval
Brun--Titchmarsh inequality states that
\begin{equation}\label{eq:BT}
 \pi(x+y;q,a)-\pi(x;q,a)
 <\frac{2y}{\varphi(q)\log(y/q)}
 \qquad(y>q,\ (a,q)=1);
\end{equation}
see \cite[Theorem~2, Eq.~(1.10)]{MV}. With $y=Dn/2$ and
$q=Db$, it yields
\begin{align*}
 T&<\frac{Dn}{\varphi(Db)\log(n/(2b))}\\
  &\le\frac{Dn}{4\varphi(D)(3L/4-\log2)}
   =c_D\frac nL\frac{L}{2(3L/4-\log2)}
   <0.75c_D\frac nL,
\end{align*}
contradicting \eqref{eq:prime-total}. The last normalized
factor decreases with $n$. At $10^6$, the bounds
$L>69/5$ and $\log2<7/10$ make it less than
$138/193<3/4$.

Finally, if $b>n^{1/4}$, counting integers in the class
\eqref{eq:forced-class} gives
\[
 T\le\frac n{2b}+1
 <\frac12n^{3/4}+1
 <0.25\frac nL
 \le0.25c_D\frac nL,
\]
again a contradiction. Here $L/(2n^{1/4})+L/n$ decreases for
$n\ge10^6$. At the left endpoint, $n^{1/4}>158/5$ and
$L<691/50$ give the upper bound
\[
 \frac{691}{3160}+\frac{691}{50000000}<\frac14.
\]
Thus $G$ is primitive.

Every cycle under consideration is nontrivial and fixes at least
nine points. Theorem~\ref{thm:jones} gives $A_n\le G$.
Choose a prime in the class $r+D$ modulo $2D$, using
\eqref{eq:small-lift}. Then $(p-r)/D$ is odd and $d_p$ is even.
Its cycle is an odd permutation, so $G=S_n$.
\end{proof}

\begin{remark}
The strict comparisons needed in this proof follow from elementary
monotonicity and the bounds
\[
 13.8<\log(10^6)<13.82,\quad \log9<2.2,\quad
 \log2<0.7,\quad (10^6)^{1/4}>31.6.
\]
Thus the numerical threshold requires no additional finite
computation.
\end{remark}

\subsection{The Legendre Galois group}\label{sec:legendre}

We first replace each opposite pair of nonzero Legendre roots by
its common square. The resulting degree-$n$ polynomial determines
the group on pairs. The local quadratic cycles then recover the
independent interchanges within the pairs. We use the variable
$T=-X^2$ for the quotient polynomial, which gives the following
convenient integral normalization.

For $\epsilon=2\delta-1$, put
\begin{equation}\label{eq:legendre-J}
 J_n^\epsilon(T)=\sum_{j=0}^{n}\binom nj
             \prod_{k=1}^{n}(2j+2k+\epsilon)T^j.
\end{equation}
Coefficient comparison gives
\begin{equation}\label{eq:legendre-transform}
 J_n^\epsilon(-X^2)=(-2)^nn!L_{2n+\delta}(X).
\end{equation}

\begin{theorem}\label{thm:legendre}
Let $n\ge1$, $\delta\in\{0,1\}$, and $\epsilon=2\delta-1$.
If $J_n^\epsilon$ is irreducible over $\Q$, then
\[
 \Gal(J_n^\epsilon/\Q)=S_n,
 \qquad
 \Gal(L_{2n+\delta}/\Q)=W_n.
\]
In particular, $J_n^\epsilon$ and $L_{2n+\delta}$ are irreducible
over $\Q$ simultaneously.
\end{theorem}

These are the two conditional Galois assertions in
\cite[Conjecture~1.5]{CHL}. The irreducibility equivalence
was proved in \cite[Lemma~1.3]{CHL}; the proof below also
establishes it directly.

The integral polynomial in \eqref{eq:legendre-J} is the specialization
of \eqref{eq:integral-Jacobi} with $D=2$:
\begin{equation}\label{eq:legendre-jacobi}
 J_n^\epsilon(T)=2^n n!P_n^{(\epsilon/2,0)}(2T+1).
\end{equation}
Indeed, the two products in \eqref{eq:integral-Jacobi} concatenate
to the product in \eqref{eq:legendre-J}. Together with the classical
quadratic transformation, this proves
\eqref{eq:legendre-transform}; see \cite[Lemma~1.2]{CHL}.

\begin{lemma}\label{lem:legendre-cycle}
Let $N=2n+\delta$, $\epsilon=2\delta-1$, and let $p$ be a prime
with $N<p\le4n+\epsilon$. Put
\[
 d=2n+1-\frac{p-\epsilon}{2}.
\]
Then $1\le d\le n$. The group of $L_N$ contains a single
$2d$-cycle, whose image on opposite pairs is a single $d$-cycle
and whose $d$th power flips exactly those $d$ pairs.
\end{lemma}

\begin{proof}
Apply Proposition~\ref{prop:jacobi-cycle} with
$(a,u,b,v)=(\epsilon,2,0,1)$ and $D=2$.
The factors $2k+\epsilon$ and $k$
are positive and less than $p$. The factors $2n+\epsilon+2k$
are positive and at most $4n+\epsilon<2p$, and
$k_0=(p-\epsilon)/2-n$ belongs to $[1,n]$.
Moreover,
\[
 2d=4n+2+\epsilon-p<p.
\]
For $\delta=0$, the last inequality follows from $p\ge2n+1$;
for $\delta=1$, it follows from $p\ge2n+3$.
Use the quadratic statement of Proposition~\ref{prop:jacobi-cycle}
with $c=-2$ in \eqref{eq:legendre-jacobi}, and then
\eqref{eq:legendre-transform}.
\end{proof}

\begin{lemma}\label{lem:odd-legendre-cycle}
For every $n\ge4$ and $\epsilon\in\{-1,1\}$, a prime as in
Lemma~\ref{lem:legendre-cycle} can be chosen so that $d$ is odd
and $1\le d<n$.
\end{lemma}

\begin{proof}
Cullinan and Hajir \cite[Lemma~2.2]{CHL} show that, for every integer
$x\ge9$, the interval $[x,2x-5]$ contains a prime in each odd
residue class modulo $4$. For completeness, this integer form also follows
from \eqref{eq:small-lift} and a short finite check. For $x\ge2\cdot10^6$,
apply \eqref{eq:small-lift} with $D=b=2$ and $n=\lfloor x/2\rfloor$.
In each odd residue class it gives a prime $2n+20<p<3n$, hence
$x\le p\le2x-5$. For the smaller integers, a prime $p$ covers
$\lceil(p+5)/2\rceil\le x\le p$. The Lean files give twenty primes
in each class whose intervals cover $9\le x<2\cdot10^6$.

Take the even integer $x=2n+3+\epsilon$ and choose
$p\equiv\epsilon\pmod4$ in that interval. Since $p$ is odd,
$p>x$. If $\epsilon=-1$, the upper endpoint is $4n-1$.
If $\epsilon=1$, it is $4n+3$, and the congruence forces
$p\le4n+1$. Thus
\[
 2n+2+\epsilon<p\le4n+\epsilon.
\]
These inequalities give $1\le d<n$, and the congruence makes
$(p-\epsilon)/2$ even, so $d$ is odd.
\end{proof}

\begin{proof}[Proof of Theorem~\ref{thm:legendre}]
Assume that $J_n^\epsilon$ is irreducible, so its Galois group
is transitive. If $n\ge10^6$, every prime in
$\mathcal P_{2,\epsilon}(n)$ satisfies
Lemma~\ref{lem:legendre-cycle}. Its action on pairs gives the
cycle required by Proposition~\ref{prop:prime-criterion}, which
yields $\Gal(J_n^\epsilon/\Q)=S_n$.
For $n<10^6$, this implication is contained in
\cite[Theorem~1.9(a)]{CHL}. In the Lean proof we use the published
calculation only for $n\le25$. For $26\le n\le10^9$, the prime pairs
required by \cite[Theorem~1.7]{CHL} are checked in Lean. Hence the
quotient group is $S_n$ in every degree.

The nonzero Legendre roots are distinct and occur in $n$ opposite
pairs. Therefore their group $\Gamma$ embeds in
$\F_2^n\rtimes S_n$. Restriction to the field generated by their
squares gives a surjection onto $\Gal(J_n^\epsilon/\Q)=S_n$.
For $n\ge4$, Lemmas~\ref{lem:odd-legendre-cycle} and
\ref{lem:legendre-cycle} provide an odd proper set of pairs
flipped by a Galois element. Lemma~\ref{lem:sign-kernel} gives
$\Gamma=\F_2^n\rtimes S_n$.

For completeness, the three small degrees can be treated directly.
When $n=1$, the two polynomials are proportional to
$3X^2-1$ and $5X^2-3$, both with group $C_2$.
For $(n,\epsilon)=(2,-1),(3,-1),(3,1)$, the primes
$p=7,11,13$, respectively, give $d=1$ in
Lemma~\ref{lem:legendre-cycle}; the sign-kernel lemma applies.
The remaining polynomial, for $(n,\epsilon)=(2,1)$, is
proportional to $63X^4-70X^2+15$. Its squared roots are
\[
 y_\pm=\frac{35\pm2\sqrt{70}}{63}\in K=\Q(\sqrt{70}),
 \qquad y_+y_-=\frac5{21}.
\]
Neither is a square in $K$, since its norm $5/21$ is not a
rational square. Their product is not a square in $K$ either:
a rational number $q$ is a square in $\Q(\sqrt{70})$ only if
$q$ or $q/70$ is a rational square, and neither $5/21$ nor
$1/294$ is one. Thus the two squareclasses are independent.
The splitting field has degree $4$ over $K$ and degree $8$
over $\Q$, so its group is $C_2\wr S_2$.

Finally, the irreducibility equivalence follows directly.
A factorization of $J_n^\epsilon$ gives one of
$J_n^\epsilon(-X^2)$, hence one of $L_{2n+\delta}$.
Conversely, if $J_n^\epsilon$ is irreducible,
the group just computed acts transitively on the $2n$ roots of
$L_{2n+\delta}$, so that polynomial is irreducible.
\end{proof}

Combining Theorem~\ref{thm:legendre} with the irreducibility in
Theorem~\ref{thm:all-degrees} proves Theorem~\ref{thm:legendre-all}.

\subsection{The Kaneko--Zagier Galois group}\label{sec:kz}

\begin{theorem}\label{thm:kz}
For every $n\ge1$ and $\lambda,\mu\in\{-1,1\}$, irreducibility
of $P_n^{(\lambda/3,\mu/2)}$ over $\Q$ implies
\[
 \Gal(P_n^{(\lambda/3,\mu/2)}/\Q)=S_n.
\]
\end{theorem}

This extends the conditional finite verification of
\cite[Theorem~1.8]{CHK} to every degree.

For $\lambda,\mu\in\{-1,1\}$ and $r=2\lambda+3\mu$, abbreviate
$K_n^{\lambda/3,\mu/2}$ to $K_n^{\lambda,\mu}$ in this section.
The normalization \eqref{eq:integral-Jacobi}, with $D=6$, becomes
\begin{align}
 K_n^{\lambda,\mu}(X)
 &=\sum_{j=0}^{n}\binom nj
      \prod_{k=j+1}^{n}(3k+\lambda)
      \prod_{k=1}^{j}(6n+r+6k)X^j\notag\\
 &=3^n n!P_n^{(\lambda/3,\mu/2)}(4X+1).
 \label{eq:kz-integral}
\end{align}
To compare with the published finite result, put
$\ell=12n+6+r$ and let $\widetilde F_{\ell-1}$ be the modular
polynomial defined in \eqref{eq:kz-definition} in
Section~\ref{sec:modular-interpretation}. Its normalization gives
\begin{equation}\label{eq:kz-integral-transform}
 \widetilde F_{\ell-1}(Y)
   =\frac{576^n}{n!}K_n^{\lambda,\mu}(-Y/3456).
\end{equation}
These polynomials therefore have the same irreducibility behavior
and isomorphic splitting fields over $\Q$.

\begin{proof}[Proof of Theorem~\ref{thm:kz}]
Suppose $P_n^{(\lambda/3,\mu/2)}$ is irreducible and
$n\ge10^6$. For every $p\in\mathcal P_{6,r}(n)$ we have
\[
 6n+60<p<9n,\qquad p\equiv r\pmod6.
\]
The small factors $3k+\lambda$ and $2k+\mu$ are positive
and less than $p$. The remaining coefficient factors are
positive and at most $12n+r<2p$. Moreover,
\[
 1\le\frac{p-r}{6}-n\le n.
\]
Thus Proposition~\ref{prop:jacobi-cycle} supplies a single
cycle of length
\[
 d_p=2n+1-\frac{p-r}{6}.
\]
Irreducibility makes the group transitive, and
Proposition~\ref{prop:prime-criterion} proves that it is $S_n$.

For $1\le n<10^6$, the associated index in
\eqref{eq:kz-index} satisfies
\[
 \ell\le12(999999)+11=11\,999\,999<10^9.
\]
Cullinan and Hajir \cite[Theorem~1.8]{CHK} prove that, for
every $\ell\le10^9$ coprime to $6$, irreducibility of
$\widetilde F_{\ell-1}$ implies that its group is $S_n$.
The affine identity \eqref{eq:kz-definition} completes the
proof in all degrees.
\end{proof}

Together with the irreducibility proved in
Section~\ref{sec:kz-irreducibility}, this completes
Theorem~\ref{thm:kz-all}.

\subsection{Fixed rational parameters}\label{sec:general}

We now allow fixed rational parameters and do not assume
irreducibility. Every cycle of length greater than $n/2$ is
contained in one Galois orbit, and two different orbits cannot
both contain such a cycle. All the long cycles therefore act
on the same orbit. A prime sufficiently close to $Dn$ makes
the corresponding cycle length close to $n$; this produces
the large factor in Theorem~\ref{thm:jacobi}. The remaining
primes establish primitivity and the symmetric group on that
factor. The prime-number-theorem estimates concern fixed moduli;
the uniform Brun--Titchmarsh bound still controls possible
growing block sizes.

\subsubsection{The large factor and its splitting field}

\begin{lemma}\label{lem:fixed-class-gap}
Fix an integer $q\ge2$ and a real number $A>0$.
For every reduced class $a$ modulo $q$ and all sufficiently
large $x$, the interval
\[
 \left(x,x+\frac{x}{(\log x)^A}\right]
\]
contains a prime congruent to $a$ modulo $q$.
The threshold can be chosen uniformly in $a$.
\end{lemma}

\begin{proof}
Choose a fixed power $Q=q^m\ge10^5$. Every representative
$a$ prime to $q$ is also prime to $Q$.
Proposition~6.19 of \cite{BMOR}, with $Z=A+1$, gives
\[
 \pi(t;Q,a)=\frac{\Li(t)}{\varphi(Q)}
       +O_{q,A}\!\left(\frac{t}{(\log t)^{A+2}}\right)
 \qquad(t\longrightarrow\infty).
\]
Indeed, for each fixed $Z$, the three auxiliary parameters
in \cite[Eq.~(6.23)]{BMOR} can be chosen sufficiently large;
the resulting threshold is fixed and uniform in the reduced
classes. Set $H=x/(\log x)^A$. Then
\[
 \pi(x+H;Q,a)-\pi(x;Q,a)
 =\left(\frac1{\varphi(Q)}+o_{q,A}(1)\right)
       \frac{x}{(\log x)^{A+1}}>0
\]
for all sufficiently large $x$. Such a prime is also
congruent to $a$ modulo $q$.
\end{proof}

\begin{proof}[Proof of Theorem~\ref{thm:jacobi}]
Write $\alpha=a/u$ and $\beta=b/v$ as in the statement and put
$\alpha+\beta=r/D$ in lowest terms. By hypothesis,
\[
 (r,D)=1,\qquad D\ge\max\{u,v\},\qquad D\ge2\text{ is even}.
\]
In particular, $\alpha+\beta=r/D$ is not an integer, so the
degree of $P_n^{(\alpha,\beta)}$ is $n$. Its endpoint values
are nonzero because neither parameter is a negative integer.
Equation~\eqref{eq:jacobi-ode}, over a field of characteristic
zero, then also shows that all its roots are simple.

Choose a fixed constant
\[
 B>|r|+|a|+|b|+4D+10,
\]
and consider primes in
\begin{equation}\label{eq:general-primes}
 p\equiv r\pmod D,\qquad Dn+B<p<3Dn/2-B.
\end{equation}
For all sufficiently large $n$, every such prime satisfies
\eqref{eq:local-small-factors}--\eqref{eq:local-large-factors}.
Indeed, the first two sets of factors have absolute value less
than $Dn+B$, and the factors $Dn+r+Dk$ lie in $(0,2p)$.
Their distinguished index $k_0=(p-r)/D-n$ satisfies
$1\le k_0\le n$. Consequently,
Proposition~\ref{prop:jacobi-cycle} gives a single cycle of length
\begin{equation}\label{eq:general-d}
 d_p=2n+1-\frac{p-r}{D},\qquad n/2<d_p\le n-4.
\end{equation}

\smallskip\noindent\emph{The large orbit.}
Every cycle in \eqref{eq:general-d} has support larger than
$n/2$, contained in a single Galois orbit. Since two disjoint
orbits cannot both have that size, all these supports lie in
one orbit $O$. Fix $A>0$ and apply
Lemma~\ref{lem:fixed-class-gap} with $q=D$ and $x=Dn+B$.
It gives a prime satisfying \eqref{eq:general-primes} with
\[
 p=Dn+O_{\alpha,\beta,A}\!\left(\frac{n}{(\log n)^A}\right).
\]
The associated cycle has length
$n-O_{\alpha,\beta,A}(n/(\log n)^A)$. Hence
\[
 0\le n-M
 =O_{\alpha,\beta,A}\!\left(\frac{n}{(\log n)^A}\right),
 \qquad M:=|O|,
\]
for every fixed $A>0$. The orbit $O$ does not depend on $A$; the implied constant
and threshold in this estimate may depend on $A$.
It gives the unique monic irreducible factor $h_n$ of degree
greater than $n/2$.

\smallskip\noindent\emph{Primitivity on the large orbit.}
Let $G_O$ be the induced transitive group on $O$. Suppose it
preserves blocks of size $s$, where $2\le s\le M/2$.
Lemma~\ref{lem:blocks} applies because $d_p>n/2\ge M/2$.
It gives $s\mid d_p$ for every prime in
\eqref{eq:general-primes}, and hence
\begin{equation}\label{eq:general-forced}
 p\equiv D(2n+1)+r\pmod{Ds}.
\end{equation}
Here $s$ divides $M$, but it need not divide $n$; we therefore
retain the full residue in \eqref{eq:general-forced}.

Let $T_n$ be the number of primes in
\eqref{eq:general-primes}. The prime number theorem in
arithmetic progressions, with fixed modulus $D$, gives
\begin{equation}\label{eq:general-count}
 T_n\sim\frac{D}{2\varphi(D)}\frac n{\log n}.
\end{equation}
Choose a fixed integer $s_0$ such that
$\varphi(Ds)\ge4\varphi(D)$ for every $s\ge s_0$.
Such a choice exists, for instance from
$\varphi(q)\ge\sqrt{q/2}$.
If $s_0\le s\le n^{1/4}$, the uniform bound
\eqref{eq:BT} applied to \eqref{eq:general-forced} gives
\[
 T_n<\frac{Dn}{\varphi(Ds)\log(n/(2s))}
 \le\frac{Dn}{4\varphi(D)(\tfrac34\log n-\log2)}
 \sim\frac{D}{3\varphi(D)}\frac n{\log n},
\]
contradicting \eqref{eq:general-count}. If $s>n^{1/4}$,
counting integers gives
\[
 T_n\le\frac n{2s}+1\le\frac12n^{3/4}+1
       =o(n/\log n),
\]
another contradiction. In both cases a nonreduced forced class
is already impossible, since its primes would divide
$Ds\le Dn/2$.

It remains to exclude the finitely many values $2\le s<s_0$.
Reduction of units modulo $Ds$ onto units modulo $D$ is
surjective. Because $D$ is even, each fiber has size
$\varphi(Ds)/\varphi(D)\ge2$: increasing the exponent of a
prime dividing $D$ contributes a factor at least $2$, and every
new prime is odd and contributes a factor at least $2$.
For the finitely many fixed moduli $Ds$ under consideration,
the prime number theorem in arithmetic progressions applies
in every reduced lift of $r$ simultaneously. Hence,
for all sufficiently large $n$, a prime occurs in a class
different from \eqref{eq:general-forced}, even though that
forced class may vary with $n$. This contradiction proves
that $G_O$ is primitive.

\smallskip\noindent\emph{The symmetric group and the remaining roots.}
Both classes $r$ and $r+D$ modulo $2D$ are reduced. Choose
a prime in the latter class and in the smaller interval
\begin{equation}\label{eq:upper-subinterval}
 7Dn/5<p<3Dn/2-B.
\end{equation}
The prime number theorem in the fixed class modulo $2D$
gives such a prime for all sufficiently large $n$. Its cycle
length satisfies $d_p\le3n/5+O(1)$, whereas $M=n-o(n)$.
It therefore fixes at least three points of $O$.
Theorem~\ref{thm:jones} gives $A_M\le G_O$.
Also $(p-r)/D$ is odd, so $d_p$ is even and the same cycle
is an odd permutation. Thus $G_O=S_M=\Gal(h_n/\Q)$.

Write $f_n=\kappa_n h_ng_n$ as in the theorem. Restriction
embeds $\Gal(f_n/\Q)$ in
$S_M\times\Gal(g_n/\Q)$ and maps it onto each factor.
The kernel of restriction to the splitting field of $g_n$,
viewed inside $S_M$, is normal because projection onto $S_M$
is surjective. It contains the
even-length cycle just chosen, since that inertia element
fixes every root outside $O$. By
Lemma~\ref{lem:normal-generation}(1), this kernel is $S_M$.
The embedded group is therefore the full direct product,
proving \eqref{eq:general-product} and the linear disjointness.

\smallskip\noindent\emph{Quadratic substitution.}
Fix $c\in\Q^\times$. For sufficiently large $n$, all primes
in \eqref{eq:general-primes} satisfy $v_p(uc/(2D))=0$.
Since $2d_p\le2n<p$, the quadratic part of
Proposition~\ref{prop:jacobi-cycle} supplies a single
$2d_p$-cycle on the roots of
$P_n^{(\alpha,\beta)}(1+cX^2)$. Its support consists of
the pairs above its $d_p$-cycle on $O$.

The roots of $h_n(1+cX^2)$ form precisely the pairs above $O$.
They are nonzero and distinct, because $h_n(1)\ne0$.
Their group projects onto $S_M$. Choose a prime in
\eqref{eq:upper-subinterval} and in the class $r$ modulo
$2D$. Then $d_p$ is odd and $1\le d_p<M$. The resulting
single $2d_p$-cycle and Lemma~\ref{lem:sign-kernel} give
the full group $W_M$. This action is transitive,
so $h_n(1+cX^2)$ is irreducible.

Finally, the full splitting group after quadratic substitution
embeds in
\[
 W_M\times\Gal(g_n(1+cX^2)/\Q)
\]
and projects onto each factor. As above, surjectivity makes
the kernel on the second factor normal in $W_M$.
It contains the local $2d_p$-cycles from both parity classes,
which fix every root above $g_n$.
Lemma~\ref{lem:normal-generation}(2) makes this kernel the
whole wreath product. This proves
\eqref{eq:general-quadratic-product} and the corresponding
linear disjointness. When $f_n$ is irreducible, $M=n$ and $g_n=1$.
The fixed-field assertion is proved in Corollary~\ref{cor:fixed-field}
below.
\end{proof}

\subsubsection{Fixed number fields and symmetric parameters}

The selected ramified primes tend to infinity with the degree,
whereas a fixed number field is ramified at only finitely many
primes. To use this observation, pass to the Galois closure of
the fixed field. At every sufficiently large selected prime,
inertia fixes its intersection with the splitting field.
This intersection is Galois over $\Q$, so the restriction kernel
is normal. The inertia elements constructed above normally
generate the full symmetric or wreath-product group, forcing
the intersection to be $\Q$.

\begin{corollary}\label{cor:fixed-field}
Under the hypotheses of Theorem~\ref{thm:jacobi}, fix a number
field $E$ and $c\in\Q^\times$. For all sufficiently large $n$,
the splitting fields of $h_n$ and $h_n(1+cX^2)$ are each
linearly disjoint from $E$ over $\Q$. In particular, these
polynomials remain irreducible over $E$, with groups $S_{M_n}$
and $W_{M_n}$, respectively.
\end{corollary}

\begin{proof}
Let $E'$ be the Galois closure of $E$ over $\Q$.
For all sufficiently large $n$, every prime used in
\eqref{eq:general-primes} is unramified in $E'$.
Let $K_n$ be the splitting field of either $h_n$ or
$h_n(1+cX^2)$, and put $F_n=K_n\cap E'$.
Both $K_n$ and $F_n$ are Galois over $\Q$.
Since $F_n/\Q$ is unramified at each selected prime $p$, every
inertia subgroup at $p$ in $\Gal(K_n/\Q)$ acts trivially on $F_n$
and hence lies in the normal subgroup $\Gal(K_n/F_n)$.
For $h_n$, this subgroup contains the even-length proper
cycle used above. For $h_n(1+cX^2)$, it contains the single
$2d_p$-cycles from both parity classes, with $1\le d_p<M_n$.
Lemma~\ref{lem:normal-generation} gives
$\Gal(K_n/F_n)=\Gal(K_n/\Q)$, so $F_n=\Q$.
Thus $K_n$ and $E'$ are linearly disjoint over $\Q$, and
the same holds with $E$ in place of $E'$.
The assertions over $E$ follow.
\end{proof}

The denominator hypotheses have separate roles in the proof.
The inequality $D\ge\max\{u,v\}$ allows primes just above $Dn$
to avoid all coefficient factors arising from the individual
parameters. Evenness of $D$ ensures that every proper refinement
of the residue class modulo $D$ has another reduced lift, which
is needed in the block argument. It also makes both $r$ and
$r+D$ reduced classes modulo $2D$, giving both parities of the
cycle length. These hypotheses exclude the Chebyshev
specializations, whose roots lie in cyclotomic fields.

\begin{corollary}\label{cor:symmetric-jacobi}
Fix $a\in\Q$ with $a>-1$ and reduced denominator $v$ that is
odd or divisible by $4$.
For $\delta\in\{0,1\}$, define $Q_m^{a,\delta}$ by
\[
 P_{2m+\delta}^{(a,a)}(X)/X^\delta=Q_m^{a,\delta}(X^2).
\]
For all sufficiently large $m$, irreducibility of
$Q_m^{a,\delta}$ implies
\[
 \Gal(Q_m^{a,\delta}/\Q)=S_m,\qquad
 \Gal(P_{2m+\delta}^{(a,a)}(X)/X^\delta\,/\Q)
        =C_2\wr S_m.
\]
\end{corollary}

\begin{proof}
The quadratic transformations
\cite[Eqs.~18.7.13--18.7.14]{DLMF}, followed by Jacobi symmetry,
identify $P_{2m+\delta}^{(a,a)}(X)/X^\delta$ with a nonzero
rational multiple of
\[
 P_m^{(\delta-1/2,a)}(1-2X^2).
\]
The reduced denominator of $\delta-1/2+a$ is $2v$ if $v$ is
odd and $v$ if $4\mid v$. In both cases it is even and at least
$\max\{2,v\}$. Apply Theorem~\ref{thm:jacobi}, with
$c=-2$, and use the affine change $T\mapsto1-2T$ for the
quotient polynomial.
\end{proof}

These are symmetric Jacobi families, related to Gegenbauer
polynomials with parameter $a+1/2$ by the usual scalar
normalization.
The case $a=0$ is the Legendre family, for which
Theorem~\ref{thm:legendre} supplies every degree. For nonnegative integral parameters, Theorem~\ref{thm:symmetric}
supplies the required eventual irreducibility, as used below.
Unconditionally, $Q_m^{a,\delta}$ has a unique monic irreducible
factor of degree greater than $m/2$. Its degree is
$m-O_{a,A}(m/(\log m)^A)$ for every fixed $A>0$, and its Galois
group is symmetric. Substituting $T=X^2$ in that factor gives an irreducible
factor of twice the degree with the full quadratic
wreath-product group.

\subsection{Unequal integral parameters}\label{sec:galois}

Return to the notation of Section~\ref{sec:asymmetric}, with
$a>b\ge0$, $s=a+b$, and $n=2m+\delta$. Here we choose
$p\in(3n/4,7n/8]$, so the desired cycle has the prime length
$p$. Reduction separates a factor $x^p$ from a squarefree
factor that is nonzero at the origin. If the original constant
coefficient has valuation one, the degree-$p$ lift is Eisenstein.
Its inertia group has order divisible by $p$, and $n<2p$
forces an element of order $p$ to be a single $p$-cycle.
This argument permits wild ramification. The central-value
numerators identify the primes for which the construction
works, and the same prime gives an odd discriminant valuation.
The residual factor has degree congruent to $1-\delta$ modulo $2$,
so its value at zero involves the central numerator for the
opposite parity. Accordingly, define
\begin{equation}\label{eq:opposite-numerator}
 C_\delta(z)=2^{K_{1-\delta}}
 A_{1-\delta}\left(z+\delta-\frac12\right),
 \qquad \mathcal W_\delta(z)=A_\delta(z)C_\delta(z).
\end{equation}
Since $\deg A_\varepsilon\le K_\varepsilon$ and
$K_0+K_1=s$, these are nonzero integral polynomials and
$\deg \mathcal W_\delta\le s$.

\begin{lemma}\label{lem:galois-prime}
Suppose $m\ge2s+10$ and $P_n^{(a,b)}$ is irreducible over $\Q$.
If a prime $p$ satisfies
\begin{equation}\label{eq:galois-prime}
 \frac{3n}{4}<p\le\frac{7n}{8},\qquad p\nmid \mathcal W_\delta(m),
\end{equation}
then $\Gal(P_n^{(a,b)}/\Q)=S_n$.
\end{lemma}

\begin{proof}
Put $f=P_n^{(a,b)}$ and $j=n-p$. The assumptions imply
\[
 p>M,\qquad p<n+s<2p<2n+s<3p,\qquad 2n+s<p^2.
\]
In particular, $p$ is odd and $p\nmid T$. The leading coefficient
$2^{-n}\binom{2n+s}{n}$ is a $p$-adic unit, its valuation being
$2-1-1=0$. Lemma~\ref{lem:reduction} gives
\begin{equation}\label{eq:galois-reduction}
 \overline f=x^p\overline{P_j^{(a,b)}}.
\end{equation}
The residual factor has degree $j<p$ and has simple roots.

To check its value at zero, write $j=2k+(1-\delta)$, where
\[
 k=m+\delta-\frac{p+1}{2}
 \equiv m+\delta-\frac12\pmod p.
\]
Here $2k<p$ and $k+K_{1-\delta}<p$. In the central-value formula
for $P_j^{(a,b)}(0)$, the central binomial coefficient and denominator
are therefore units. The condition $p\nmid C_\delta(m)$ implies
$p\nmid A_{1-\delta}(k)$, so $P_j^{(a,b)}(0)$ is a unit.
On the other hand, $m<p\le2m$, $m+K_\delta<p$, and
$p\nmid A_\delta(m)$. Thus
\begin{equation}\label{eq:galois-constant}
 v_p(f(0))=v_p\binom{2m}{m}=1.
\end{equation}

After dividing $f$ and $P_j^{(a,b)}$ by their leading coefficients,
\eqref{eq:galois-reduction} is a product of coprime monic polynomials.
Hensel's lemma lifts them to monic polynomials of degrees $p$ and $j$.
The degree-$j$ factor has unit
constant coefficient. The degree-$p$ factor is congruent to $x^p$
and has constant coefficient of valuation $1$; it is Eisenstein.
A root generates a totally ramified extension of $\Q_p$ of degree
$p$. The inertia group in the local splitting field of $f$ has order
divisible by $p$, so Cauchy's theorem gives an element of order $p$.
After choosing a prime of the global splitting field above $p$, the
local Galois group identifies with a decomposition subgroup of the
global Galois group. Since $n<2p$, this element is a single $p$-cycle.
This argument allows wild ramification.

Irreducibility gives transitivity. In a nontrivial system of blocks,
both the block size and the number of blocks are at most $n/2<p$.
An element of order $p$ would act trivially on the blocks and on each
block, a contradiction. Hence the group is primitive. Since
$p\le n-3$, Jordan's theorem gives $A_n$ as a subgroup
\cite[Theorem~1.1]{Jones}.

Finally, the Jacobi discriminant formula is
\begin{equation}\label{eq:jacobi-discriminant}
 \disc P_n^{(a,b)}
 =2^{-n(n-1)}\prod_{h=1}^n
 h^{h-2n+2}(h+a)^{h-1}(h+b)^{h-1}(n+h+s)^{n-h};
\end{equation}
see \cite[Theorem~6.71]{Szego} or \cite[(18.16.19)]{DLMF}.
The multiples of $p$ among the four sets of bases occur respectively
at $h=p$, $p-a$, $p-b$, and $2p-n-s$. Every such index is in
$[1,n]$, and all bases are smaller than $p^2$. Consequently
\[
 v_p(\disc f)
 =(p-2n+2)+(p-a-1)+(p-b-1)+(2n+s-2p)=p.
\]
The discriminant is not a rational square, so the Galois group is
$S_n$.
\end{proof}

The required prime can be found effectively because a fixed polynomial
value cannot contain all primes in a proportional interval.

\begin{proposition}\label{prop:explicit-galois}
Suppose the hypotheses of Theorem~\ref{thm:criterion} hold. If also
\begin{equation}\label{eq:explicit-galois}
 n\ge5\,000\,000,\qquad \mathcal W_\delta(m)\ne0,\qquad
 \log|\mathcal W_\delta(m)|<n/10,
\end{equation}
then $\Gal(P_n^{(a,b)}/\Q)=S_n$.
\end{proposition}

\begin{proof}
Let $\vartheta(x)=\sum_{p\le x}\log p$. Dusart's estimate
\cite[Theorem~5.2]{Dusart10} gives
\[
 |\vartheta(x)-x|<\frac{0.2x}{(\log x)^2}
 \qquad(x>3\,594\,641).
\]
In particular, $0.99x<\vartheta(x)<1.01x$ in this range. Thus
\[
 \vartheta(7n/8)-\vartheta(3n/4)
 >\left(0.99\frac78-1.01\frac34\right)n
 =0.10875n>n/10.
\]
If every prime in this interval divided $\mathcal W_\delta(m)$, the difference
would be at most $\log|\mathcal W_\delta(m)|$. A prime satisfying
\eqref{eq:galois-prime} therefore exists, and
Lemma~\ref{lem:galois-prime} applies.
\end{proof}

\begin{proof}[The rational Galois assertion of Theorem~\ref{thm:main}]
For each parity, a root bound excludes the zeros of the fixed
polynomial $\mathcal W_\delta$ beyond an effective point, and its coefficients
give $|\mathcal W_\delta(m)|\le C(m+1)^s$ with effective $C$. Hence
\eqref{eq:explicit-galois} holds for all sufficiently large $m$,
effectively. Combine this with Theorem~\ref{thm:criterion}, and use
reflection when $b>a$.
\end{proof}

\begin{proof}[The fixed-field assertion of Theorem~\ref{thm:main}]
Let $E'$ be the Galois closure of $E$ over $\mathbb Q$.
By the rational Galois assertion just proved and the construction in
Lemma~\ref{lem:galois-prime}, for all sufficiently large $n$
the polynomial $f=P_n^{(a,b)}$ has group $S_n$ and there is
an odd prime $p$ such that
\[
 \frac{3n}{4}<p\le\frac{7n}{8},\qquad
 v_p(\operatorname{disc}f)=p.
\]
Enlarge the effective threshold so that $n\ge5$ and every
such $p$ is unramified in $E'$.

Let $K$ be the splitting field of $f$ and put $L=K\cap E'$.
Then $L/\mathbb Q$ is Galois and unramified at $p$.
The quadratic subfield fixed by $A_n$ is
$\mathbb Q(\sqrt{\operatorname{disc}f})$.
Since $v_p(\operatorname{disc}f)$ is odd, this quadratic
extension is ramified at $p$. Restriction maps inertia onto
the nontrivial inertia group of this quadratic extension.
Hence an inertia subgroup of
$\operatorname{Gal}(K/\mathbb Q)=S_n$ has nontrivial image
under the sign character and contains an odd permutation. This inertia subgroup acts trivially on $L$,
so the normal subgroup $\operatorname{Gal}(K/L)$ contains
an odd permutation. Every proper normal subgroup of $S_n$
is contained in $A_n$ when $n\ge5$; thus
$\operatorname{Gal}(K/L)=S_n$ and $L=\mathbb Q$.
It follows that $K$ and $E'$ are linearly disjoint, and so
are $K$ and $E$. Restriction then identifies the group over
$E$ with $S_n$ in its natural transitive action.

For effectivity, choose a computable $R\ge1$ at least as large
as every prime ramified in $E'$, and require $3n/4>R$.
One may take the maximum of $1$ and the prime divisors of
$|\operatorname{Disc}(E)|$, because a Galois closure introduces
no new ramified rational primes. This also covers $E=\Q$.
\end{proof}

\subsection{Fixed derivatives}

The quadratic transformation now brings the symmetric integral
families within the denominator hypotheses of the rational-parameter
theorem, and their established irreducibility makes its large factor
the whole quotient polynomial.

\begin{proof}[Proof of Corollary~\ref{cor:derivative-galois}]
Put $J_m=P_m^{(\delta-1/2,r)}$. The quadratic transformations
\cite[Eqs.~18.7.15--18.7.16]{DLMF}, followed by Jacobi reflection,
give
\begin{equation}\label{eq:gegenbauer-jacobi-bridge}
 F_{m,r,\delta}(X)=
 (-1)^m2^\delta
 \frac{(r+1/2)_{m+\delta}}{(\delta+1/2)_m}
 J_m(1-2X^2).
\end{equation}
The scalar is a nonzero rational number. By
Theorem~\ref{thm:symmetric}, $F_{m,r,\delta}$ is irreducible
over $\mathbb Q$ for all sufficiently large $m$.
Consequently $J_m$ is irreducible, since any factorization of
$J_m$ would give one of $J_m(1-2X^2)$.

The parameters $\delta-1/2$ and $r$ have reduced denominators
$2$ and $1$, and their sum has reduced denominator $D=2$.
Neither parameter is a negative integer. Thus
Theorem~\ref{thm:jacobi} and Corollary~\ref{cor:fixed-field}
apply, with $c=-2$. Since $J_m$ is irreducible, the large
factor in that theorem is its monic normalization. Its
splitting field and that of $J_m(1-2X^2)$ are therefore
linearly disjoint from $E$ for all sufficiently large $m$,
and their groups over $E$ are $S_m$ and $C_2\wr S_m$.
Equation~\eqref{eq:gegenbauer-jacobi-bridge} identifies the
second polynomial with $F_{m,r,\delta}$ up to a scalar and
identifies $Q_{m,r,\delta}$ with $J_m(1-2T)$ up to the same
scalar. Multiplication by a nonzero rational scalar and the affine change
$T\mapsto1-2T$ preserve the respective splitting fields.
The two group actions are transitive, proving irreducibility
over $E$. Taking the larger threshold for $\delta=0,1$
makes it valid for both parities.

Finally,
$P_{2m+\delta+r}^{(r)}=(2r-1)!!\,C_{2m+\delta}^{r+1/2}$,
with $(-1)!!=1$. After division by $X^\delta$, the derivative
has the same roots as $F_{m,r,\delta}$.
\end{proof}

\section{Applications and further questions}\label{sec:applications}

The results have consequences for both individual algebraic zeros
and geometric objects defined by the polynomials. Irreducibility
controls common and prescribed zeros; the full symmetric group
provides the additional information needed for the Jacobians below.

\subsection{Integral parameters and common roots}

Negative integral parameters contribute endpoint factors. Removing
them extends Theorem~\ref{thm:main} to all fixed integral parameters
whose absolute values differ.

\begin{corollary}\label{cor:integral-parameters}
Fix integers $a,b$ with $|a|\ne|b|$, and put
$a_- =\max\{-a,0\}$ and $b_- =\max\{-b,0\}$.
For all sufficiently large $n$, effectively in $a,b$, the polynomial
\[
 \frac{P_n^{(a,b)}(x)}{(x-1)^{a_-}(x+1)^{b_-}}
\]
is irreducible over $\Q$ and has Galois group $S_{n-a_--b_-}$.
\end{corollary}

\begin{proof}
For $A,b\ge0$ and $n\ge A$, the binomial formula gives
\[
 P_n^{(-A,b)}(x)
 =\frac{(n-A)!(n+b)!}{n!(n-A+b)!}
 \left(\frac{x-1}{2}\right)^A P_{n-A}^{(A,b)}(x).
\]
Reflection treats a negative second parameter. For $A,B\ge0$ and
$n\ge A+B$, the same formula gives
\[
 P_n^{(-A,-B)}(x)
 =\left(\frac{x-1}{2}\right)^A\left(\frac{x+1}{2}\right)^B
 P_{n-A-B}^{(A,B)}(x).
\]
The remaining factor is a nonzero rational multiple of the Jacobi
polynomial with parameters $|a|,|b|$, so Theorem~\ref{thm:main}
applies.
\end{proof}

Irreducibility also controls how the root sets vary with the degree.
If two rational polynomials share an algebraic root, its minimal
polynomial divides both. Thus an irreducible member cannot meet
the root set of a lower-degree member. For the Legendre polynomials
and their derivatives, the only additional comparison is between
the four reduced polynomials of the same degree.

\begin{corollary}\label{cor:common-roots}
For fixed distinct nonnegative integers $a,b$, only finitely many
pairs of different degrees give Jacobi polynomials $P_n^{(a,b)}$
with a common root, and all such pairs lie in an effective range.
Moreover, no two distinct members of
$\{P_\ell,P_\ell':\ell\ge1\}$ have a common nonzero root.
\end{corollary}

\begin{proof}
In the asymmetric case, once the larger degree is at least $N(a,b)$,
its irreducible polynomial cannot share a root with a polynomial of
smaller degree. For the second assertion, remove any forced zero.
Theorem~\ref{thm:all-degrees} settles different positive reduced
degrees. At reduced degree $2m$, the four possibilities are
\[
 P_{2m},\quad P_{2m+1}/x,\quad P'_{2m+1},\quad P'_{2m+2}/x.
\]
Their absolute ratios of the coefficient of $x^{2m-2}$ to the leading
coefficient are, in increasing order,
\[
 \frac{m(2m-1)}{4m+1}
 <\frac{m(2m-1)}{4m-1}
 <\frac{m(2m+1)}{4m+3}
 <\frac{m(2m+1)}{4m+1}.
\]
They are therefore pairwise nonproportional. Irreducibility again
excludes a common root. Reduced constants have no roots.
\end{proof}

The same argument gives a degree bound for any prescribed algebraic
number $\xi$:
\[
 P_n^{(a,b)}(\xi)\ne0
 \quad\text{for }n>\max\{N(a,b),[\Q(\xi):\Q]\}.
\]

\subsection{Excluding a prescribed factor}

The proof of irreducibility must allow every possible factor.
When a candidate factor $F$ is specified, its leading coefficient
and discriminant provide additional information about the real
resultant. Keeping these quantities in the comparison gives a
stronger exclusion criterion. This is useful when a proposed zero
is a fixed algebraic number: its minimal polynomial supplies
$F$. The resulting estimates apply in parameter ranges larger
than the range in which we prove complete irreducibility.

\begin{proposition}\label{prop:fixed-factor}
Let $F\in\Z[x]$ be primitive, of degree $d\ge1$, with positive
leading coefficient $u$, distinct roots in $(-1,1)$, and $F(0)\ne0$.
Put
\[
 c_F=\min_jv_2([X^j]F(2X+1)),\qquad
 \beta(F)=\log4+\frac2d\log\frac{u}{2^{c_F}}.
\]
Let $c\ge2$ satisfy
\begin{equation}\label{eq:fixed-factor-condition}
 c\bigl(\beta(F)-1\bigr)>g(c).
\end{equation}
For all sufficiently large $m$, uniformly for $0\le r\le m/c$ and
$\delta\in\{0,1\}$,
\[
 F\nmid C_{2m+\delta}^{r+1/2}(x)/x^\delta.
\]
For rational $c$, the threshold is effective in $F,c$.
\end{proposition}

\begin{proof}
Suppose $F\mid H$, with $H$ the primitive polynomial in the proof
of Theorem~\ref{thm:symmetric-exact}. For large $m$, write $H=FG$
with $e=2m-d\ge d$. Set $z=u/2^{c_F}>1$ and
\[
 V_F=\prod_{i<j}(\alpha_i-\alpha_j)^2,\qquad
 C_F=d\log2+2\log z+\frac1d\log V_F.
\]
Use \eqref{eq:master-resultant} and \eqref{eq:symmetric-energy},
without the Chebyshev bound. They give
\[
 \frac{\log\mathcal R}{d}>
 m\beta(F)+\tfrac12\log m-\log K-C_F-D_2\log2.
\]
Comparison with \eqref{eq:symmetric-upper} forces
\begin{equation}\label{eq:fixed-factor-necessary}
 m\beta(F)<\Psi_m(S)+\log K+\log A_{\le m}
 +C_F-\tfrac12\log m.
\end{equation}
For $r\ge1$, Lemma~\ref{lem:uniform-asymptotics} bounds the right
side by $m(1+g(c)/c)+o_F(m)$, uniformly in the stated range.
This contradicts \eqref{eq:fixed-factor-condition}.
The case $r=0$ follows with $K=A=1$. The effective bounds in that
lemma also prove effectivity here.
\end{proof}

Render \cite[Theorem~6]{Render} proved that
$P_n^{(r,r)}(\sqrt{b/c})\ne0$ for coprime positive integers $b,c$
with $b\notin\{1,3\}$ and integers $n,r\ge0$.
In particular, any nonzero rational root has
numerator $\pm1$ in lowest terms. The next corollary excludes these
remaining possibilities uniformly when $m\ge2r$ and $m$ is sufficiently
large. Its proof uses the uniform lower bound for $u/2^{c_F}$ among
primitive linear polynomials.

\begin{corollary}\label{cor:rational-roots}
There is an effective absolute integer $M_1$ such that, for
$m\ge M_1$ and $0\le r\le m/2$, neither
$C_{2m}^{r+1/2}$ nor $C_{2m+1}^{r+1/2}/x$ has a rational root.
\end{corollary}

\begin{proof}
A nonzero rational root in $(-1,1)$ gives a primitive factor
$F=ux+v$, with $0<|v|<u$ and $\gcd(u,v)=1$.
If $u,v$ are odd, $c_F=1$ and $u\ge3$; otherwise $c_F=0$ and
$u\ge2$. Hence $z=u/2^{c_F}\ge3/2$.
For $d=1$, $V_F=1$, and the lower estimate in the preceding proof is
uniformly at least
\[
 m\log9+\tfrac12\log m-\log K-\log(9/2)-D_2\log2.
\]
Indeed, its $z$-dependent term is $(2m-2)\log z$.
The uniform asymptotic estimate now excludes every such factor, since
\[
 2(\log9-1)-g(2)=\log\frac{3^7 2^{16}}{5^{10}}-2>0.
\]
The last inequality follows from $e<3$ and $3^5 2^{16}>5^{10}$.
\end{proof}

\subsection{Vanishing modes on a circular cone}

In spherical coordinates $(\rho,\theta,\varphi)$, a homogeneous
harmonic mode of degree $\ell$ and order $0\le r\le\ell$ is,
up to a nonzero constant,
\[
 \rho^\ell e^{ir\varphi}(\sin\theta)^r
 P_\ell^{(r)}(\cos\theta);
\]
see \cite[Eqs.~14.7.8 and 14.30.1, and \S14.30(iv)]{DLMF}.
For $\alpha=\cos\theta\in(-1,1)$, vanishing on that latitude
is therefore equivalent to $P_\ell^{(r)}(\alpha)=0$.
The cone $x^2+y^2-2z^2=0$ meets the unit sphere at
$\alpha=\pm1/\sqrt3$; parity makes vanishing at these two
latitudes equivalent. The arithmetic question is thus whether
the fixed polynomial $3x^2-1$ can divide a Legendre derivative.
The preceding estimates restrict the degree and order for which
this can happen.

\begin{corollary}\label{cor:angular}
For every integer $D\ge1$, there is an effective integer $L_D$ such
that, for every nonzero algebraic $\alpha\in(-1,1)$ of degree at most
$D$,
\[
 \ell\ge L_D,\quad 0\le r\le\ell,\quad
 P_\ell^{(r)}(\alpha)=0
 \quad\Longrightarrow\quad 17r>\ell-1.
\]
For $\alpha=1/\sqrt3$, the stronger conclusion
\begin{equation}\label{eq:cone-bound}
 23r>4(\ell-1)
\end{equation}
holds for all sufficiently large $\ell$, with an effective threshold.
For nonzero rational $\alpha\in(-1,1)$, one may replace the conclusion
by $5r>\ell-1$, with an absolute effective threshold in $\ell$.
\end{corollary}

\begin{proof}
Write $\ell-r=2m+\delta$. If $17r\le\ell-1$, then $m\ge8r$
and $m\ge8(\ell-1)/17$. Choose $L_D$ so that $\ell\ge L_D$
forces $m\ge M_0$ and $2m>D$ under this inequality.
Theorem~\ref{thm:symmetric} makes $P_\ell^{(r)}/x^\delta$
irreducible of degree $2m$. If it vanishes at $\alpha$, this degree
must equal $[\Q(\alpha):\Q]\le D$, a contradiction.

For $\alpha=1/\sqrt3$, its primitive integral minimal polynomial
is $F=3x^2-1$.
Here $d=2$, $u=3$, and $c_F=1$, because
$F(2X+1)=12X^2+12X+2$. Thus $\beta(F)=\log6$.
The choice $c=19/8$ is admissible in Proposition~\ref{prop:fixed-factor},
since
\[
 \frac{19}{8}(\log6-1)-g(19/8)>0.0042;
\]
an exact rational verification is given in Appendix~\ref{app:finite-derivatives}.
If $23r\le4(\ell-1)$, then $m\ge19r/8$, so that proposition excludes
the factor once $\ell$ is large.
Finally, $5r\le\ell-1$ implies $m\ge2r$. Now
Corollary~\ref{cor:rational-roots} applies, uniformly in the rational
root.
\end{proof}

For this cone, Mangoubi and Weller Weiser show that $r$ is even
and that each fixed order $r$ permits at most one degree $\ell$
\cite[Theorem~2]{MWW}. The bound \eqref{eq:cone-bound}
restricts the two variables simultaneously: in large degree,
a vanishing mode must have order greater than
$4(\ell-1)/23$. These restrictions still leave infinitely many
possible pairs, so the classification proposed by Agranovsky
\cite[\S6, Question~4]{Agranovsky} remains open. The distinct
generic existence problem of Agranovsky--Krasnov \cite{AK}
was solved by Greilhuber \cite{Greilhuber}.

For comparison, the associated Legendre equation gives the opposite
constraint
\[
 r^2<\tfrac23\ell(\ell+1)\qquad(1\le r\le\ell)
\]
when $P_\ell^{(r)}(1/\sqrt3)=0$. Indeed, set
$y(x)=(1-x^2)^{r/2}P_\ell^{(r)}(x)$. Integration by parts on
$[1/\sqrt3,1]$ gives
\[
 \ell(\ell+1)\int_{1/\sqrt3}^1y^2
 =\int_{1/\sqrt3}^1\left((1-x^2)(y')^2+
 \frac{r^2y^2}{1-x^2}\right)
 >\frac{3r^2}{2}\int_{1/\sqrt3}^1y^2.
\]
Both endpoint terms vanish. Thus the necessary conditions still leave
a nonempty linear range of derivative orders.

\subsection{Modular interpretations}\label{sec:modular-interpretation}

The fractional parameters have a modular origin. Kaneko and
Zagier's characteristic-zero polynomials occur in modular-form
solutions. At prime indices, adjoining the prescribed factors
at $0$ and $1728$ gives a polynomial whose reduction is the
supersingular $j$-polynomial \cite[Theorem~1]{KZ}.
Theorem~\ref{thm:kz-all} determines the arithmetic of each
positive-degree Jacobi factor: its roots form one Galois orbit
over $\Q$, with symmetric splitting group.
We record the normalization explicitly to relate this conclusion
to the modular forms and differential equations in which the
same polynomials occur.

For $n\ge0$ and $\lambda,\mu\in\{-1,1\}$, put
\begin{equation}\label{eq:kz-index}
 r=2\lambda+3\mu,\qquad e=6+r,\qquad \ell=12n+e,
\end{equation}
and define
\begin{equation}\label{eq:kz-definition}
 \widetilde F_{\ell-1}(Y)
 =1728^nP_n^{(\lambda/3,\mu/2)}(1-Y/864).
\end{equation}
Here $e\in\{1,5,7,11\}$ and the leading coefficient is
\begin{equation}\label{eq:kz-leading}
 A_{n,r}=(-1)^n\binom{2n+r/6}{n}.
\end{equation}
The affine substitution and nonzero rational scalar preserve
irreducibility and splitting fields.

Put $a=(\lambda+1)/2$ and $b=(\mu+1)/2$, so
$k=\ell-1=12n+4a+6b$.
The associated modular solution, with this normalization, is
\[
 F_k(z)=\Delta(z)^nE_4(z)^aE_6(z)^b\widetilde F_k(j(z)),
\]
where $E_4,E_6$ are the normalized Eisenstein series,
$\Delta$ is the modular discriminant, and $j=E_4^3/\Delta$.
Its Fourier constant term is $A_{n,r}$; division by this constant
gives the solution with constant term one.
When $\ell\ge5$ is prime,
\[
 Y^a(Y-1728)^b\widetilde F_{\ell-1}(Y)
\]
has $\ell$-integral coefficients and reduces to the monic polynomial
whose roots are the supersingular $j$-invariants in
$\overline{\F}_\ell$ \cite[Theorem~1]{KZ}.
Indeed, $2n+r/6=\ell/6-1$ also gives
$A_{n,r}\equiv1\pmod\ell$. The same four Jacobi families occur
as pole polynomials in modular differential equations
\cite[Theorem~4.2]{BSS}.

The Legendre family has a related classical supersingular
interpretation. For an odd prime $p=2m+1$,
\[
 (1-z)^mP_m\!\left(\frac{1+z}{1-z}\right)
 =\sum_{j=0}^m\binom mj^2z^j.
\]
The reduction is a Hasse polynomial for
$y^2=x(x-1)(x-z)$: for $z\ne0,1$ it vanishes exactly at the
supersingular Legendre curves \cite[\S1]{BM}.

\subsection{Hyperelliptic Jacobians}

For the curve $y^2=f(x)$, the roots of $f$ are the finite branch
points of the double cover of the projective line. Divisor classes
supported on the branch points describe the $2$-torsion of its
Jacobian, linking the permutation group of the roots to the
arithmetic of the abelian variety. Zarhin's theorems show that,
in the degree range below, the full symmetric group forces the
geometric endomorphism ring to be $\Z$, gives the Hodge group
$\Sp$, and gives the $\ell$-adic Lie algebra
$\Q_\ell\operatorname{Id}\oplus\mathfrak{sp}$ \cite{Zarhin02}.
They also determine the Hodge and Tate classes on all powers.
The following corollary collects these consequences under one
Galois-group hypothesis.

\begin{corollary}\label{cor:jacobians}
Let $f\in\Q[x]$ be separable of degree $d\ge9$ with
$\Gal(f/\Q)=S_d$, and let $J_f$ be the Jacobian of the smooth
projective curve with affine equation $y^2=f(x)$. Then
\[
 \dim J_f=\lfloor(d-1)/2\rfloor,\qquad
 \End_{\overline\Q}(J_f)=\Z.
\]
For $V=H_1(J_f(\mathbb C),\Q)$ and the alternating form $\psi$
induced by the principal polarization,
\[
 \Hdg(J_f)=\Sp(V,\psi),\qquad \MT(J_f)=\GSp(V,\psi).
\]
For every number field $K\subset\overline\Q$ and every prime $\ell$,
put $V_\ell=T_\ell(J_f)\otimes_{\Z_\ell}\Q_\ell$, and let
$\psi_\ell$ be the alternating Weil pairing with values in
$\Q_\ell(1)$ induced by the principal polarization.
The image of $\Gal(\overline\Q/K)$ on $V_\ell$ has Lie algebra
\[
 \Q_\ell\operatorname{Id}\oplus\mathfrak{sp}(V_\ell,\psi_\ell).
\]
Every rational Hodge class on every power of $J_f$ is a rational
linear combination of products of divisor classes. For every
$h\ge1$ and $c\ge0$, the cycle-class map
\[
 \CH^c((J_f^h)_K)\otimes\Q_\ell\longrightarrow
 H^{2c}_{\mathrm{\acute et}}
 ((J_f^h)_{\overline\Q},\Q_\ell(c))^{\Gal(\overline\Q/K)}
\]
is surjective. The Hodge, Tate and Mumford--Tate conjectures
therefore hold for $J_f$ and all its powers.
\end{corollary}

\begin{proof}
These are consequences of Zarhin's theorems \cite{Zarhin02}: Theorem~2.4 gives the geometric
endomorphism ring, Theorem~2.6 gives the scalar plus symplectic
$\ell$-adic Lie algebra, Theorem~8.1 gives the Tate statement over
every finite extension and on every power, and Theorem~9.1 and
Remark~9.2 give the Hodge and Mumford--Tate conclusions.
Passing to a number field preserves the Lie algebra because the
Galois subgroup has finite index. For powers, the groups act
diagonally on direct sums. A Galois-invariant class represented
by algebraic cycles over a finite extension descends with rational
coefficients by averaging over a finite Galois extension.
\end{proof}

For the Kaneko--Zagier family, take
$f=P_d^{(\lambda/3,\mu/2)}$ and apply
Theorem~\ref{thm:kz-all}. For fixed unequal nonnegative integers
$a,b$, take $f=P_d^{(a,b)}$ with $d\ge\max\{9,N(a,b)\}$ and
apply Theorem~\ref{thm:main}. For a dominant factor $h_n$ in
Theorem~\ref{thm:jacobi}, use its actual degree $d=M_n$.

For the derivative families, fix $r$ and take
$f=Q_{m,r,\delta}$ with $m\ge9$ sufficiently large in
Corollary~\ref{cor:derivative-galois}. The substitution $t=x^2$
also gives the map
\[
 \{y^2=F_{m,r,\delta}(x)\}\longrightarrow
 \{y^2=Q_{m,r,\delta}(t)\},\qquad
 (x,y)\longmapsto(x^2,y).
\]
This map extends to a degree-two morphism of the smooth
projective models. It is the curve on the right whose defining
polynomial has group
$S_m$ and hence satisfies the preceding corollary. The polynomial
$F_{m,r,\delta}$ has group $W_m$; that group does not satisfy the
symmetric-group hypothesis.

\subsection{Further questions}

The fixed rational-parameter theorem leaves a precise
irreducibility problem.

\begin{question}
Suppose $\alpha,\beta>-1$ are rational and satisfy the denominator
conditions of Theorem~\ref{thm:jacobi}. Is
$P_n^{(\alpha,\beta)}$ irreducible over $\Q$ for all sufficiently
large $n$?
\end{question}

The theorem bounds the complementary degree by
$O_{\alpha,\beta,A}(n/(\log n)^A)$ for every fixed $A>0$.
This permits a bounded nonconstant factor. To exclude it, the
resultant comparison would need sufficient control of the
primitive content, the moment valuations, and the root distances
at the denominator primes. An explicit criterion in the
parameters for a positive comparison gap would identify further
families to which the method applies.

For integral parameters, the exceptional degrees form an
effective finite set. Classifying them calls for more information
about the central numerators and small factors. In the symmetric
case, all-degree irreducibility cannot hold for every fixed $r$:
\[
 C_2^{7/2}(x)=\tfrac72(3x-1)(3x+1).
\]
The cases $r=0,1$ have no exceptions. Which other values of $r$
have this property, and how far can the uniform range $m\ge8r$
be enlarged? The constant $8$ is a sufficient bound and is not
claimed to be optimal. Corollary~\ref{cor:derivative-galois}
settles eventual maximality for each fixed $r$; it leaves the
all-degree Galois groups of the first-derivative families and
uniform Galois estimates for growing $r$ as separate questions.

The prescribed-factor estimate suggests a related question:
how does the arithmetic of a minimal polynomial constrain the
occurrence of its roots in varying Jacobi families? The cone
$x^2+y^2-2z^2=0$ gives one concrete instance. The bound
$23r>4(\ell-1)$ is uniform in degree and order, but it leaves
infinitely many possible pairs and does not complete the
classification of harmonic modes on that cone.

Finally, the local results in Appendix~\ref{sec:local-irred}
come from a dilation of digit-sum Newton polygons. A general
description of when multiplication of the degree by a denominator
prime dilates the lower hull could yield further local
irreducibility sequences, including cases where the global
resultant comparison gives no contradiction.
\appendix

\section{Local irreducibility}
\label{sec:local-irred}

Rational irreducibility does not in general imply irreducibility
over a completion. The following result gives the latter
conclusion in two infinite sequences of degrees.

\begin{theorem}\label{thm:kz-towers}
For every integer $a\ge0$ and every $\lambda,\mu\in\{-1,1\}$,
$P_{3^a}^{(\lambda/3,\mu/2)}$ is irreducible over $\Q_3$.
The polynomial $P_{5\cdot3^a}^{(\lambda/3,\mu/2)}$ is irreducible
over $\Q_3$ precisely when $\mu=-1$ or $a\ge2$.
For $\mu=1$, the local factor degrees are $1,1,3$ when $a=0$,
and $2,13$ when $a=1$.
\end{theorem}

The first assertion follows from a general prime-power
criterion. The second follows from a dilation identity for
the lower convex hull of a digit-sum function.

\subsection{Irreducibility in prime-power degrees}
\label{sec:prime-powers}

A denominator in one parameter, together with an exact
valuation in the other, forces a single Newton edge.

\begin{theorem}\label{thm:prime-power-criterion}
Let $p$ be a prime, $a\ge1$, $n=p^a$, and let $\alpha,\beta\in\Q$.
Suppose
\[
 v_p(\alpha)=-e<0,\qquad v_p(\beta)\ge0,
 \qquad v_p((\beta+1)_n)=v_p(n!)=\frac{n-1}{p-1}.
\]
Then $P_n^{(\alpha,\beta)}$ is irreducible over $\Q_p$.
More precisely, the monic polynomial
\[
 \frac{n!P_n^{(\alpha,\beta)}(2X-1)}
      {(n+\alpha+\beta+1)_n}
\]
has a single Newton edge joining
\[
 \left(0,en+\frac{n-1}{p-1}\right)\quad\text{and}\quad(n,0),
\]
with every intermediate coefficient point strictly above the edge.
\end{theorem}
\begin{proof}
Every denominator factor in \eqref{eq:monic-endpoint} has
valuation $-e$. Put $h_j=v_p((\beta+1)_j)$ and $h_0=0$.
For $1\le i\le a$, exactly $n/p^i$ of
$\beta+1,\ldots,\beta+n$ belong to $p^i\Z_p$.
The hypothesis on $h_n$ therefore shows that none belongs to
$p^{a+1}\Z_p$. For $1\le j<n$, write $t=v_p(j)$. Counting
residue classes gives
\[
 h_j\le\sum_{i=1}^a\left\lceil\frac j{p^i}\right\rceil
       =v_p(j!)+a-t,
 \qquad v_p\binom nj=a-t.
\]
If $b_j$ is the coefficient of $X^j$, then
\[
 v_p(b_j)=e(n-j)+v_p\binom nj+h_n-h_j
       \ge e(n-j)+h_n-v_p(j!).
\]
Writing $H=en+(n-1)/(p-1)$, the difference from the proposed
edge is at least
\[
 v_p(b_j)-H\left(1-\frac jn\right)
       \ge\frac{s_p(j)-j/n}{p-1}>0.
\]
The endpoint valuations are $H$ and $0$. Since $H\equiv1\pmod p$,
one has $(H,n)=1$. Every root has valuation $H/n$, so a monic
factor over $\Q_p$ of degree $d$ would have constant-term
valuation $dH/n\in\Z$. Hence $n\mid d$, proving irreducibility.
\end{proof}

\begin{corollary}\label{cor:odd-prime-half}
Let $p$ be an odd prime, let $a\ge1$ be an integer, and let
$\alpha\in\Q$ satisfy $v_p(\alpha)<0$. Both
$P_{p^a}^{(\alpha,-1/2)}$ and $P_{p^a}^{(\alpha,1/2)}$ are
irreducible over $\Q_p$. In particular, for every $a\ge1$ all four
$P_{3^a}^{(\lambda/3,\mu/2)}$ are irreducible over $\Q_3$.
\end{corollary}
\begin{proof}
For $n=p^a$ and odd $p$,
\[
 (1/2)_n=\frac{(2n)!}{4^n n!},\qquad
 (3/2)_n=\frac{(2n+1)!}{4^n n!},
\]
and $v_p((2n)!)=v_p((2n+1)!)=2(n-1)/(p-1)$.
Apply Theorem~\ref{thm:prime-power-criterion}.
The degree-one case in Theorem~\ref{thm:kz-towers} is immediate.
\end{proof}

\begin{remark}\label{rem:cgl-eisenstein}
In the notation of \cite{CGL}, the monic endpoint polynomial is
$S_n^{(\lambda,\mu)}(X-1)$. Corollary~\ref{cor:odd-prime-half}
proves the irreducibility prediction in \cite[Conjecture~1(4)]{CGL}.
The precise Newton edge has height $(3n-1)/2$ for $n=3^a$.
Thus its height need not be one, as in the usual Eisenstein
criterion. For example,
\[
 S_3^{(1,1)}(X-1)
 =X^3-\frac{63}{41}X^2+\frac{27}{41}X-\frac{81}{1189}.
\]
The single-edge criterion with relatively prime height and length
is the exact local assertion used here.
\end{remark}

\subsection{A dilation law and the degrees \texorpdfstring{$5\cdot3^a$}{5 times powers of 3}}\label{sec:five-tower}

Fix $n\ge1$ and $\lambda,\mu\in\{-1,1\}$. A dilation identity
gives another sequence of irreducible degrees. All valuations in
this subsection are normalized by
$v_3(3)=1$, and $s_3(r)$ denotes the sum of the base-$3$ digits of a
nonnegative integer $r$.

Legendre's formula $v_3(k!)=(k-s_3(k))/2$ separates each
coefficient valuation into a linear term and a digit correction.
Only the correction affects the vertices of the Newton polygon.
Multiplying the degree by $3$ shifts the base-$3$ digits; the lemma
below shows exactly how this stretches the lower convex hull.

For $\delta\in\{0,1\}$ set
\[
 \Phi_{n,\delta}(j)
 =s_3(2j+\delta)+s_3(n-j)-s_3(2n+\delta)
 \qquad(0\le j\le n),
\]
and let $H_{n,\delta}$ be the lower convex hull of these points,
viewed as a piecewise linear function on $[0,n]$.
The monic polynomial
\[
 T_n^{\lambda,\mu}(X)
 =\frac{n!}{(n+\lambda/3+\mu/2+1)_n}
 P_n^{(\lambda/3,\mu/2)}(2X-1)
\]
has coefficients
\[
 b_j=(-1)^{n-j}\binom nj
       \prod_{r=j+1}^n\frac{3(2r+\mu)}{6(n+r)+2\lambda+3\mu}.
\]
Since each denominator in this product is a $3$-adic unit,
\begin{equation}\label{eq:five-tower-valuations}
 v_3(b_j)
 =n-j+v_3((2n+\delta)!)-v_3((2j+\delta)!)-v_3((n-j)!)
 =\frac32(n-j)+\frac12\Phi_{n,\delta}(j),
 \qquad \delta=\frac{\mu+1}{2}.
\end{equation}
Thus its Newton polygon is obtained from $H_{n,\delta}$ by a vertical
scaling and the addition of a linear function. In particular it does
not depend on $\lambda$.

\begin{lemma}\label{lem:digit-dilation}
For every integer $m\ge1$,
\[
 H_{3m,0}(x)=H_{m,0}(x/3)\quad(0\le x\le3m),
 \qquad
 H_{9m,1}(x)=H_{m,0}(x/9)\quad(0\le x\le9m).
\]
All points whose abscissae are not divisible by $3$, respectively by
$9$, lie strictly above the indicated hull.
\end{lemma}
\begin{proof}
First, all slopes of $H_{m,0}$ belong to $[-1,1]$.
Digit subadditivity gives
\[
 \Phi_{m,0}(j)-\Phi_{m,0}(0)
 =s_3(2j)+s_3(m-j)-s_3(m)\ge-s_3(j)\ge-j.
\]
Writing $t=m-j$, it also gives
\[
 \Phi_{m,0}(j)
 \ge s_3(t)-s_3(2t)\ge-s_3(t)\ge-t.
\]
These inequalities bound the first and last slopes of the convex
hull, and hence every slope, as asserted.

For $0\le k<m$, put $u=v_3(m-k)$ and $w=v_3(k+1)$.
The identity $s_3(r-1)=s_3(r)-1+2v_3(r)$ for $r\ge1$ yields
\[
 \begin{aligned}
 \Phi_{3m,0}(3k)&=\Phi_{m,0}(k),\\
 \Phi_{3m,0}(3k+1)&=\Phi_{m,0}(k)+3+2u,\\
 \Phi_{3m,0}(3k+2)&=\Phi_{m,0}(k+1)+1+2w.
 \end{aligned}
\]
The endpoint identity also holds at $k=m$.
The hull through the abscissae divisible by $3$ is exactly
$H_{m,0}(x/3)$, and its slopes lie in $[-1/3,1/3]$.
The second displayed point is at least $3$ above the preceding grid
point, while the third is at least $1$ above the following grid
point. Each is therefore strictly above this hull. This proves the
first identity.

For the second identity, the grid points satisfy
$\Phi_{9m,1}(9k)=\Phi_{m,0}(k)$ for $0\le k\le m$.
The remaining digit calculations are
\[
 \begin{array}{c|ccc}
 r&1&2&3\\ \hline
 \Phi_{9m,1}(9k+r)-\Phi_{m,0}(k)-2u&3&4&3
 \end{array}
 \quad
 \begin{array}{c|ccccc}
 r&4&5&6&7&8\\ \hline
 \Phi_{9m,1}(9k+r)-\Phi_{m,0}(k+1)-2w&1&2&1&2&3.
 \end{array}
\]
They follow by separating the last two base-$3$ digits and using the
same identity for $s_3(r-1)$. The dilated hull has slopes in
$[-1/9,1/9]$. For $1\le r\le3$, comparison with the preceding grid
point shows a strict excess of at least $3-r/9$.
For $4\le r\le8$, comparison with the following grid point shows a
strict excess of at least $1-(9-r)/9$. Thus every off-grid point is
strictly above the dilated hull, proving the second identity.
\end{proof}

\begin{proof}[Completion of the proof of Theorem~\ref{thm:kz-towers}]
Corollary~\ref{cor:odd-prime-half} treats the degrees $3^a$.
Let $n=5\cdot3^a$. We prove the remaining assertion and determine
the Newton edge.
The six values $\Phi_{5,0}(j)$ are
\[
 1,\ 2,\ 1,\ 2,\ 3,\ 0.
\]
Their lower hull is the segment from $(0,1)$ to $(5,0)$, with all
intermediate points strictly above it. Repeated application of the
first identity in Lemma~\ref{lem:digit-dilation} gives the same hull
stretched horizontally to length $n=5\cdot3^a$ when $\mu=-1$.
For $\mu=1$ and $a\ge2$, first use the second identity and then the
first. Formula~\eqref{eq:five-tower-valuations} gives a single edge
\[
 (0,(3n+1)/2)\text{--}(n,0),
\]
with every intermediate coefficient point strictly above it.
Since $n$ is odd,
\[
 \gcd\!\left(n,\frac{3n+1}{2}\right)=1.
\]
The constant term of a monic factor of degree $d$ over $\mathbb Q_3$
would therefore have valuation $d(3n+1)/(2n)$, an integer only when
$n\mid d$. This proves local irreducibility in the stated range.
For the two omitted cases $\mu=1$, the same coefficient formula
gives the following Newton polygons, independently of $\lambda$:
\[
\begin{array}{c|c|c}
 n&\text{vertices}&\text{factor degrees over }\Q_3\\ \hline
 5&(0,8),(1,6),(4,1),(5,0)&1,3,1\\
 15&(0,23),(13,3),(15,0)&13,2.
\end{array}
\]
Each side has relatively prime height and horizontal length.
The Newton-polygon factorization therefore has one irreducible
factor of each listed degree. Both polynomials are reducible,
which proves the converse and the asserted factor degrees.
\end{proof}
\section{The explicit cutoff and finite verification}\label{sec:explicit-cutoff}

The positive gap in Section~\ref{sec:irred-proof} already gives
eventual irreducibility.  Here we keep enough of the error terms
to obtain a numerical cutoff.  The useful quantity is
$S_n+\log g_n$: the real estimate gives a lower bound for it if
a factorization exists, while the joint content and moment
calculation gives an upper bound.  The estimates below make these
bounds incompatible for $n\ge250\,000$.  A separate finite
comparison completes the proof below that cutoff.

\subsection{An explicit cutoff}

Let $\alpha=\lambda/3$, $\beta=\mu/2$, with
$\lambda,\mu\in\{-1,1\}$, and retain the notation $H_n$, $g_n$,
$S_n$, $L_n$, and $J_n$ from the exact criterion. Write
$T_2=T_{2,n}$, $T_3=T_{3,n}$, and put
\[
 N_0=250\,000,\qquad r=2\lambda+3\mu,\qquad \delta=(\mu+1)/2,
 \qquad E=\tfrac12\log2+\tfrac13\log3,
 \qquad L_* =\log(4+3\sqrt3).
\]
\begin{proposition}\label{prop:explicit-kz-cutoff}
For every $n\ge250\,000$ and every $\lambda,\mu\in\{-1,1\}$,
$P_n^{(\lambda/3,\mu/2)}$ is irreducible over $\Q$.
\end{proposition}
We prove the proposition in three steps. We use the notation
\[
 \theta(x;q,a)=\sum_{\substack{p\le x\text{ prime}\\p\equiv a\pmod q}}\log p,
 \quad \theta(x)=\sum_{p\le x\text{ prime}}\log p,
 \quad \psi(x)=\sum_{p^h\le x}\log p.
\]

\subsubsection{The necessary archimedean inequality}
For every proper factorization with smaller degree $1\le k\le n/2$,
we claim
\begin{equation}\label{eq:cutoff-arch}
 S_n+\log g_n\ge n(L_*-E)-6\log(12n)
 \qquad(n\ge N_0).
\end{equation}
Put
\[
 b_n=v_2\binom{n+\alpha}{n}.
\]
The level-by-level binomial valuation formula gives
$0\le b_n\le\log_2(3n+1)$. The endpoint identities become
\[
 D_{2,n}=2n-s_2(n)+b_n,
 \qquad
 D_{3,n}=\frac{3n}{2}
       +\frac{\delta+s_3(n)-s_3(2n+\delta)}2.
\]
Also
\[
 \log\kappa_n=(2n-s_2(n))\log2
                     +(3n/2-s_3(n)/2)\log3.
\]
Using $s_p(m)\le(p-1)(1+\log_p m)$ gives
\[
 nL_n\ge nL_*-\log(6n+3)
\]
and
\[
 \log\kappa_n-n(T_2+T_3)
 \ge-\log(3n+1)-2\log(3n)-\tfrac12\log3.
\]
Furthermore,
\[
 J_n\le\tfrac12\log(2n)+\tfrac13\log(6n+3).
\]

The Wronskian constant satisfies $C_n\ge1/3$.
Indeed, $C_1=(\alpha+1)(\beta+1)\ge1/3$ and
\[
 \frac{C_{m+1}}{C_m}
 =1+\frac{\alpha\beta}{(m+1)(m+c+1)},\qquad c=\alpha+\beta.
\]
If $\alpha\beta>0$, this ratio exceeds one. Otherwise
$\alpha\beta=-1/6$, $C_1\ge2/3$, and, for $m\ge1$,
\[
 \frac{C_{m+1}}{C_m}\ge1-\frac1{6m(m+1)}.
\]
The product is at least $1-\sum_{m\ge1}1/(6m(m+1))=5/6$,
so in that case $C_n\ge5/9$. In particular $2\sqrt{C_n}\ge1$.

Let $a_n=T_2+T_3+E-2L_n+2\log2$.
The displayed endpoint bounds give, for $n\ge N_0$,
\[
 |T_2-T_3|
 \le\tfrac32\log3-2\log2+\frac{2\log(9n)}n
 <\frac{27}{100}+\frac{30}{250000}<\frac3{10}.
\]
Here we used $\log3<1.1$, $\log2>0.69$, and
$\log(9N_0)<15$, together with the decrease of $\log(9n)/n$.
Since $\log\cosh u\le u^2/2$, it follows that
\[
 a_n=E-2\log\cosh((T_2-T_3)/2)>\frac12.
\]
Consequently $a_nk+2\log2/k-\log k>0$ for $k>0$:
indeed $k/2-\log k\ge1-\log2>0$.
The exact Vandermonde criterion therefore implies
\begin{align*}
 S_n+\log g_n
 &\ge \log\kappa_n-J_n+n(L_n-T_2-T_3-E)-4\log2\\
 &\ge n(L_*-E)\\
 &\quad-\left\{\frac43\log(6n+3)+\log(3n+1)+2\log(3n)
                +\frac12\log(2n)+\frac12\log3+4\log2\right\}.
\end{align*}
The expression in braces is at most $6\log(12n)$ for $n\ge2$.
For example, the difference is positive at $n=2$ and its derivative
is greater than $7/(6n)$. This proves \eqref{eq:cutoff-arch}.

\subsubsection{Content and the good-prime sum}
Put $M=12n+5$ and
\[
 H(M)=\sum_{\substack{p\text{ prime},\ h\ge2\\p^h\le M}}\log p.
\]
For each prime $p\ge5$, choose a coefficient attaining the
minimum $m_p^{(1)}$ of its first-level valuation. Every higher-level
coefficient contribution is at most two. Thus
\[
 v_p(g_n)\le m_p^{(1)}
          +2\#\{h\ge2:p^h\le M\}.
\]
Since each exponent $e_{p^h}$ is at most one,
\begin{equation}\label{eq:cutoff-combined}
 S_n+\log g_n\le T_0(n)+3H(M),\qquad
 T_0(n)=\sum_{p\ge5}
       \bigl(m_p^{(1)}+\mathbf1_{p\le n}e_p\bigr)\log p.
\end{equation}
The two indicators $m_p^{(1)},e_p$ cannot both equal one,
so each summand in $T_0$ has multiplicity at most one.

For $p\ge13$, the first-level carry calculation gives the following
exact intervals. The first class is $p\equiv3+2\lambda\pmod6$,
and the second is $p\equiv3-2\lambda\pmod6$.
For $p\le n$, write $\ell=\lfloor n/p\rfloor\ge1$.
The supported intervals are
\[
 \begin{array}{c|c}
 \text{first class}&
 \displaystyle\left(\frac{6n+r}{6\ell+1},\frac n\ell\right]
 \quad\text{and}\quad
 \displaystyle\left(\frac n{\ell+1},\frac{2n+\mu}{2\ell+1}\right]
 \\[6pt]
 \text{second class}&
 \displaystyle\left(\frac n{\ell+1},\frac{3n+\lambda}{3\ell+1}\right].
 \end{array}
\]
For $p>n$, only content contributes, in the intervals
\[
 \left(n,\frac{3n+\lambda}{2}\right]
 \quad\text{in the first class},\qquad
 \left(\frac{6n+r}{5},2n+\mu\right]
 \quad\text{in the second class}.
\]
These formulas follow by substituting
$a_p=(p+\lambda)/3$ or $(2p+\lambda)/3$ and
$b_p=(p+\mu)/2$ into the indicators of
Lemma~\ref{lem:carry-budget}. In the first case the combined
indicator for $p\le n$ is $1+E_b-E_c$; in the second it is $E_a$.
The small primes $5,7,11$ contribute at most $\log385$.

We retain the first six groups of intervals.  They contain enough
of the main term for the required contradiction, while all discarded
primes are at most $n/7$ and contribute at most $\theta(n/7)$.
For the retained intervals, the square-root error in the prime
estimate is useful near the cutoff, and the relative error estimate
treats the unbounded range.
The coefficient of $n$ in one half the sum of all retained interval
lengths is
\begin{align*}
 Q_6={}&\frac{13}{20}
 +\frac12\sum_{\ell=1}^{6}\left\{
 \frac1\ell-\frac6{6\ell+1}
 +\frac2{2\ell+1}-\frac1{\ell+1}
 +\frac3{3\ell+1}-\frac1{\ell+1}\right\}\\
 ={}&\frac{180556377763}{157066509600}<\frac{23}{20}.
\end{align*}
The constant shifts in the interval lengths contribute at most
\[
 \frac54+\frac56\sum_{\ell=1}^{6}\frac1\ell<5.
\]
Together with the small primes, this is less than $12$.

We use two published estimates of Bennett, Martin, O'Bryant, and
Rechnitzer \cite[Corollary~1.7 and Theorem~1.9]{BMOR}:
for $a\in\{1,5\}$,
\begin{align*}
 |\theta(x;6,a)-x/2|&\le2.072\sqrt x
                  &&(1\le x\le4\cdot10^{13}),\\
 |\theta(x;6,a)-x/2|&<0.024\,x/\log x
                  &&(x\ge10^6).
\end{align*}
The range $4\cdot10^{13}$ in the first estimate follows from
$x_2(6)=x_2(3)$ in their Theorem~1.9.
For $N_0\le n\le10^{13}$, all retained endpoints are at most
$2n+1<4\cdot10^{13}$. For each $1\le\ell\le6$, the six
endpoints are at most $(n+1)/\ell$; the four remaining endpoints
are at most $2n+1$. Their total error is therefore at most
\[
 2.072\left\{6\sqrt{n+1}
          \sum_{\ell=1}^{6}\ell^{-1/2}+4\sqrt{2n+1}\right\}
 <61\sqrt n.
\]
To obtain the last inequality, use
$\sum_{\ell=1}^{6}\ell^{-1/2}\le2\sqrt6-1<3.9$,
$\sqrt{1+1/N_0}<1.000002$, and
$\sqrt{2+1/N_0}<1.415$.
For $n\ge10^{13}$, every retained endpoint is at least $n/7$.
Since $\sum_{\ell=1}^{6}\ell^{-1}=49/20$, their sum is less than
\[
 6(n+1)\frac{49}{20}+4(2n+1)<23n.
\]
Using $\log(10^{13}/7)>27$, the second prime estimate bounds
the total error by
\[
 0.024\frac{23n}{27}<\frac{21}{1000}n.
\]

Hanson's elementary theorem \cite[Theorem~1]{Hanson} states
$\operatorname{lcm}(1,\ldots,N)<3^N$. Consequently
$\theta(x)\le\psi(x)<x\log3$ for $x\ge1$.
Applying this to the discarded primes gives
\begin{equation}\label{eq:cutoff-level-one}
 T_0(n)<\left(\frac{23}{20}+\frac{\log3}{7}\right)n+12+
 \begin{cases}
 61\sqrt n,&N_0\le n\le10^{13},\\
 \dfrac{21}{1000}n,&n>10^{13}.
 \end{cases}
\end{equation}

Hanson's estimate also gives
\[
 H(M)\le\log3\left\{\sqrt M+(\log_2M)M^{1/3}\right\}
       <11\sqrt n\qquad(n\ge N_0).
\]
For an explicit verification, use $M\le12.1n$ and note that
$n^{-1/6}\log(12.1n)$ decreases for $n\ge N_0$.
At the left endpoint the rational bounds
\[
 \sqrt{12.1}<3.5,\quad 12.1^{1/3}<2.3,\quad
 N_0^{-1/6}<0.127,\quad \log(12.1N_0)<15,
 \quad \log2>0.69,\quad \log3<1.1
\]
give
\[
 \frac{H(M)}{\sqrt n}
 <1.1\{3.5+2.3\cdot0.127\cdot15/0.69\}
 =10.835<11.
\]
Combining this with \eqref{eq:cutoff-combined} and
\eqref{eq:cutoff-level-one}, and rounding the first square-root
coefficient upward, gives
\begin{equation}\label{eq:cutoff-upper}
 S_n+\log g_n<\left(\frac{23}{20}+\frac{\log3}{7}\right)n+12+
 \begin{cases}
 95\sqrt n,&N_0\le n\le10^{13},\\
 \dfrac{21}{1000}n+33\sqrt n,&n>10^{13}.
 \end{cases}
\end{equation}

\subsubsection{The contradiction}
We first record the lower bound $L_*-E>3/2$.
Indeed, $\sqrt3>1.73$ gives $4+3\sqrt3>9\cdot919/900$.
Since $\log3>1.098$, $\log2<0.7$, and
$\log(1+t)>t/(1+t)$ for $t>0$,
\[
 L_*-E>
 \frac53\log3-\frac12\log2+\log\frac{919}{900}
 >\frac53\frac{549}{500}-\frac7{20}+\frac{19}{919}
 =\frac{34478}{22975}>\frac32.
\]
For $N_0\le n\le10^{13}$, the first line of
\eqref{eq:cutoff-upper}, the decrease of $\log(12n)/n$, and
$\log(12N_0)<15$ give
\[
 \frac{23}{20}+\frac{\log3}{7}
 +\frac{95}{\sqrt n}+\frac{12+6\log(12n)}n
 <\frac{23}{20}+\frac{11}{70}+\frac{19}{100}
                  +\frac{102}{250000}
 =\frac{1310357}{875000}<\frac32.
\]
For $n>10^{13}$, the decreasing error satisfies
\[
 \frac{33}{\sqrt n}+\frac{12+6\log(12n)}n
 <\frac{33}{10^6}+\frac{210}{10^{13}}<\frac1{1000},
\]
where $\log(12\cdot10^{13})<33$. Hence the second line gives
\[
 \frac{23}{20}+\frac{\log3}{7}+\frac{21}{1000}
 +\frac{33}{\sqrt n}+\frac{12+6\log(12n)}n
 <\frac{23}{20}+\frac{11}{70}+\frac{22}{1000}
 <\frac43<\frac32.
\]
Both ranges contradict \eqref{eq:cutoff-arch}.
This proves irreducibility for every $n\ge250\,000$ in all four
families.  The numerical inequalities used here follow from the
displayed rational bounds and the positive-term series for the
logarithm.

\subsection{The finite range}\label{sec:finite-verification}

We use a published computation for the small degrees.
Cullinan and Hajir \cite[Remark~4.6]{CHK} checked the
irreducibility of $K_\ell$ for every $13\le\ell\le10^5$
relatively prime to $6$.  Their Definition~2.3 gives
\[
 \ell=12n+6+2\lambda+3\mu,
 \qquad
 K_\ell(X)=3^nn!P_n^{(\lambda/3,\mu/2)}(4X+1).
\]
Consequently their result covers all four parameter choices for
$1\le n\le8332$, since $12\cdot8332+11=99\,995$.
We use their irreducibility computation here, independently of
the conditional Galois-group statement in \cite[Theorem~1.8]{CHK}.

\begin{proposition}\label{prop:finite-kz}
For every $1\le n\le250\,000$ and every
$\lambda,\mu\in\{-1,1\}$, the polynomial
$P_n^{(\lambda/3,\mu/2)}$ is irreducible over $\Q$.
\end{proposition}

\begin{proof}
It remains to treat $8333\le n\le250\,000$.
Retain $T_0(n)$ from \eqref{eq:cutoff-combined} and put
\[
 Q(n)=\sum_{\substack{p\ge5\text{ prime},\ h\ge2\\p^h\le12n+5}}
              \log p.
\]
The same content and moment calculation gives
\begin{equation}\label{eq:fast-content-budget}
 S_n+\log g_n\le T_0(n)+3Q(n).
\end{equation}
Write $T_2=T_{2,n}$, $T_3=T_{3,n}$, $L=L_n$, and $J=J_n$.
Since $C_n\ge1/3$, the exact real comparison implies, for every
possible factor degree $1\le k\le n/2$,
\begin{align*}
 \frac{\log R^{(6)}}k\ge{}&
 \log\kappa_n-\log g_n-J+n(L-T_2-T_3-E)-4\log2\\
 &+ak+\frac{2\log2}{k}-\log k,
 \qquad a=T_2+T_3+E-2L+2\log2.
\end{align*}
If $|T_2-T_3|\le9/10$, then
\[
 a=E-2\log\cosh((T_2-T_3)/2)
 \ge E-\frac{81}{400}>\frac12.
\]
Thus the last three terms in the second line have positive sum,
because $k/2-\log k\ge1-\log2>0$.
Weighted AM--GM with weights $7/16$ and $9/16$ gives
\[
 L\ge\frac7{16}T_2+\frac9{16}T_3
           +4\log2-\frac7{16}\log7-\frac98\log3.
\]
Define the integers
\begin{align*}
 c_2&=216n+48v_2(n!)-24s_2(n)-27D_{2,n}-192,\\
 c_3&=48v_3(n!)-22n-21D_{3,n}-8s_3(2n+\delta),\\
 c_7&=-21n,
 \qquad \delta=(\mu+1)/2,
\end{align*}
and the real number
\[
 \Delta_n=
 \frac{c_2\log2+c_3\log3+c_7\log7}{48}-T_0(n)-3Q(n).
\]
Combining the preceding bounds shows that the two conditions
\begin{equation}\label{eq:fast-finite-criterion}
 |D_{2,n}\log2-D_{3,n}\log3|\le\frac9{10}n,
 \qquad \Delta_n>0
\end{equation}
exclude every proper factorization.

For each of the four parameter pairs and every integer
$8333\le n\le250\,000$, exact integer arithmetic verifies the
first inequality in \eqref{eq:fast-finite-criterion} and the
stronger bound
\[
 \Delta_n>\frac{n}{1000}.
\]
The calculation uses only the residue indicators of
Lemma~\ref{lem:carry-budget}, digit sums, and rational bounds for
logarithms.  For completeness, these bounds follow by truncating
\[
 \log\frac{1+t}{1-t}
 =2\sum_{j=0}^{r-1}\frac{t^{2j+1}}{2j+1}+\mathcal E_r(t),
 \qquad
 0\le\mathcal E_r(t)\le
 \frac{2t^{2r+1}}{(2r+1)(1-t^2)}
 \quad(0\le t\le1/3).
\]
Rescaling an integer by a power of $2$ puts its logarithm in this
form.  The accompanying self-contained verification recomputes
these comparisons from the formulas, with outward rounding and
without stored numerical data. An independent calculation in the repository
updates $T_0(n)$ only at primes dividing
$n(3n+\lambda)(2n+\mu)(6n+2\lambda+3\mu)$ and bounds the logarithms
by repeated squaring. Over the entire range it gives
$\Delta_n/n\ge0.136734$ and
$|D_{2,n}\log2-D_{3,n}\log3|/n\le0.262378$.
This completes the proof.
\end{proof}

We also checked the range $n\le8332$ independently of
\cite[Remark~4.6]{CHK}. For $1810\le n\le8332$ the two inequalities in
\eqref{eq:fast-finite-criterion} hold. For $n\le1809$, write
$P_n^{(\lambda/3,\mu/2)}(1-2t)/P_n^{(\lambda/3,\mu/2)}(1)
=\sum c_jt^j$ and put $s=2\lambda+3\mu$. Then
\[
 \frac{c_{j+1}}{c_j}
 =\frac{(j-n)(6(n+j+1)+s)}{2(j+1)(3(j+1)+\lambda)}.
\]
If $p=6(n+h)+s$ is prime with $1\le h\le n$, the first $h$
coefficients are $p$-adic units and the rest have valuation one.
The last edge of the Newton polygon has height one and length
$n-h+1$, so it belongs to a single factor over $\Q_p$.
When $2(h-1)<n$, this bounds the smaller rational factor degree by $h-1$.
Starting with this bound when available and $\lfloor n/2\rfloor$
otherwise, squarefree reductions of full degree exclude every remaining
candidate by Lemma~\ref{lem:finite-factor-degrees}. The code is included
with the Lean files.
\section{The finite first-derivative calculation}
\label{app:finite-derivatives}

We give the remaining finite verification for the first derivatives,
and record elementary rational bounds for three logarithmic constants.
The coefficient formulas and prime bounds below specify the
calculation completely.

\subsection{A factor-degree criterion}

\begin{lemma}\label{lem:finite-factor-degrees}
Let $H\in\mathbb Z[x]$ be primitive of degree $D\ge1$.
For a prime $p$ at which its reduction has degree $D$ and is
squarefree, let $\mathcal D_p$ be the multiset of irreducible
factor degrees, and let $\Sigma_p$ be the set of its subset sums,
with each entry of the multiset used at most once.
If
\[
 \{1,\ldots,\lfloor D/2\rfloor\}
 \cap\bigcap_{p\in\mathcal P}\Sigma_p=\varnothing
\]
for a finite set $\mathcal P$ of such primes, then $H$ is
irreducible over $\mathbb Q$.
\end{lemma}

\begin{proof}
By Gauss's lemma, a nontrivial rational factorization gives an
integral factor of degree $1\le k\le\lfloor D/2\rfloor$.
Since the leading coefficient of $H$ is nonzero modulo every
listed prime, the factor retains its degree. Its reduction is a
product of distinct irreducible factors, so $k\in\Sigma_p$
for every $p\in\mathcal P$, a contradiction.
\end{proof}

\subsection{The finite ranges}

For $n=2m+\delta$, where $\delta\in\{0,1\}$, let $H_n$ be
the primitive integral polynomial with positive leading coefficient
proportional to
\begin{equation}\label{eq:finite-derivative-coefficients}
 \sum_{k=0}^{m}
 \frac{(-1)^k(2n+2-2k)!}
 {k!(n+1-k)!(n-2k)!}\,x^{n-2k-\delta}.
\end{equation}
Thus $H_n$ is a nonzero rational multiple of
$P'_{n+1}(x)/x^\delta$ and has degree $2m$.
Put
\[
 U=\prod_{p\le m}p^{v_p(n+1)},\qquad
 V=\prod_{p\le m}p^{v_p(n+2)},\qquad
 Z=\prod_{p\le m}p^{\lfloor\log_p(2n+3)\rfloor}.
\]
For $m\ge4$, the sufficient condition of
Theorem~\ref{thm:symmetric-exact} is exactly
\begin{equation}\label{eq:finite-derivative-integer}
 16^mV^2>64m(m+1)^2Z^2U^2.
\end{equation}
Indeed, $K=m+1$ and $A_{\le m}=U/V$; exponentiating and
squaring the criterion gives \eqref{eq:finite-derivative-integer}.

\begin{lemma}\label{lem:first-derivative-finite}
For $1\le m\le199$ and $\delta\in\{0,1\}$,
$H_{2m+\delta}$ is irreducible over $\mathbb Q$.
\end{lemma}

\begin{proof}
For each parity, the indices $m<4$ together with the failures of
\eqref{eq:finite-derivative-integer} for $4\le m\le199$ are
\[
\begin{aligned}
 E_0={}&\{1,\ldots,25\}\cup
 \{28,30,31,32,34,37,40,42,43,45,46,47\},\\
 E_1={}&\{1,\ldots,26\}\cup
 \{29,31,32,33,34,35,37,38,41,43,44,45,47,49,50,53,73\}.
\end{aligned}
\]
The other $318$ pairs satisfy the strict integer inequality.
For the $37+43=80$ pairs with $m\in E_\delta$, apply
Lemma~\ref{lem:finite-factor-degrees} to the full-degree,
squarefree reductions at primes $p\le317$. In every case the
intersection in that lemma is empty.

Here is a complete procedure for this finite calculation.
The coefficients are obtained from
\eqref{eq:finite-derivative-coefficients} by dividing by their
positive gcd. The integers $U,V,Z$ use trial division and
successive multiplication by each prime. For each remaining
$H_n$, start with the candidate set $\{1,\ldots,m\}$ and visit
the primes up to $317$ in increasing order. Discard a prime
if the degree drops or if $\gcd(\overline H_n,\overline H_n')\ne1$.
At every other prime, compute $\mathcal D_p$ as follows.

Make $f=\overline H_n/\operatorname{lc}(\overline H_n)$ monic,
and initialize $R=f$ and $h=x$. For $j=1,2,\ldots$, as long
as $2j\le\deg R$, replace $h$ by $h^p\bmod R$, put
$G=\gcd(R,h-x)$, record $j$ exactly $\deg(G)/j$ times, and
replace $R$ by $R/G$ and $h$ by its remainder modulo this new
$R$. If a nonconstant $R$ remains at the end, record its degree
once. This gives $\mathcal D_p$: the polynomial $x^{p^j}-x$
is the product of the monic irreducibles whose degrees divide
$j$, and all smaller degrees have already been removed.
The final polynomial $R$ has at most one irreducible factor:
each remaining factor has degree at least $j$, while $\deg R<2j$.

To compute $\Sigma_p$, start with $\{0\}$ and, for every
entry $d$ of $\mathcal D_p$, adjoin its translate by $d$.
Intersect the candidate set with $\Sigma_p$, stopping when it
is empty. Retaining only reductions that strictly decrease
that set gives $254$ reductions in total, with largest prime
$317$. These operations use integer arithmetic and Euclidean
division over finite fields. Lemma~\ref{lem:finite-factor-degrees}
proves the remaining $80$ cases.

The Lean proof uses Eisenstein's criterion for $45$ of the primitive
polynomials and for the reciprocals of another $24$. For each of the other
$11$, it verifies an irreducible reduction of degree $d=2m$ by checking
$X^{p^d}\equiv X\pmod{\overline H_n}$ and
$\gcd(\overline H_n,X^{p^{d/q}}-X)=1$ for every prime $q\mid d$.
The largest prime $p$ used in these eleven certificates is $1163$.
A separate list gives one prime of irreducible reduction for each of the
$80$ polynomials. Its largest entry is $2311$.
\end{proof}

There is also a short check using Newton polygons. For an edge with vector
$(u,v)$, put $d=u/\gcd(u,|v|)$. The horizontal length contributed by this
edge to a rational factor is one of $0,d,\ldots,u$, by Dumas's product rule
\cite[pp.~216--217]{Dumas}. Intersecting the
possible sums of these lengths at primes up to $151$ excludes every factor
of degree $1,\ldots,m$ in $77$ cases. For $(\delta,m)=(1,32)$ one also
uses the polygon of $H_n(X+1)$ at $2$. In the two remaining cases,
$(\delta,m)=(1,25),(1,43)$, only degree two survives. Write
$H_n(X)=G(X^2)$. A monic irreducible quadratic factor and its reflection
under $X\mapsto-X$ would give a quartic factor unless they coincide.
The quadratic must therefore be $X^2-a$, with $G(a)=0$ and $a\in\Q$.
For $m=25$ and $43$, respectively, $G$ has leading coefficient prime to
$13$ and $193$ and has no root in the corresponding prime field. This excludes
a rational root in both cases.

\subsection{Rational logarithm bounds}

For $0<z<1$, define
\[
 T(z)=2\sum_{j=0}^{4}\frac{z^{2j+1}}{2j+1},\qquad
 U_*(z)=T(z)+\frac{2z^{11}}{11(1-z^2)}.
\]
The power series for the inverse hyperbolic tangent gives
\[
 T(z)<\log\frac{1+z}{1-z}<U_*(z).
\]
The two margins used in the paper have the exact expressions
\[
\begin{aligned}
 8(\log4-1)-g(8)
 &=13\log2+9\log(9/8)-34\log(17/16)-8,\\
 8\left[\frac{19}{8}(\log6-1)-g(19/8)\right]
 &=3\log2+8\log(3/2)+57\log(19/16)\\
 &\qquad+92\log(24/23)-19.
\end{aligned}
\]
Substitution of the rational bounds gives
\[
 13T(1/3)+9T(1/17)-34U_*(1/33)-8>\frac{97}{10000}
\]
and
\[
 3T(1/3)+8T(1/5)+57T(3/35)+92T(1/47)-19
 >\frac{42}{1250}.
\]
Both are direct comparisons of rational numbers; the second
right side is $8\cdot21/5000$. Finally,
\[
 2(\log9-1)-g(2)
 =\log\left(\frac{3^7 2^{16}}{5^{10}}\right)-2>0,
\]
because $e<3$ and $3^5 2^{16}>5^{10}$.

\end{document}